\documentclass[reqno, 11pt]{amsart}

\makeatletter
\let\origsection=\section \def\section{\@ifstar{\origsection*}{\mysection}} 
\def\mysection{\@startsection{section}{1}\z@{.7\linespacing\@plus\linespacing}{.5\linespacing}{\normalfont\scshape\centering\S}}
\makeatother  

\usepackage{mathrsfs}
\usepackage{amsmath}
\usepackage{amssymb}
\usepackage{amsthm}
\usepackage{mathtools}
\usepackage{letterswitharrows}
\usepackage{tikz-cd}
\usetikzlibrary{intersections,spath3}
\usepackage{thmtools}
\usepackage{thm-restate}
\usepackage[only,llbracket,rrbracket]{stmaryrd}
\usepackage[linesnumbered,ruled,vlined]{algorithm2e}
\usepackage{enumitem}
\setenumerate{label={\normalfont (\roman*)}}

\usepackage[utf8]{inputenc}
\usepackage[T1]{fontenc}
\usepackage{lmodern}
\usepackage[babel]{microtype}
\usepackage[english]{babel}
\usepackage{relsize}

\usepackage{graphicx}
\usepackage{subcaption}
\usepackage{import}

\usepackage{svg}
\usepackage[extractname=filename]{svg-extract}

\usepackage{geometry}
\usepackage{xcolor} 
\colorlet{darkishRed}{red!80!black}
\colorlet{darkishBlue}{blue!60!black}
\colorlet{darkishGreen}{green!60!black}
\definecolor{lucasBlue}{HTML}{A4D4F5}
\definecolor{lucasDarkGreen}{HTML}{8CBF65}
\definecolor{lucasLightRed}{HTML}{F5B3B9}
\definecolor{lucasOrange}{HTML}{EBC07A}

\definecolor{okorange}{RGB}{230,159,0}

\definecolor{okskyblue}{RGB}{86,180,233}

\definecolor{okgreen}{RGB}{0,158,115}

\definecolor{okyellow}{RGB}{240,228,66}

\definecolor{okblue}{RGB}{0,114,178}

\definecolor{okvermillion}{RGB}{213,94,0}

\definecolor{okpurple}{RGB}{204,121,167}

\colorlet{oktintorange}{okorange!45}

\colorlet{oktintskyblue}{okskyblue!45}

\colorlet{oktintgreen}{okgreen!45}

\colorlet{oktintyellow}{okyellow!45}

\colorlet{oktintblue}{okblue!45}

\colorlet{oktintvermillion}{okvermillion!45}

\colorlet{oktintpurple}{okpurple!45}

\colorlet{oktintblack}{black!15}

\usepackage{hyperref}
\hypersetup{
	colorlinks,
	linkcolor={blue!60!black},
	citecolor={green!60!black},
	urlcolor={blue!60!black},
    hypertexnames={false},
}
\usepackage[abbrev, msc-links]{amsrefs}
\usepackage[nameinlink, capitalise, noabbrev]{cleveref}
\crefformat{enumi}{#2#1#3}
\crefformat{equation}{#2(#1)#3}
\crefname{mainresult}{Theorem}{Theorems}
\usepackage{doi}

\renewcommand{\PrintDOI}[1]{\doi{#1}}

\let\setminus=\smallsetminus

\renewcommand{\leq}{\leqslant}
\renewcommand{\geq}{\geqslant}

\newtheorem{theorem}{Theorem}[section] 

\newtheorem{corollary}[theorem]{Corollary}
\newtheorem{lemma}[theorem]{Lemma}

\newtheorem{observation}[theorem]{Observation}

\newcounter{claimcounter}[theorem]
\newtheorem{claim}[claimcounter]{Claim}
\newtheorem*{claim*}{Claim}
\crefname{claim}{Claim}{Claims}
\newcounter{casecounter}[theorem]
\newtheorem{case}[casecounter]{Case}
\newtheorem*{case*}{Case}
\crefname{case}{Case}{Cases}

 \newcounter{subcasecounter}[casecounter]
 \newtheorem{subcase}[subcasecounter]{Subcase}

\newcounter{subclaimcounter}[claimcounter]
\newtheorem*{subclaim*}{Subclaim}

\theoremstyle{definition}

\crefname{mainexample}{Example}{Examples}
\newtheorem{example}[theorem]{Example}
\crefname{example}{Example}{Examples}

\newtheorem{remark}[theorem]{Remark}

\crefname{routine}{Routine}{Routines}

\crefname{subroutine}{Subroutine}{Subroutines}

\crefname{subsubroutine}{Subsubroutine}{Subsubroutines}

\crefname{step}{Step}{Steps}
\theoremstyle{remark}

\crefname{subsection}{Subsection}{Subsections}

\newcommand{\COMMENT}[1]{{}}

\let\epsilon=\varepsilon
\let\theta=\vartheta
\let\rho=\varrho
\let\phi=\varphi

\makeatletter

\def\calCommandfactory#1{%
  \expandafter\def\csname c#1\endcsname{\mathcal{#1}}}
\def\frakCommandfactory#1{%
  \expandafter\def\csname frak#1\endcsname{\mathfrak{#1}}}
   
\newcounter{ctr}
\loop
  \stepcounter{ctr}
  \edef\X{\@Alph\c@ctr}
  \expandafter\calCommandfactory\X
  \expandafter\frakCommandfactory\X
\ifnum\thectr<26
\repeat

\usepackage{etoolbox}

\newbool{arXiv}
\booltrue{arXiv} 

\newcommand{\arXivOrNot}[2]{\ifbool{arXiv}{{#1}}{{#2}}}

\newbool{pdfBool}
\booltrue{pdfBool} 

\newcommand{\pdfOrNot}[2]{\ifbool{pdfBool}{{#1}}{{#2}}}

\makeatletter
\newcommand\thankssymb[1]{\textsuperscript{\@fnsymbol{#1}}}
\makeatother

\newcounter{mylabelcounter}

\makeatletter
\newcommand{\labelText}[2]{%
#1\refstepcounter{mylabelcounter}%
\immediate\write\@auxout{%
  \string\newlabel{#2}{{1}{\thepage}{{\unexpanded{#1}}}{mylabelcounter.\number\value{mylabelcounter}}{}}%
}%
}

\makeatletter
\newcommand\footnoteref[1]{\protected@xdef\@thefnmark{\ref{#1}}\@footnotemark}
\makeatother

\usetikzlibrary{fit,shapes.geometric,calc,decorations.pathmorphing, math, backgrounds}

\usepackage{stmaryrd} 
\usepackage{tabularx}
\usepackage{varwidth}

\definecolor{cMaroon}{HTML}{93152a}
\newcommand{\defn}[1]{{\color{cMaroon}{\emph{#1}}}}

\newcommand{\cof}[1]{\operatorname{cf}({#1})}

\newcommand{\TwStar}{T_{\aleph_0}\star1}
\newcommand{\TkStar}{T_{\kappa}\star1}
\newcommand{\TkM}{T_{\kappa}^{\#}} 

\newcommand{\Above}{\mathord{\uparrow}}

\newcommand{\successor}[1]{\mathrm{succ}({#1})}
\newcommand{\upcl}[1]{\lfloor {#1} \rfloor}
\newcommand{\downcl}[1]{\lceil {#1} \rceil}

\usepackage{accents,graphicx}
\newcommand{\oacc}[4]{\text{$\lineskiplimit=0pt #1{\raisebox{#2}[\height][\depth]{\scalebox{.65}{$\circ$}}}{\mathstrut\smash[#3]{#4}}$}}

\newcommand{\odowncl}[1]{\lceil \oacc\accentset{-.2ex}{t}{#1} \rceil}
\newcommand{\oupcl}[1]{\lfloor \oacc\underaccent{.2ex}{b}{#1} \rfloor}

\newcommand{\Set}[1]{{\left\lbrace {#1} \right\rbrace}}
\def\set#1:#2{\Set{{#1} \colon {#2}}}

\newcommand{\cardinality}[1]{{\left\lvert {#1} \right\rvert}}

\title{The unavoidable minor in graphs of uncountable edge-connectivity}

\author[M.~Pitz]{Max Pitz}

\author[L.~Real]{Lucas Real\thankssymb{1}}
\thanks{\thankssymb{1} Supported by S{\~a}o Paulo Research Foundation (FAPESP) through grant number 2025/00669-5.}

\author[R.~Schaut]{Roman Schaut\thankssymb{2}}
\thanks{\thankssymb{2} Supported by a State Graduate Funding Program scholarship of the University of Hamburg.}

\address{University of Hamburg, Department of Mathematics, Bundesstrasse 55 (Geomatikum), 20146 Hamburg, Germany}
\email{\{max.pitz, roman.schaut, lucas.real\}@uni-hamburg.de}

\begin{document}

\begin{abstract}
    For any given uncountable cardinal $\kappa$, we show that every $\kappa$--edge-connected graph admits a minor embedding of the graph $\TkStar$ (the $\kappa$--regular tree with a new apex vertex adjacent to all its nodes), which is itself $\kappa$--edge-connected.
    This extends a previous result of Kurkofka regarding unavoidable minors of $\aleph_0$--edge-connected graphs to all uncountable cardinalities $\kappa$.
\end{abstract}



\vspace*{-12pt}
\maketitle

\section{Introduction}
Which substructures must every graph of a given connectivity contain?
This question goes back to Halin and Jung~\cite{HalinJung1963}, who asked, for every $k\in\mathbb{N}$, for the minimal elements of the class of $k$--connected graphs with respect to a given containment relation.
As one might suspect, the length of the resulting list of unavoidable graphs depends on how weak this relation is, and the minor relation typically yields the shortest list.
For finite graphs, Halin and Jung determined the minor-minimal $k$--connected graphs for all $k \leq 4$~\cite{HalinJung1963}*{\S2}, and Oporowski, Oxley and Thomas~\cite{Oporowski} found the unavoidable minors of all sufficiently large $3$--connected and internally $4$--connected graphs; see~\cite{bibel}*{Section~9.4} for an overview of the finite theory. 
For infinite graphs, vertex- and edge-connectivity are treated in Chapters~5 and~6 of the forthcoming book~\cite{pitzbook}, which also surveys the results below.

We follow the terminology of~\cite{bibel}.
In particular, a \defn{minor embedding} of a graph $H$ in a graph $G$ is understood as a disjoint family $\set{B_v}:{v\in V(H)}$ of connected sets of vertices of $G$ such that there is a $B_u$--$B_v$ edge in $G$ whenever $uv\in E(H)$.
Each set $B_v$ is called a \defn{branch set} or \defn{bag} corresponding to $v\in V(H)$, while the subgraph of $G$ induced by $\bigcup_{v\in V(H)}B_v$ is said to be an \defn{inflation} or \defn{minor model} of $H$.
The existence of such an embedding is expressed by saying that $H$ is a \defn{minor} of $G$.

For infinite graphs, the answer depends on the size of the graph.
Halin~\cite{halin78}*{Theorem~9.1} showed that every $k$--connected graph of uncountable regular size $\kappa$ contains a subdivision of $K_{k,\kappa}$, building on work of Dirac~\cite{Dirac1977} for $k=2$ and of Halin and Jung~\cite{HalinJung1963}*{Satz~6} for rayless graphs.
For countable graphs, Oporowski, Oxley and Thomas~\cite{Oporowski}*{(5.2)} showed that every essentially $k$--connected infinite graph contains as a minor either $K_{k,\aleph_0}$ or one of finitely many graphs built from rays arranged along a finite tree, and Gollin and Heuer~\cite{enddefiningsequences} extended this picture to $k$--connected sets of arbitrary infinite size.
Finally, for a (possibly infinite) cardinal $\kappa$, we call two vertices $u,v \in V(G)$ \defn{$\kappa$--vertex-connected} if they are joined by $\kappa$ many internally disjoint paths, and a graph is \defn{$\kappa$--connected} if its vertices are pairwise $\kappa$--vertex-connected.
If $\kappa$ is infinite, the answer becomes simple: every $\kappa$--connected graph contains a subdivision of $K_\kappa$~\cite{halin1967}.

For edge-connectivity, the picture changes completely.
Here, we recall that $G$ is said to be \defn{$\kappa$--edge-connected} if $\cardinality{V(G)} \geq 2$ and $G - F$ is connected for every set of edges $F\subseteq E(G)$ of size $<\kappa$; vertices belonging to different components of $G - F$ are said to be \defn{separated} by this set of edges.
High edge-connectivity does not force high vertex-connectivity: the graph $\TkStar$, obtained from the $\kappa$--regular tree by adding an apex vertex adjacent to all its nodes, is $\kappa$--edge-connected (see \cref{sec:examples}), yet deleting the apex together with any node of the tree disconnects it.
So the results above say nothing about edge-connectivity.
For $\kappa = \aleph_0$, Kurkofka~\cite{kurkofka2022} showed that every $\aleph_0$--edge-connected graph contains either $\TwStar$ or the \emph{Farey graph} as a minor, and Knappe and Kurkofka~\cite{knappe2024immersion} showed that the halved Farey graph is unavoidable for strong immersions.
Finding the unavoidable minors for \textit{uncountable} edge-connectivity was the core of Problem~7.3 in~\cite{kurkofka2022}.
Our work addresses this problem in full: for uncountable $\kappa$ the Farey-type obstruction disappears, and $\TkStar$ is the \textit{unique} minor-minimal $\kappa$--edge-connected graph:

\begin{theorem}\label{thm:mainThm}
    For any given uncountable cardinal $\kappa$, every $\kappa$--edge-connected graph contains $\TkStar$ as a minor.
\end{theorem}

As in the countable case, where Kurkofka~\cite{kurkofka2022}*{Theorem~3.3} observed that the two minor-minimal graphs are not the only unavoidable ones for the topological minor relation, our result is best possible in this sense: \cref{example:LackOfSpreadingPair} provides a $\kappa$--edge-connected graph which does not contain a minor embedding of $\TkStar$ whose branch sets have size $<\kappa$, and therefore no subdivision of $\TkStar$ either.

\subsection{Proof sketch and overview}\label{sec:OrganizingPaper}

For a fixed uncountable cardinal $\kappa$, a minor embedding of $\TkStar$ in $G$ could be easily obtained if the following criterion is met: there is a (connected) set $A \subseteq V(G)$ of size $< \kappa$ and some component of $G - A$ in which one can find a minor embedding of $T_{\kappa}$ in which every branch set meets $N(A)$. 
Then, \cref{thm:mainThm} follows once we declare $A$ as the branch set for an apex vertex for this tree minor. 
This approach motivates a brief overview in \cref{sec:RankFunctions} of the hierarchical (rank) description of graphs without embeddings of $T_{\kappa}$.

However, the $\kappa$--edge-connected graph $\TkM$ introduced in \cref{example:interleaved} (and depicted in \cref{fig:TkM}) does not satisfy this criterion, although it can be easily seen to contain a minor embedding of $\TkStar$ whose branch set for the corresponding apex vertex is of size $\kappa$.
The core of \cref{sec:NormalTrees} is then to extract further structural information about such graphs $G$ to which the above methodology does not apply for finding a model of $\TkStar$.
In particular, we will conclude that if $\kappa$ is regular, then the graph $\TkM$ is the only reason why the above methodology can fail.
This is now derived in a two-step approach.

After revisiting the theory of normal partition trees, our first step is to establish the existence of a suitable $T$--graph minor of $G$ which reflects several of its connectivity properties.
In particular, as detailed throughout \cref{subsection:NPTconstruction}, the proof of our main \cref{thm:mainThm} is reduced to a suitable family of $\kappa$--edge-connected $T$--graphs.

In a second step, and as the main engine of the paper, \cref{thm:sufficientCondition} within \cref{sec:sufficientCondition} then formalizes a sufficient condition on $T$--graphs for constructing a minor of $\TkM$, and thus embedding $\TkStar$ as a byproduct.  
Combined with this sufficient condition, the reduction of \cref{sec:NormalTrees} will then establish our main \cref{thm:mainThm} immediately when $\kappa$ is a regular cardinal. 
The same strategy almost works for general values of $\kappa$ as well, and a uniform proof could be obtained provided that $G$ does not contain a minor embedding of the graph $\mathbb{T}_{\kappa}$ from \cref{example:TreeOfCliques} (and depicted in \cref{fig:TreeOfCliques}).
Then, the final \cref{sec:singCard} isolates this specific $\cof{\kappa}$--edge-connected obstruction (which admits $\TkStar$ as a minor) and shows how its construction can be combined with the previous main approach.

\subsection*{Acknowledgements and use of AI}
The authors wrote this paper and used the AI assistant Claude (Anthropic) only afterwards to proofread the manuscript and to provide language and grammar corrections.

\section{Examples of highly edge-connected graphs}\label{sec:examples}

If $\kappa$ is a given infinite cardinal, this paper fixes the notations $K_{\kappa}$ and $T_{\kappa}$ for the clique of cardinality $\kappa$ and the $\kappa$--regular (rooted) tree, respectively.  
The aim of this section is to provide examples of $\kappa$--edge-connected graphs that, as mentioned in the above proof sketch, will also play a role in the later sections.
As part of the statement of our main \cref{thm:mainThm}, we first verify that $\TkStar$ is indeed $\kappa$--edge-connected.

\begin{example}[The $\kappa$--regular tree with an apex vertex]
    We recall that \defn{$\TkStar$} is obtained from $T_{\kappa}$ after adding an apex vertex $\mathrm{p}$ adjacent to all its nodes, as illustrated by \cref{figure:Tk1}.
    If $t\in T_{\kappa}$ is some node and if $\mathrm{succ}(t)$ denotes the set of all successors of $t$ in $T_\kappa$, define $P_s \coloneq tsp$ for $s \in \mathrm{succ}(t)$.
    Then the $\kappa$ internally disjoint paths $\set{P_s}:{s \in \mathrm{succ}(t)}$ verify that $t$ and $\mathrm{p}$ cannot be separated by fewer than $\kappa$ edges in $\TkStar$.
    In particular, no two vertices of $\TkStar$ can be separated by fewer than $\kappa$ edges, witnessing the $\kappa$--edge-connectivity of this graph.
\end{example}

As motivated in the proof sketch from the introduction, under certain assumptions we are able to construct a minor of the following $\kappa$--edge-connected graph:

\begin{example}[Interleaved copies of $T_{\kappa}$]\label{example:interleaved}
    The graph \defn{$\TkM$} is defined as the supergraph of $T_\kappa$ on the same vertex set $V(\TkM)=V(T_\kappa)$ where every node $t \in T_\kappa$ is adjacent to every neighbor of its successors in $T_\kappa$. 
    More precisely, $\TkM$ is defined on the edge set $$E(T_{\kappa}) \cup \set{ts}:{t < s \in T_\kappa \text{ and }\mathrm{dist}_T(t,s) = 2}.$$
    In order to verify that $\TkM$ is $\kappa$--edge-connected, denote by $r$ the root of $T_\kappa$ and let $t \neq r \in T_{\kappa}$ be an arbitrary node.
    It suffices to prove that no fewer than $\kappa$ many edges can separate $t$ from the root $r$.
    Let $s$ be the unique predecessor of $t$ in $T_{\kappa}$ and define $P_y \coloneq syt$ for every $y \in \successor{t}$.
    Then the set $\set{P_y}:{y \in \successor{t}}$ provides $\kappa$ internally vertex disjoint $s$--$t$ paths in $\TkM$, so that these two nodes cannot be separated in $\TkM$ by a set of fewer than $\kappa$ many edges.
    By induction on the height of $t$ in $T_\kappa$, no fewer than $\kappa$ many edges can separate $s$ from the root $r$ and thus $t$ and $r$ are $\kappa$--edge-connected.
\end{example}

\begin{figure}
    \centering
    \begin{minipage}{0.45\textwidth}
    \centering
    \begin{tikzpicture}[
        scale=0.25,
        TreeNode/.style = {circle, fill=black, inner sep=1.5pt},
        ApexVertex/.style = {circle, fill=black, inner sep=2pt},
        TreeEdge/.style = {black, line width=0.6pt},
        ApexEdge/.style = {gray, line width=0.4pt}
    ]
        \node[left] at (32,31) {$\mathrm{p}$};
        \node[TreeNode] (A) at (22,12) {}; 
        \node[TreeNode] (B) at (12,18) {};
        \node[TreeNode] (C) at (20,22) {};
        \node[TreeNode] (D) at (22,22) {};
        \node[TreeNode] (E) at (10,22) {};
        \node[TreeNode] (F) at (12,22) {};
        \node[TreeNode] (G) at (14,22) {};
        \node[TreeNode] (H) at (24,22) {};
        \node[TreeNode] (I) at (22,18) {};
        \node[TreeNode] (J) at (30,22) {};
        \node[TreeNode] (K) at (32,22) {};
        \node[TreeNode] (L) at (34,22) {};
        \node[TreeNode] (M) at (32,18) {};

        \node[ApexVertex] (X) at (30,30) {};

        \draw[TreeEdge,thick] (A) to (B);
        \draw[TreeEdge,thick] (A) to (I);
        \draw[TreeEdge,thick] (A) to (M);

        \draw[TreeEdge,thick] (E) to (B);
        \draw[TreeEdge,thick] (F) to (B);
        \draw[TreeEdge,thick] (G) to (B);

         \draw[TreeEdge,thick] (C) to (I);
        \draw[TreeEdge,thick] (D) to (I);
        \draw[TreeEdge,thick] (H) to (I);

        \draw[TreeEdge,thick] (J) to (M);
        \draw[TreeEdge,thick] (K) to (M);
        \draw[TreeEdge,thick] (L) to (M);

        \draw[ApexEdge,dashed] (X) to[in=50, out=270] (A);
        \draw[ApexEdge,dashed] (X) to[in=0, out=180] (B);
        \draw[ApexEdge,dashed] (X) to[in=0, out=270] (C);
        \draw[ApexEdge,dashed] (X) to[in=0, out=270] (D);
        \draw[ApexEdge,dashed] (X) to[in=0, out=180] (E);
        \draw[ApexEdge,dashed] (X) to[in=0, out=180] (F);
        \draw[ApexEdge,dashed] (X) to[in=0, out=180] (G);
        \draw[ApexEdge,dashed] (X) to[in=0, out=270] (H);
        \draw[ApexEdge,dashed] (X) to[in=0, out=270] (I);
        \draw[ApexEdge,dashed] (X) to[in=0, out=0] (J);
        \draw[ApexEdge,dashed] (X) to[in=0, out=0] (K);
        \draw[ApexEdge,dashed] (X) to[in=0, out=0] (L);
        \draw[ApexEdge,dashed] (X) to[in=0, out=0] (M);

    \end{tikzpicture}
    \caption{The graph $\TkStar$ is obtained from a $\kappa$--regular tree after adding an apex vertex $\mathrm{p}$.}
    \label{figure:Tk1}
    \end{minipage}
    \hfill
    \begin{minipage}{0.45\textwidth}
    \centering
    \begin{tikzpicture}[
        scale=0.1,
        vertex/.style = {circle, fill=black, inner sep=1.5pt},
        TreeEdge/.style = {black, line width=0.6pt},
        ApexEdge/.style = {gray, line width=0.4pt},
        ApexEdge2/.style = {gray, line width=0.6pt}
    ]

    \foreach \name/\x/\y in {
        S/45/14,
        Q/25/30,
        R/65/30,
        P/14/44,
        O/34/44,
        N/54/44,
        M/74/44,
        A/10/50,
        B/14/50,
        C/18/50,
        D/30/50,
        E/34/50,
        F/38/50,
        G/50/50,
        H/54/50,
        I/58/50,
        J/70/50,
        K/74/50,
        L/78/50
    }
     \node[vertex] (\name) at (\x,\y) {};   

     \foreach \x/\y in {
        Q/A,
        Q/B,
        Q/C,
        R/G,
        R/H,
        R/I}
        \draw[ApexEdge,dashed] (\x) to[in=-100,out=150] (\y);

         \foreach \x/\y in {
        Q/D,
        Q/E,
        Q/F,
        R/J,
        R/K,
        R/L}
        \draw[ApexEdge,dashed] (\x) to[in=-80,out=30] (\y);
    
    \foreach \x/\y in {
        A/P,
        B/P,
        C/P,
        D/O,
        E/O,
        F/O,
        G/N,
        H/N,
        I/N,
        J/M,
        K/M,
        L/M,
        P/Q,
        O/Q,
        N/R,
        M/R,
        Q/S,
        R/S}
    \draw[TreeEdge, thick] (\x) to (\y);

    \draw[ApexEdge2,dashed] (S) to[in=-20,out=110] (P);
    \draw[ApexEdge2,dashed] (S) to[in=-20,out=110] (O);

    \draw[ApexEdge2,dashed] (S) to[in=200,out=70] (N);
    \draw[ApexEdge2,dashed] (S) to[in=200,out=70] (M);
     
    \end{tikzpicture}
    \caption{The graph $\TkM$ is obtained from the rooted $\kappa$--regular tree $T_{\kappa}$ after adding edges between nodes and successors of their successors.}
    \label{fig:TkM}
    \end{minipage}
\end{figure}

We now verify that the graph $\TkM$ defined in the above \cref{example:interleaved} indeed contains $\TkStar$ as a minor.
This can be observed as follows: if $P\subseteq V(\TkM)$ denotes the set of vertices of even height in the fixed tree $T_{\kappa}$ as in \cref{example:interleaved}, then this set induces a connected subgraph of $\TkM$ (by definition of $\TkM$) and trivially has all other vertices of $T_{\kappa}$ (namely, those of odd height) in its neighborhood. 
Then let $Q\coloneq V(\TkM) \setminus P$ be the set of all vertices of odd height in $T_{\kappa}$.
The induced subgraph $\TkM[Q]$ then has $\kappa$ many components, each isomorphic to $T_{\kappa}$ and given by $\set{y \in \TkM}:{y \geq_T s} \setminus P$ for every $s$ of height $1$ in $T_{\kappa}$.
In particular, any component of $\TkM[Q]$ together with the connected set $P$ describes a minor embedding of $\TkStar$ in $\TkM$.

\begin{figure}[ht]
    \centering
    \begin{tikzpicture}[scale=0.5, transform shape]

  \fill[rounded corners=1.1pt, black!18] (0.27,4.03) -- (2.72,1.5) -- (2.73,1.76) -- (0.27,5.5) -- cycle;
 \fill[rounded corners=1.1pt, black!18] (2.72,1.28) -- (9.12,-4) -- (9.12,-4) -- (2.75,1.79) -- cycle;
 \fill[rounded corners=1.1pt, black!18] (13.22,1.5) -- (9.12,-3.8) -- (9.12,-3.8) -- (13.15,3) -- cycle;
\fill[rounded corners=1.1pt, black!18] (9.15,1.04) -- (9.02,-3.8) -- (9,-3.8) -- (8.78,1.54) -- cycle;
 \fill[rounded corners=1.1pt, black!18] (5,4.76) -- (2.77,1.5) -- (2.75,1.74) -- (4.45,6.28) -- cycle;
\fill[rounded corners=1.1pt, black!18] (-0.54,7.37) -- (0.26,4.46) -- (0.26,5.02) -- (-0.27,8.62) -- cycle;
\fill[rounded corners=1.1pt, black!18] (7.74,4.18) -- (9,1.5) -- (9,2) -- (8.28,5.72) -- cycle;
\fill[rounded corners=1.1pt, black!18] (9.46,6.26) -- (9,2.25) -- (9,1.5) -- (10,4.74) -- cycle;
\fill[rounded corners=1.1pt, black!18] (1,7) -- (0.26,4.46) -- (0.27,4.74) -- (0.8,9.55) -- cycle;
\fill[rounded corners=1.1pt, black!18] (3,4.46) -- (2.76,1.26) -- (2.75,1.52) -- (2.5,5.48) -- cycle;
\fill[rounded corners=1.1pt, black!18] (2.93,8) -- (2.77,4.74) -- (2.72,5.28) -- (2.43,9.02) -- cycle;
\fill[rounded corners=1.1pt, black!18] (13.4,5.03) -- (13.24,1.77) -- (13.19,2.31) -- (13,6.01) -- cycle;
  \fill[black, rounded corners=3.3pt] (2.5,2) rectangle (3,1);
  \fill[black, rounded corners=3.3pt] (8.75,-3.5) rectangle (9.29,-4.09);
  \fill[black, rounded corners=2.1pt] (0,5.48) rectangle (0.51,4);
  \fill[black, rounded corners=3.3pt] (2.5,6) rectangle (3.04,4);
  \fill[black, rounded corners=3pt] (4.45,6.55) rectangle (5,4);
  \fill[black, rounded corners=3.3pt] (8.74,2.48) rectangle (9.25,1);
  \fill[black, rounded corners=3.3pt] (7.74,6.01) rectangle (8.28,4.01);
  \fill[black, rounded corners=3.3pt] (9.46,6.55) rectangle (10,4);
  \fill[black, rounded corners=3.1pt] (-0.54,8.99) rectangle (0,6.99);
  \fill[black, rounded corners=3.3pt] (0.51,9.56) rectangle (1.08,7);
  \fill[black, rounded corners=3.3pt] (2.43,9.55) rectangle (3,6.99);
  \fill[black, rounded corners=3.3pt] (12.95,3.02) rectangle (13.49,1.02);
  \fill[black, rounded corners=3.3pt] (12.95,6.55) rectangle (13.48,4);
  \node[draw=none, node font=\LARGE] at (9.61,-4.33) {$\emptyset$};
  \node[draw=none, node font=\LARGE] at (2.35,0.65) {$p$};
\node[draw=none, node font=\LARGE] at (5.13,3.61) {$q$};
\end{tikzpicture}
    
    \caption{Illustration of the graph $\mathbb{T}_{\kappa}$ for a given singular cardinal $\kappa$ in terms of the tree $(\cof{\kappa}^{<\omega},\sqsubseteq)$. 
    In this diagram, rectangles of different heights correspond to cliques of different sizes.}
    \label{fig:TreeOfCliques}
\end{figure}
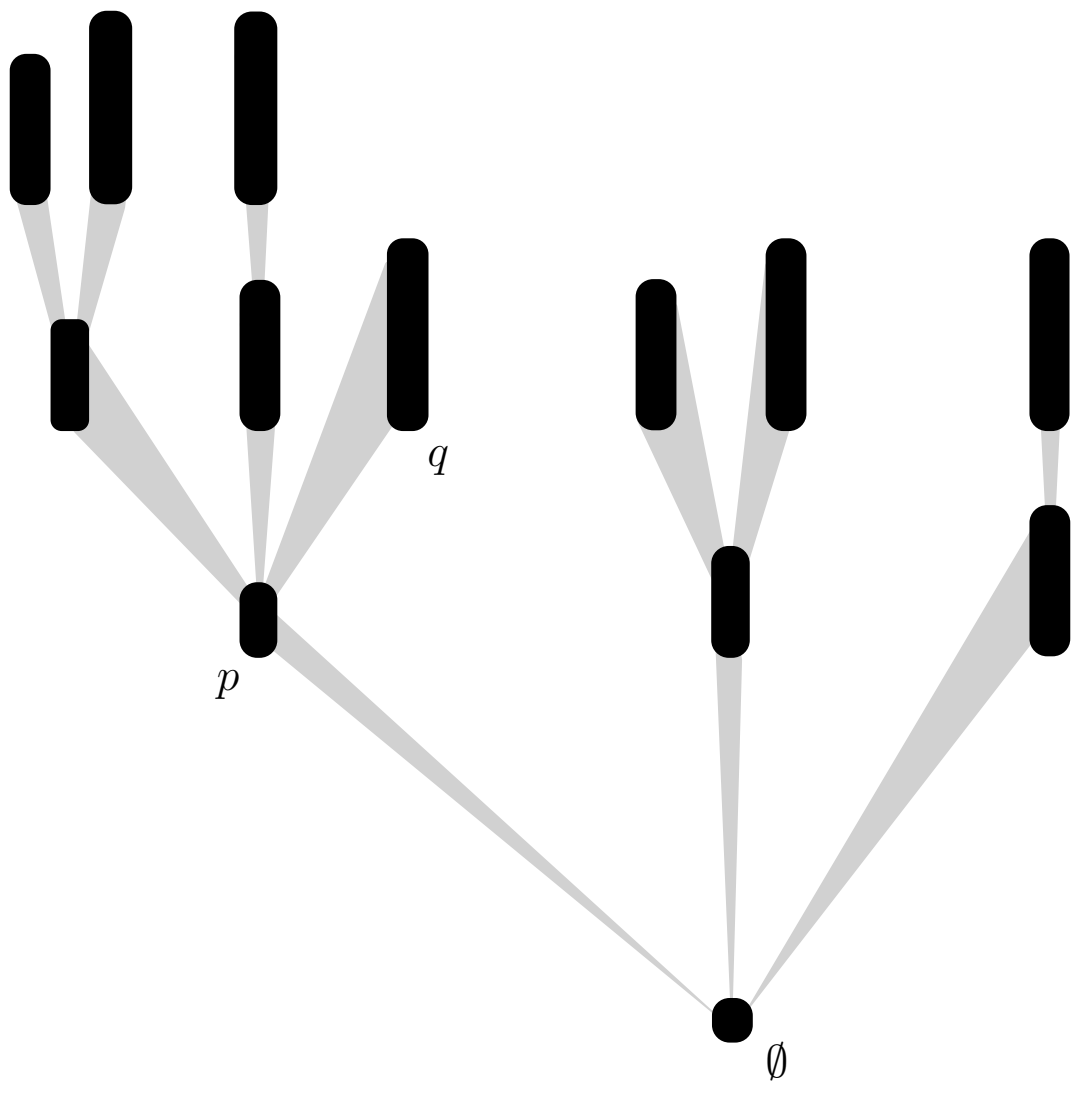

Lastly, the next example is particular to the setting in which $\kappa$ is a singular cardinal, but provides a priori only a $\cof{\kappa}$--edge-connected graph.
Despite that, it still contains $\TkStar$ as a minor and will appear as an intermediate step throughout \cref{sec:singCard} for reducing the main \cref{thm:mainThm} to a better-tailored statement:

\begin{example}[Tree of cofinal cliques]\label{example:TreeOfCliques}
    Suppose that $\kappa$ is a singular cardinal and let $\Set{\kappa_{\alpha}}_{\alpha < \cof{\kappa}} \subseteq \kappa$ be a strictly increasing sequence of cardinals, all larger than $\cof{\kappa}$, cofinal in $\kappa$.
    Then, we recall that the set $\cof{\kappa}^{<\omega}$ of finite (strictly increasing) sequences in $\cof{\kappa}$ has a natural ordering $\sqsubseteq$, where $p\sqsubseteq q$ if, and only if, $p$ is an initial segment of $q$.
    In particular, $(\cof{\kappa}^{<\omega},\sqsubseteq)$ describes a graph-theoretical tree if one fixes the empty sequence as its root, and where successors are defined by the extension via one additional element.

    The \defn{tree of cofinal cliques} is then the graph \defn{$\mathbb{T}_{\kappa}$} obtained by inflating $p \in \cof{\kappa}^{<\omega}$ by a clique $K_p$ of size $\kappa_{\ell(p)}$, where $\ell(p)$ is the last element of $p$ or $1$ if $p$ is the empty sequence, and then adding all edges of the form $\set{uv}:{u\in K_t,v \in K_s}$ whenever $st \in E(\cof{\kappa}^{<\omega})$.
    See also \cref{fig:TreeOfCliques} for an illustration.
\end{example}

The graph $\mathbb{T}_{\kappa}$ contains $\TkStar$ as a minor: fix a set $P \subseteq V(\mathbb{T}_{\kappa})$ obtained by picking an arbitrary vertex from the clique $K_p$ for every $p\in \cof{\kappa}^{<\omega}$. 
In particular, the definition of $E(\mathbb{T}_{\kappa})$ ensures that $P$ is connected and every vertex from $Q \coloneq V(\mathbb{T}_{\kappa}) \setminus P$ is adjacent to some vertex of $P$.

Furthermore, since every vertex in $Q$ is a vertex of some clique $K_p$, the subgraph $\mathbb{T}_{\kappa}[Q]$ is connected and $\kappa$--regular: if a given vertex $v \in Q$ is an element of the clique $K_p$, then $v$ is adjacent to every vertex from $\set{K_q}:{q\in \mathrm{succ}(p)}$, while $\sup\set{\kappa_{\ell(q)}}:{q \in \mathrm{succ}(p)} = \kappa$ by definition of $(\cof{\kappa}^{<\omega},\sqsubseteq)$.
In summary, $\mathbb{T}_{\kappa}[Q]$ is a $\kappa$--regular graph and thus contains $T_\kappa$ as a subgraph.
Combining this subgraph with the set $P$ as an apex vertex bag describes a minor embedding of $\TkStar$ in $\mathbb{T}_{\kappa}$.

\section{Rank functions and tree minors}\label{sec:RankFunctions}

For a fixed uncountable cardinal $\kappa$, recall the \defn{$\kappa$--rank} as introduced by Bürger and Kurkofka in \cite{BuergerKurkofka2022}:
assign to every graph whose vertex set is of size $<\kappa$ the $\kappa$--rank $0$.
Given an ordinal $\alpha > 0$ and a graph $G$ which has no $\kappa$--rank $<\alpha$, assign $\alpha$ as the $\kappa$--rank of $G$ if there exists a set $X \subseteq V(G)$ of size $<\kappa$ for which every component of $G - X$ has $\kappa$--rank $<\alpha$.
Together with a result of Seymour and Thomas in \cite{seymour1993binarytree}, Bürger and Kurkofka deduced in \cite{BuergerKurkofka2022} that the graphs with a $\kappa$--rank are exactly those graphs that do not admit $T_\kappa$ as a minor; for a textbook account, see also \cite{pitzbook}*{Section~9.5}.

\begin{theorem}\cite{BuergerKurkofka2022}*{Theorem~3.2}\label{thm:kappaRank}
    Let $\kappa$ be an uncountable cardinal.
    Let $G$ be a graph.
    The following are equivalent:
    \begin{itemize}

        \item $G$ has a $\kappa$--rank;
        
        \item $G$ contains no $T_\kappa$ minor.
        
    \end{itemize}
\end{theorem}

For our purposes, however, we will also need a localized version of the $\kappa$--rank, which assigns a rank to some subgraphs of a fixed global graph $G$.
The aim is then to characterize when a fixed graph $G$ contains a minor embedding of $T_{\kappa}$ whose branch sets meet a prescribed set of vertices $A\subseteq V(G)$.

As in the previous hierarchical definition, fix an uncountable cardinal $\kappa$ and assign \defn{$(A,\kappa)$--rank} $0$ to every subgraph $H \subseteq G$ with $\cardinality{V(H)\cap A} < \kappa$.
Given an ordinal $\alpha >0$ and a subgraph $H \subseteq G$ without an $(A,\kappa)$--rank $<\alpha$, assign the $(A,\kappa)$--rank $\alpha$ to $H$ if there exists a set $X \subseteq V(H)$ of size $<\kappa$ for which every component of $H-X$ already has $(A,\kappa)$--rank $<\alpha$.

When $\kappa$ is clear from the context, we will from now on write \defn{$A$--rank} instead of $(A,\kappa)$--rank.
By adapting the proof of \cref{thm:kappaRank} in \cite{pitzbook}*{Section~9.5}, the following result is the core of this section:

\begin{theorem}\label{lemma:Arank}
    Let $\kappa$ be an uncountable cardinal.
    Let $G$ be a fixed graph and $A \subseteq V(G)$ a set of vertices.
    Let $H \subseteq G$ be a subgraph.
    The following are equivalent:

    \begin{itemize}
            \item $H$ has no $(A,\kappa)$--rank;

            \item $H$ admits a minor embedding of $T_{\kappa}$ in which every bag intersects $A$.
    \end{itemize}
\end{theorem}

The proof of the above theorem relies on the following routine property regarding rank functions:

\begin{lemma}\label{lemma:RankReducing}
    If $G$ has an $A$--rank $\alpha$, then every subgraph of $G$ has $A$--rank at most $\alpha$.
\end{lemma}

\begin{proof}
    Let $H$ be a subgraph of $G$.
    We prove the statement by induction on the $A$--rank of $G$. 
    Observe that the statement is trivially true when $G$ has $A$--rank $0$: in this case, $G$ (and thus any of its subgraphs) contains fewer than $\kappa$ vertices from $A$.

    Suppose now that $G$ has $A$--rank $\alpha > 0$ and let $X \subseteq V(G)$ be a set of size $<\kappa$ such that the components of $G-X$ have $A$--rank $<\alpha$. 
    Since $H$ is a subgraph of $G$, the set $\set{H \cap C}:{C \text{ is a component of }G-X}$ defines a partition of $H-X$ into subgraphs which, by the induction hypothesis, have $A$--rank smaller than $\alpha$.
    Therefore, $H$ has $A$--rank at most $\alpha$.
\end{proof}

\begin{proof}[Proof of~\cref{lemma:Arank}]
    We can assume without loss of generality that $H = G$ and that $G$ is connected.
    Let us first assume that $G$ admits an $A$--rank, and let us assume further, towards a contradiction, that $G$ admits a minor embedding $\mathcal{W}$ of $T_\kappa$ in which every bag intersects $A$.
    Starting with $C_0 \coloneq G$ and inductively picking a rank reducing set $X_i$ in $C_i$ and a component $C_{i+1}$ of $C_{i} - X_i$ that still contains a subembedding of $\mathcal{W}$ isomorphic to $T_\kappa$, we obtain an infinite decreasing sequence $C_0 \supsetneq C_{1} \supsetneq \dots$ of graphs of strictly decreasing $A$--rank, a contradiction.

    For the other direction, let us assume now that $G$ admits no $A$--rank.
    For every cardinal $\mu\leq \kappa$, we say that a vertex $v\in V(G)$ is \textit{$\mu$--dominating} if, for every $X \subseteq V(G)$ of size $\cardinality{X} < \mu$, there exists at least one component of $G-X$ without an $A$--rank and in which $v$ has a neighbor. 
    We then claim the following:

    \begin{claim}\label{lem:Arank:claim}
        Fix a set of vertices $X' \subseteq V(G)$ with $\cardinality{X'} < \kappa$ and let $C$ be a component of $G - X'$ without an $A$--rank.
        Then, for every uncountable regular cardinal $\mu \leq \kappa$, there exists a $\mu$--dominating vertex in $C$.
    \end{claim}

    \begin{proof}[Proof of~\cref{lem:Arank:claim}]\renewcommand{\qedsymbol}{$\blacksquare$}
        Assume for a contradiction that, for some uncountable regular cardinal $\mu$ with $\cardinality{X'} < \mu \leq \kappa$, the component $C$ as above admits no $\mu$--dominating vertex. 
        Hence, for every $v\in C$, there is a set $X_v\subseteq V(G)$ with $\cardinality{X_v} < \mu$ and such that $v$ has no neighbors in any component of $G - X_v$ with no $A$--rank.
        We will now recursively define a sequence $X_0 \subseteq X_1 \subseteq X_2 \subseteq \dots$ of sets of vertices of $G$, all of size $< \mu$, as follows:

        \begin{itemize}
            \item $X_0\coloneq X'\cup \Set{v_*}$, where $v_*\in C$ is arbitrary. 
            Note that $\cardinality{X_0} \leq \cardinality{X'} + 1 < \mu$ by assumption.
            
            \item $X_{n+1}\coloneq X_n\cup \bigcup_{v\in D_n}X_v$, where the vertices from $D_n\coloneq C\cap X_n$ are not $\mu$--dominating by hypothesis. 
            Since $\cardinality{X_n} < \mu$ by induction, it follows from the regularity of $\mu$ that $\cardinality{X_{n+1}} < \mu$ as well.
        \end{itemize}

        Let $Y \coloneq \bigcup_{n\in\mathbb{N}}X_n$; the regularity of $\mu$ ensures that $\cardinality{Y} < \mu$. 
        In particular, since $C$ has no $A$--rank and $X' \subseteq Y$, there exists a component $D$ of $G-Y$ contained in $C$ and without an $A$--rank.
        On the other hand, the connectedness of $C$ implies the existence of a path $P \subseteq C$ starting at $v_*$ and having an endpoint in $D$. 
        Now because $P \cap X' = \emptyset$ and $v_* \in Y$, there is some $v \in Y \cap P \subseteq Y \setminus X'$ which has a neighbor in $D$.
        But then $v \in D_n$ for some index $n \in \mathbb{N}$ such that $v \in X_n$, contradicting the fact that $X_v \subseteq X_{n+1}$. 
        Indeed, by the choice of $X_v$, the vertex $v$ should have no neighbors in any component of $G - X_v$ with no $A$--rank and, in particular, in the one containing $D$.
    \end{proof}

    We shall now recursively construct a minor embedding $\set{V_i}:{i < \kappa}$ of $T_{\kappa}$ such that every bag intersects $A$. 
    To that end, fix an enumeration $\Set{t_i}_{i<\kappa}$ of $T_{\kappa}$ with the property that, for every $i>0$, the node $t_i$ is a successor of some node $t_j$ with $j<i$. 
    We moreover assume $T_{\kappa}$ to be rooted in $t_0$.

    Let $\theta \coloneq \aleph_0$ if $\kappa$ is regular and $\theta \coloneq \cof{\kappa}$ otherwise.
    Since $G$ admits no $A$--rank, it contains at least one vertex $v_0\in A$. 
    Moreover, applying \cref{lem:Arank:claim} to $X' = \emptyset$, it follows that for each uncountable regular cardinal $\mu \leq \kappa$ there exists a vertex in $V(G)$ that is $\mu$--dominating. 
    Then let $K_0 \subseteq V(G)$ be a set of size $\leq \theta$ which contains $v_0$ and, for every uncountable regular cardinal $\mu \leq \kappa$, a vertex that is $\mu$--dominating. 
    Next, let $V_0$ be a connected set of vertices in $G$ of size at most $\theta$ that contains $K_0$.

    Now let $i<\kappa$ and suppose that $\set{V_j}:{j<i}$ has already been constructed and corresponds to a minor embedding of the subtree of $T_{\kappa}$ induced by $\Set{t_j}_{j<i}$. 
    Since $\cardinality{V_j} \leq \theta$ for every $j<i$ by the inductive assumption, the set $X_i \coloneq \bigcup_{j<i} V_j$ has size $<\kappa$. 
    For the index $j < i$ for which $t_j$ is the predecessor of $t_i$ in $T_{\kappa}$ the set $V_j$ contains, by construction, a vertex $u_i \in V_j$ that is $\mu_i$--dominating for some uncountable regular cardinal $\mu_i > \cardinality{X_i}$.
    Then, by the definition of $\mu_i$--dominating vertices, there exists some component $C$ of $G - X_i$ with no $A$--rank and some $w_i \in C$ neighboring $u_i$; and since $C$ has no $A$--rank, fix also $v_i \in A \cap V(C)$.
    
    Applying \cref{lem:Arank:claim} to $X' = X_i$, it follows that there exists a set $K_i \subseteq V(C)$ of cardinality at most $\theta$ which contains at least one $\mu$--dominating vertex for each uncountable regular $\mu \leq \kappa$. 
    Let $V_i \subseteq V(C)$ be any connected set of size at most $\theta$ that contains $K_i \cup \Set{w_i,v_i}$. 
    Then the set $\set{V_j}:{j \leq i}$ is indeed a minor embedding of the subtree of $T_{\kappa}$ induced by $\Set{t_j}_{j \leq i}$: after all, the edge $u_iw_i$ witnesses an edge between $V_i$ and its predecessor bag $V_j$.
    At the end of this recursive process, $\set{V_i}:{i<\kappa}$ provides a minor embedding of $T_{\kappa}$ and every bag $V_i$ contains the vertex $v_i \in A$.
\end{proof}

\section{A reduction for \texorpdfstring{$T$--}{T–}graphs}\label{sec:NormalTrees}

As outlined in the proof sketch, our approach towards \cref{thm:mainThm} will be to first reduce its statement from the class of all $\kappa$--edge-connected graphs to a proper subclass of graphs whose structure can be roughly represented by an (order) tree.
After all, since our goal is to find a minor of $\TkStar$ (which by nature has a tree-like structure), this first reduction will come in handy in the upcoming sections.

 Our main tool to that aim is the \textit{normal partition tree} definition as introduced by Brochet and Diestel in \cite{BrochetDiestel}. 
 In our context, the goal is to construct such a partition tree over a given $\kappa$--edge-connected graph in order to obtain a $\kappa$--edge-connected contraction minor, which shall further satisfy some suitable additional properties.
 However, this methodology first requires some terminology on \textit{order trees}, namely, partially ordered sets of the form $(T,\leq)$ that have a unique $\leq$--minimum element (called its \defn{root}) and such that \defn{$\odowncl{t}$}${} \coloneq \set{s\in T}:{s<t}$ is well-ordered for every $t\in T$. For general background on $T$--graphs and normal partition trees we refer the reader to \cite{pitzbook}.
 
In particular, following the notation from \cite{pitzbook}, the order type \defn{$\mathrm{ht}(t)$} of $\odowncl{t}$ is called the \defn{height} of the node $t$. 
Moreover, we denote the \defn{up-closure} of $t$ as the set \defn{$\upcl{t}$}${} \coloneq \set{s\in T}:{s\geq t}$, denote the \defn{down-closure} of $t$ as \defn{$\downcl{t}$}${} \coloneq \set{s\in T}:{s\leq t}$ and call the elements from \defn{$\successor{t}$}${} \coloneq \set{s \in \upcl{t}\setminus \Set{t}}:{\text{there exists no }y \in \upcl{t} \text{ with } t < y < s}$ the \defn{successors} of $t$.
The \defn{height} of $T$ corresponds to the ordinal \defn{$\mathrm{ht}(T)$}${} \coloneq \sup_{t\in T}(\mathrm{ht}(t)+1)$. 
For a given $t\in T$, we also define \defn{$\oupcl{t} \coloneq \upcl{t}\setminus \Set{t}$} and, when $\mathrm{ht}(t)$ is a successor ordinal, the corresponding maximum element \defn{$\mathrm{pred}(t)$} of $\odowncl{t}$ is called the \defn{predecessor} of $t$.  
For a subset $X\subseteq T$, we also fix the notation \defn{$\upcl{t}_X$} as a shorthand for $\upcl{t}\cap X$.
We say that $X$ is \defn{convex} or a \defn{subtree} of $T$ if, for every pair $s,t \in X$ and $s' \in T$ such that $s \leq s' \leq t$, we also have $s' \in X$.

A totally ordered subset $C \subseteq T$ under $\leq$ is called a \defn{chain}, and it is a \defn{branch} of $T$ if it is $\subseteq$--maximal regarding this property. 
We also denote \defn{$\Above C$}${} \coloneq \set{t\in T\setminus C}:{t>s\text{ for every }s\in C}$, while we refer to the convex chains of $T$ as its \defn{intervals}.
Typical examples of intervals are sets of the form \defn{$[s,t)$}${} \coloneq \set{s'\in T}:{s\leq s'<t}$ and \defn{$[s,t]$}${} \coloneq \set{s'\in T}:{s\leq s'\leq t}$ for some pair $s,t \in T$ with $s \leq t$, as well as $\upcl{t}_C$ if $C \subseteq T$ is a branch and $t \in C$.
A non-empty set $X \subseteq T$ is \defn{down-closed} if it satisfies $\downcl{t} \subseteq X$ for every $t \in X$, i.e., if it is convex and contains the root of $T$.
On the other hand, $A \subseteq T$ is called an \defn{antichain} if it comprises pairwise incomparable elements under $\leq$.
For example, the $i$th \defn{level} $\set{t\in T}:{\mathrm{ht}(t) = i}$ is an antichain in $T$ for each $i<\mathrm{ht}(T)$.

Note that order trees of height at most $\omega$ are precisely the partial orders arising as the tree-orders of rooted graph-theoretical trees.
In this case, the well-known König's theorem states that, when infinite, such trees have a node with infinitely many successors or an infinite branch. 
For arbitrary order trees, a similar dichotomy reads as follows:

 \begin{lemma}\label{claim:antichainKappa}
    If $T$ is an order tree of size at least $\kappa$ but height $< \kappa$, then there exists an antichain in $T$ of size $\kappa$.
\end{lemma}

\begin{proof}
    If $\kappa$ is regular, then some level of $T$ is of size at least $\kappa$ and provides an antichain of size $\kappa$.
    If $\kappa$ is singular, the same argument shows that, for every regular cardinal $\mu < \kappa$, some level of $T$ has cardinality at least $\mu$.
    However, it follows from the classical work of Erd\H{o}s and Tarski in \cite{ErdosTarski1943} that
        \[
            \psi \coloneq \min{\set{\mu}:{\text{there exists no antichain in $T$ of size }\mu}}
        \]
    is a regular cardinal.
    Then, the singularity of $\kappa$ implies that $\psi > \kappa$ and, thus, the existence of an antichain of size $\kappa$ in $T$.
\end{proof}

Subsequent arguments in this paper might refer to different order trees $(T,\leq)$ and $(T',\preceq)$ on perhaps the same ground set: in these cases, the above fixed notation is decorated with subscripts such as $\downcl{\cdot}_{\leq}$ and $\downcl{\cdot}_{\preceq}$ for distinguishing between the corresponding partial orders.

\subsection{Normality definition and properties}\label{subsection:NPTdefinition}

As a main tool for the forthcoming arguments in this paper, this section recalls the definitions of \textit{$T$--graphs} and \textit{normal partition trees} as first introduced by Brochet and Diestel in \cite{BrochetDiestel}.
Following the conventions above, let $(T,\leq)$ be a given order tree. 
A graph $G$ defined on $V(G)= T$ is a \defn{$T$--graph} if the following properties are verified:

\begin{enumerate}[label=$(\square_{\arabic*})$]
    \item all edges of $G$ have comparable endpoints regarding $\leq$; and
    
    \item\label{item:CofinalNeighborhoods} the neighborhood $N(t)\cap \odowncl{t}$ of every node $t\in T$ is cofinal in the chain $\odowncl{t}$.
\end{enumerate}

Additionally, for a given cardinal $\kappa$, the graph $G$ is said to satisfy the

\begin{enumerate}[label=$(\diamond_{\arabic*})$]
    \item \label{prop:DownClosedNeighborhood} \defn{down-closed neighborhood property} if, for every $t\in T$, the intersection $\upcl{t}_{N(t) \cup \Set{t}}$ is a convex subset of $T$ (with minimum element $t$);

      \item \label{prop:shortNeighborhoodProp} \defn{$\kappa$--short neighborhood property} if, for every $t \in T$ and convex $\Delta\subseteq \upcl{t}_{N(\odowncl{t})}$ of size at least $\kappa$, there is an antichain $A \subseteq \Delta$ with $\cardinality{\upcl{a}_{\Delta}} < \kappa$ for every $a \in A$ but $\sum_{a\in A} \cardinality{\upcl{a}_{\Delta}} = \cardinality{\bigcup_{a\in A}\upcl{a}_{\Delta}} \geq \kappa$.
\end{enumerate}

When $\kappa$ is understood from the context (for instance, if it is known that $G$ is $\kappa$--edge-connected), we shall refer to the second property \cref{prop:shortNeighborhoodProp}  just as the \textit{short neighborhood property}. 
Combined with the down-closed neighborhood property, and since the above second item \cref{prop:shortNeighborhoodProp} cannot hold for a $\kappa$--sized chain, it follows that the convex set $\upcl{t}_{N(\odowncl{t})}$ cannot have an interval of size $\geq \kappa$. 
Incidentally, this also provides $\kappa$ as a bound for the height of $T$: after all, if there is a node $s\in T$ of height $\kappa$, any neighbor $t\in N(s)\cap \odowncl{s}$ is such that $N(t)$ contains the $\kappa$--sized interval $(t,s)$.

Naturally, a set $X \subseteq T$ is said to be \defn{connected} if it induces a connected graph in the underlying $T$--graph $G$. 
When $T$ has height at most $\omega$, a $T$--graph is simply a graph that admits $T$ as a normal spanning tree.
Following Brochet and Diestel, a \defn{normal partition tree} $(T,\mathcal{V})$ of a graph $G$ is a partition $\mathcal{V} \coloneq \set{V_t}:{t\in T}$ of $V(G)$, indexed by the elements of an order tree $T$, satisfying the following axioms:

\begin{itemize}
    \item $G[V_t]$ is a connected subgraph of $G$ for every $t \in T$;
    \item the contraction minor $\Dot{G}\coloneq G/\mathcal{V}$ (where $V_t$ is contracted to the corresponding node $t$, and any resulting parallel edges are ignored) is a $T$--graph.
\end{itemize}

Below, we summarize four useful properties of $T$--graphs that highlight the behavior of $T$ as a normal spanning \textit{order} tree; see \cite{BrochetDiestel}*{\S$2$} for details.

\begin{lemma}\label{lemma:PropertiesTGraphs}
    Let $(T, \leq)$ be an order tree and let $G$ be a $T$--graph.
    Then,
    \begin{itemize}
        \item The interval $[s,t)$ induces a connected subgraph of $G$ for every $s,t\in T$ with $s<t$;
        
        \item If $s,t\in T$ are incomparable, then $s$ and $t$ belong to distinct components of $G - (\odowncl{s}\cap \odowncl{t})$;
        
        \item Every connected subgraph $C$ of $G$ has a unique $\leq$--minimal element $\min C$;
        
        \item If $X \subseteq T$ is down-closed, then every component $C$ of $G - X$ is spanned by $\upcl{\min C}$.
    \end{itemize}
\end{lemma}

In a slightly more general setting, we might consider the case in which, rather than being a partition of the whole vertex set $V(G)$, the family $\mathcal{V} = \set{V_t}:{t\in T}$ comprises just pairwise disjoint connected sets of vertices of $G$. 
Following the work of Lorenz and Pitz in \cite{lorenz2025obstructions}, $(T,\mathcal{V})$ is said to be a \defn{normal semi-partition tree} of $G$ if the following condition is verified: if $P$ is a path in $G$ which meets $\bigcup \mathcal{V}$ precisely in its endpoints $x$ and $y$, then the nodes $s,t\in T$ such that $x\in V_s$ and $y\in V_t$ are comparable with respect to $\leq$.
As an equivalent way of stating this property, we have:

\begin{remark}\label{remark:normality}
    $(T,\mathcal{V})$ is normal in $G$ if, and only if, the contraction minor $(\bigcup \mathcal{V})/\mathcal{V}$ is a $T$--graph and, for every component $D$ of $G - (\bigcup \mathcal{V})$, the set $\set{t\in T}:{V_t\cap N(D)\neq \emptyset}$ is a chain in $T$. 
\end{remark}

While $N(\cdot)$ refers to the usual neighborhood of vertices and sets of vertices of $G$, for $t \in V(G/\mathcal{V})$ or $X\subseteq V(G/\mathcal{V})$ we denote the neighborhood of $t$ and of $X$ in $G/\mathcal{V}$ by \defn{$\Dot{N}(t)$} and \defn{$\Dot{N}(X)$} respectively.
Furthermore, we extend the notation of $\Dot{N}(\cdot)$ to sets $Y \subseteq V(G)$, where $\Dot{N}(Y) \coloneq  \set{t\in T}:{N(Y) \cap V_t \neq \emptyset}$.
Both notations coincide when $G$ \textit{is} the contraction minor $G/\mathcal{V}$.

\subsection{A standard normal partition tree construction}\label{subsection:NPTconstruction}
In this subsection we work towards the following lemma, which is the core tool for supporting the proof sketch of our main result as outlined in the introduction.

\begin{lemma}\label{lem:reductionToTgraphs}
    Let $\kappa$ be an uncountable cardinal and let $G$ be a $\kappa$--edge-connected graph.
    Then either 
    
    \begin{itemize}
        \item there exists a set $X \subseteq V(G)$ of size $< \kappa$ such that some component of $G-X$ has no $(N(X),\kappa)$--rank, or
        
        \item $G$ admits a normal partition tree $(T, \cV)$ with $\mathrm{ht}(T)\leq \kappa$, whose branch sets have size $<\kappa$ and such that the contraction minor $G/ \cV$ is $\kappa$--edge-connected and satisfies the down-closed and short neighborhood properties.
    \end{itemize}
\end{lemma}

Indeed, both alternatives of \cref{lem:reductionToTgraphs} are possible. 
For instance, the clique $K_{\kappa}$ of course addresses the first bullet point, since the neighborhood of every vertex contains an isomorphic copy of $T_{\kappa}$ and, according to \cref{lemma:Arank}, can have no corresponding localized rank. 
On the other hand, the next example concludes that the small subsets of $\TkM$ have \textit{dispersed} neighborhoods, despite the natural copy of $T_{\kappa}$ in $\TkM$ being a normal spanning tree for this graph and addressing the down-closed and short neighborhood properties:

\begin{example}[Revisited interleaved tree]\label{example:LackOfSpreadingPair}
     Let $\TkM$ be the graph formally defined in \cref{example:interleaved} together with its fixed normal spanning tree $T\subseteq \TkM$ isomorphic to $T_{\kappa}$. 
     We will show that, for every set $X \subseteq T$ of size $<\kappa$, any component $C$ of $G-X$ admits an $N(X)$--rank. 
     Indeed, by replacing $X$ by $\downcl{X}$ if necessary, we might assume that $X$ is down-closed regarding the tree-order $\leq$ of $T$: after all, heights of nodes in $T$ are finite (so that $\cardinality{\downcl{X}} = \cardinality{X}<\kappa$) and (rankless) components of $G-\downcl{X}$ are contained in (rankless) components of $G-X$. 
     Considering this simplification, and relying on \cref{lemma:PropertiesTGraphs}, $C$ can now be written as $C = \upcl{v}$ for some $v\in T$. 
     In particular, the nodes from $N(X)\cap C$ are contained in the subtree of $T_{\kappa}$ described by $N(\odowncl{v}) \cap \upcl{v} = \Set{v} \cup \successor{v}$, which is simply a star of $\kappa$ many leaves and centered at $v$. 
     Hence, $C$ has $N(X)$--rank at most $1$, since each connected component of $C\setminus \Set{v}$ contains at most one of such leaves by the normality of $T_{\kappa}$ on $\TkM$.
\end{example}

To address the second bullet point, however, we revisit a technique introduced independently by Robertson, Seymour and Thomas in \cite{RSTExcludingSubdivisions} and by the work of Brochet and Diestel in \cite{BrochetDiestel} for constructing normal partition trees.
This is a recursive depth-search construction of a sequence $\Set{(T_i,\mathcal{V}_i)}_{i < \kappa}$ of normal semi-partition trees of $G$, in which $\cardinality{V_t} \leq \mathrm{ht}(t) + \theta(\kappa)$ for every $i < \kappa$ and $t \in T_i$. 
For each $i<\kappa$, the subgraph of $G$ spanned by $\bigcup \mathcal{V}_i$ will be denoted by $G_i$.
Then, the ongoing algorithm considers a component $D$ of $G - G_i$ and assesses, for some node $t\in T_i$, whether it admits an $N(A_t)$--rank corresponding to the connected set $A_t\coloneq \bigcup_{s< t}V_s\subseteq G_i$ (of size $<\kappa$).
If this is not the case, the proof of \cref{lem:reductionToTgraphs} is finished by the first bullet point as an output.
Otherwise, $G_i$ is extended to a supergraph $G_{i+1}$ meeting $D$ on, among other vertices, elements of an $N(A_t)$--rank reducing set $X_{t,D} \subseteq D$. 
At the end of this procedure, $T\coloneq \bigcup_{i<\kappa}T_i$ will fit as an index tree for a normal partition $\mathcal{V}\coloneq \bigcup_{i<\kappa}\mathcal{V}_i$ of $G$ witnessing the second item of \cref{lem:reductionToTgraphs}.
In order to detail this routine, we set $I_{\kappa}\coloneq \Set{\kappa}$ and $\theta(\kappa)\coloneq \aleph_0$ if $\kappa$ is regular. 
For singular values of $\kappa$, the set $I_{\kappa}\subseteq \kappa$ is a fixed cofinal sequence of cardinality $\theta(\kappa)\coloneq \cof{\kappa}$ and comprising uncountable regular cardinals.
As the main inductive assumptions regarding the sequence $\Set{(T_i,\mathcal{V}_i)}_{i<\kappa}$, we will ensure that, given $i<j<\kappa$,

\begin{enumerate}[label=$(P_{\arabic*})$]
    \item\label{item:NPTheight} $T_i$ has height at most $i$;
    
    \item $T_i$ is the restriction of $T_j$ to all its nodes of height less than $i$;
    
    \item\label{item:ExtendingEachOther} $G_i \subseteq G_j$ and $\mathcal{V}_i \subseteq \mathcal{V}_j$. 
    In particular, due to the previous item, there is no ambiguity when writing $\mathcal{V}_i \coloneq \set{V_t}:{t\in T_i}$;
    
    \item\label{item:NPTnarrow} $V_t$ is connected and $\cardinality{V_t}\leq \cardinality{\mathrm{ht}(t)}+\theta(\kappa)$ for every $t\in T_i$;
    
    \item\label{item:NPTcomponent}  
    for every component $D$ of $G-G_i$, the set $\Dot{N}(D) = \set{t \in T_i}:{V_t \cap N(D) \neq \emptyset}$ is cofinally contained in a branch $B_D \subseteq T_i$ of order type $i$ and we have $G_{i+1} \cap D \neq \emptyset$. 
    In particular, due to \cref{remark:normality}, $(T_i,\mathcal{V}_i)$ is normal in $G$.
\end{enumerate}

We start this recursive construction by letting $T_1\coloneq \Set{r}$, while the corresponding graph $G_1$ is a single fixed vertex $v_* \in V(G)$ and so that $V_r \coloneq \Set{v_*}$.
If $j$ is a limit ordinal, then the definition of $(T_j, \cV_j)$ and $G_j$ is given by the standard limit object obtained from the previously defined normal partition trees $(T_i, \cV_i)$ of $G_i$, $i<j$.

\begin{lemma}\label{lemma:LimitConstruction}
    Let $j\leq \kappa$ be a limit ordinal and suppose $(T_i, \mathcal{V}_i)$ and $G_i$ are defined for every ordinal $i<j$ and satisfy \cref{item:NPTheight}--\cref{item:NPTcomponent}.
    Then $(\mathcal{V}_j \coloneq \bigcup_{i<j}\mathcal{V}_i, T_j \coloneq \bigcup_{i<j} T_i)$ defines a normal partition tree for the graph $G_j \coloneq \bigcup_{i<j} G_i$ and satisfies \cref{item:NPTheight}--\cref{item:NPTcomponent}. 
\end{lemma}

\begin{proof}
    It is immediate from the definition that $T_j$ has height at most $\sup_{i<j} \mathrm{ht}(T_i) \leq j$, and that $\set{V_t}:{t\in T_j}$ is indeed a partition of the graph $G_j = \bigcup_{i<j} G_i$.
    In addition, $T_i$ is clearly a restriction of $T_j$ to its nodes of height at most $i$ for every $i<j$, while $V_t$ is connected and $\cardinality{V_t} \leq \cardinality{\mathrm{ht}(t)} + \theta(\kappa)$ for every $t \in T_j$ by inductive assumption. 
    Similarly, in the contraction minor $\Dot{G_j} \coloneq G_j/\mathcal{V}_j$, any given edge has both endpoints in $G_i/\mathcal{V}_i$ for some large enough index $i<j$.
    In particular, both endpoints are comparable in the order of $T_i$ and, thus, in the order of $T_j$.
    Analogously, for any $t\in T_j$ of limit height there exists an index $i<j$ for which $t \in T_i$ is also a node of the same limit height.
    It follows that the neighborhood $\Dot{N}(t)$ of $t$ in $\Dot{G_i}$ is cofinal in $\odowncl{t}$ by inductive assumption, and hence the neighborhood of $t$ in $\Dot{G_j}$ is cofinal in $\odowncl{t}$ as well.
    Finally, the normality of $(T_i,\mathcal{V}_i)$ in $G$ is ensured for every $i<j$ by induction, from which, by the definition of $\mathcal{V}_j$ as $\bigcup_{i<j} \mathcal{V}_i$, it follows that $(T_j,\mathcal{V}_j)$ is normal in $G$.
    Incidentally, if $D$ is a component of $G - G_j$ and $B_j \subseteq T_j$ denotes the chain in which $\Dot{N}(D)$ is cofinal, then $B_j$ must have order type $j$: otherwise, $D$ would already be a component of $G-G_i$ for the order type $i<j$ of $B_j$, thus contradicting the inductive assumption that $D\cap G_{i+1}\neq \emptyset$.
    To summarize, $G_j/\mathcal{V}_j$ is indeed a $T_j$--graph and the normal partition tree $(T_j,\mathcal{V}_j)$ satisfies properties \cref{item:NPTheight}--\cref{item:NPTnarrow}, where \cref{item:NPTcomponent} now follows from \cref{remark:normality}.
\end{proof}

Given $j < \kappa$, it remains to describe the construction of $(T_{j+1},\mathcal{V}_{j+1})$ and $G_{j+1}$, given that $(T_{i},\mathcal{V}_i)$ and $G_i$ are defined for every $i \leq j$.
To that end, let $D$ be a component of $G - G_j$ and let $B_D$ be the corresponding branch in $T_j$ of order type $j$ and which contains $\Dot{N}(D)$  cofinally (due to \cref{item:NPTcomponent}).
By induction on $j$, we may additionally assume that, for each pair of nodes $s,t \in B_D\subseteq T_j$ with $s \leq t$, we have fixed a set $X_{s,t} \subseteq V(G)$ of size $<\kappa$ and which is endowed with a well-ordering $\preceq_{s,t}$ of order type $\cardinality{X_{s,t}}$.   
Then $T_{j+1}$ is the extension of $T_j$ obtained by adding to $T_j$, for each connected component $D$ of $G-G_j$, a node $t_D$ satisfying $t_D > u$ for every $u \in B_D$.

Now, let $t \in B_D \cup \Set{t_D}$.
Since $(T_j,\mathcal{V}_j)$ is a normal partition tree for $G_j\coloneq G[\bigcup \mathcal{V}_j]$, it follows that $A_{t} \coloneq \bigcup_{s< t}V_s$ is a connected set of vertices of $G$ whose size is bounded by $\cardinality{j}+\theta(\kappa)<\kappa$ due to item \cref{item:NPTnarrow}.
We might then assume that $D$ admits an $N(A_t)$--rank: otherwise, there exists a component $C$ of $G - A_t$ with $C \supseteq D$ which by \cref{lemma:RankReducing} admits no $N(A_t)$--rank, hence addressing the first bullet point of \cref{lem:reductionToTgraphs}.
If this rank is $0$, set $X_{t,t_D} = \emptyset$.
Otherwise, fix a corresponding separator $X_{t,t_D} \subseteq D$ of size $< \kappa$ for which all components of $D - X_{t,t_D}$ have $N(A_t)$--rank smaller than the $N(A_t)$--rank of $D$. 
In this case, let $\preceq_{t,t_D}$ denote any well-ordering of $X_{t,t_D}$ whose order type is $\cardinality{X_{t,t_D}}$.
To define $V_{t_D}$, we now recall the following result from previous work of the authors:

\begin{lemma}[\cite{OrderCompatible}, Lemma 2.1]\label{lemma:AuxiliaryPaths}
    If $H$ is a $\kappa$--edge-connected graph, then, for every uncountable regular\footnote{Compared to the original statement from \cite{OrderCompatible}, the reader should keep in mind that the $\lambda$--edge-connectivity of two vertices $u$ and $v$ can be witnessed by a family of $\lambda$ edge-disjoint $u$--$v$ paths of the same length.} cardinal $\lambda \leq \kappa$ and every pair of distinct vertices $u,v \in V(H)$, there is a sequence $S(u,v,\lambda) \coloneq (u_i)_{0 \leq i \leq n}$ with $u_0 = u$, $u_n = v$, such that $u_i$ and $u_{i+1}$ are at least $\lambda$--vertex-connected for every $0\leq i < n$. 
\end{lemma}

Now fix an arbitrary $u_{t_D} \in D$ and, because $G$ is $\kappa$--edge-connected, for each $\lambda \in I_{\kappa}$ let $S(u_{t_D},v_*,\lambda) \subseteq V(G)$ be a finite set of vertices as claimed by \cref{lemma:AuxiliaryPaths} (where $v_*$ is the unique vertex in $V_r$). 
The set $X_{t_D} \coloneq \set{\min_{\preceq_{s,t}}X_{s,t}\cap D}:{s\leq t\leq t_D}$ has size at most $\cardinality{j}+\aleph_0$, while $Y_{t_D} \coloneq D \cap \bigcup_{\lambda \in I_{\kappa}}S(u_{t_D},v_*,\lambda)$ has size at most $\theta(\kappa)$.
Finally, for each $t \in \Dot{N}(D) \cap B_D$, we fix a vertex $v_t^D \in N(V_t) \cap D$ witnessing the fact that $N(D) \cap V_t \neq \emptyset$. 
We then denote $Z_{t_D} \coloneq \set{v_t^D}:{t\in \Dot{N}(D)\cap B_D}$ and observe that $\cardinality{Z_{t_D}} \leq \cardinality{B_D} \leq \cardinality{j}$ by item \cref{item:NPTheight}.
Since $D$ is connected, we can then fix a (non-empty) connected superset 
\begin{center}
    \begin{varwidth}{\linewidth}
        \begin{enumerate}[label=$(*)$, leftmargin=*]
            \item\label{item:VstarSet} $V_{t_D} \supseteq X_{t_D} \cup Y_{t_D} \cup Z_{t_D}$
        \end{enumerate}
    \end{varwidth}
\end{center}
contained in $D$ and of size at most $\cardinality{j} + \theta(\kappa)$.
As we shall see later, the definition of $\set{X_t}:{t \in T}$ is tailored to ensure the short neighborhood property of the final normal partition tree $(T,\mathcal{V})$, while $\set{Y_t}:{t \in T}$ will provide the $\kappa$--edge-connectivity of the arising $T$--graph $\Dot{G}$. 
Finally, the family $\set{Z_t}:{t \in T}$ ensures the down-closed neighborhood property of $(T,\mathcal{V})$. 
This finishes the description of $\mathcal{V}_{j+1} \coloneq \set{V_t}:{t \in T_{j+1}}$ as a normal partition tree of $G_{j+1} \coloneq \bigcup \mathcal{V}_{j+1}$.

Indeed, $\Dot{G}_{j+1}\coloneq G_{j+1}/\mathcal{V}_{j+1}$ is a $T_{j+1}$--graph since $G_j/\mathcal{V}_j$ is a $T_j$--graph and because, if $D$ is a component of $G - G_j$ as previously considered, then:

\begin{enumerate}
    \item the edges of $\Dot{G}_{j+1}$ having one endpoint in $V_{t_D}$ have the other endpoint in $\bigcup_{t \in B_D} V_t$, so that the corresponding endpoints in $\Dot{G}_{j+1}$ are comparable in $T_{j+1}$ by construction, and
    
    \item  item \cref{item:NPTcomponent} claims that the neighborhood of $t_D$ in $\Dot{G}_{j+1}$ is cofinal in $B_D$.
\end{enumerate}

In its turn, the connectedness of $D$ ensures that every component of $G - G_{j+1}$ contained in $D$ must have a neighbor in $V_{t_D}$.
Moreover, since $t_D$ is the maximum element of the branch in $T_{j+1}$ (of order type $j+1$), properties \cref{item:NPTheight} and \cref{item:NPTcomponent} are verified by induction. 
In fact, any component $D'$ of $G - G_{j+1}$ is obtained as a component of $D - V_{t_D}$ for some $D \subseteq G - G_j$ as above, so that $N(D') \subseteq N(D) \cup V_{t_D}$ and $N(D') \cap V_{t_D} \neq \emptyset$, thus implying $\Dot{N}(D') \subseteq B_D \cup \Set{t_D}$ and $t_D \in \Dot{N}(D')$.
In particular, $\Dot{N}(D')$ is cofinal in the branch $B_{D} \cup \Set{t_D}$ of order type $j+1$ in $T_{j+1}$.
In the remainder of this section, we will show that $(T,\mathcal{V}) \coloneq (T_{\kappa},\mathcal{V}_{\kappa})$ as defined in \cref{lemma:LimitConstruction} is indeed a normal partition tree of $G$ satisfying the second bullet point from \cref{lem:reductionToTgraphs}.
To that end, we first summarize how the nodes of $T$ and the components of $G - G_i$ interact.

\begin{remark}\label{remark:DescriptionNodes}
    Every node $s \in T \setminus \Set{r}$ at height $i \coloneq \mathrm{ht}(s)$ can be identified with the component $D_s$ of $G - G_i$ for which $s = t_{D_s}$.
    These components are nested along $T$ as follows:
    \begin{itemize}
        \item if $s < s'$, then $D_{s'} \subseteq D_s$ and

        \item if $s$ and $s'$ are incomparable, then $D_s \cap D_{s'} = \emptyset$.
    \end{itemize}

   The first statement follows by induction on $\mathrm{ht}(s')$ together with property \cref{item:NPTcomponent}, while the second statement follows from the normality of $T_{j+1}$ for $j \coloneq \max\Set{\mathrm{ht}(s),\mathrm{ht}(s')}$.
\end{remark}

\begin{lemma}\label{lemma:DownClosedNeighborhoods}
    The graph $\Dot{G}_\kappa \coloneq \bigcup_{j < \kappa} \Dot{G}_j$ is a $T$--graph with the down-closed and short neighborhood properties.
\end{lemma}

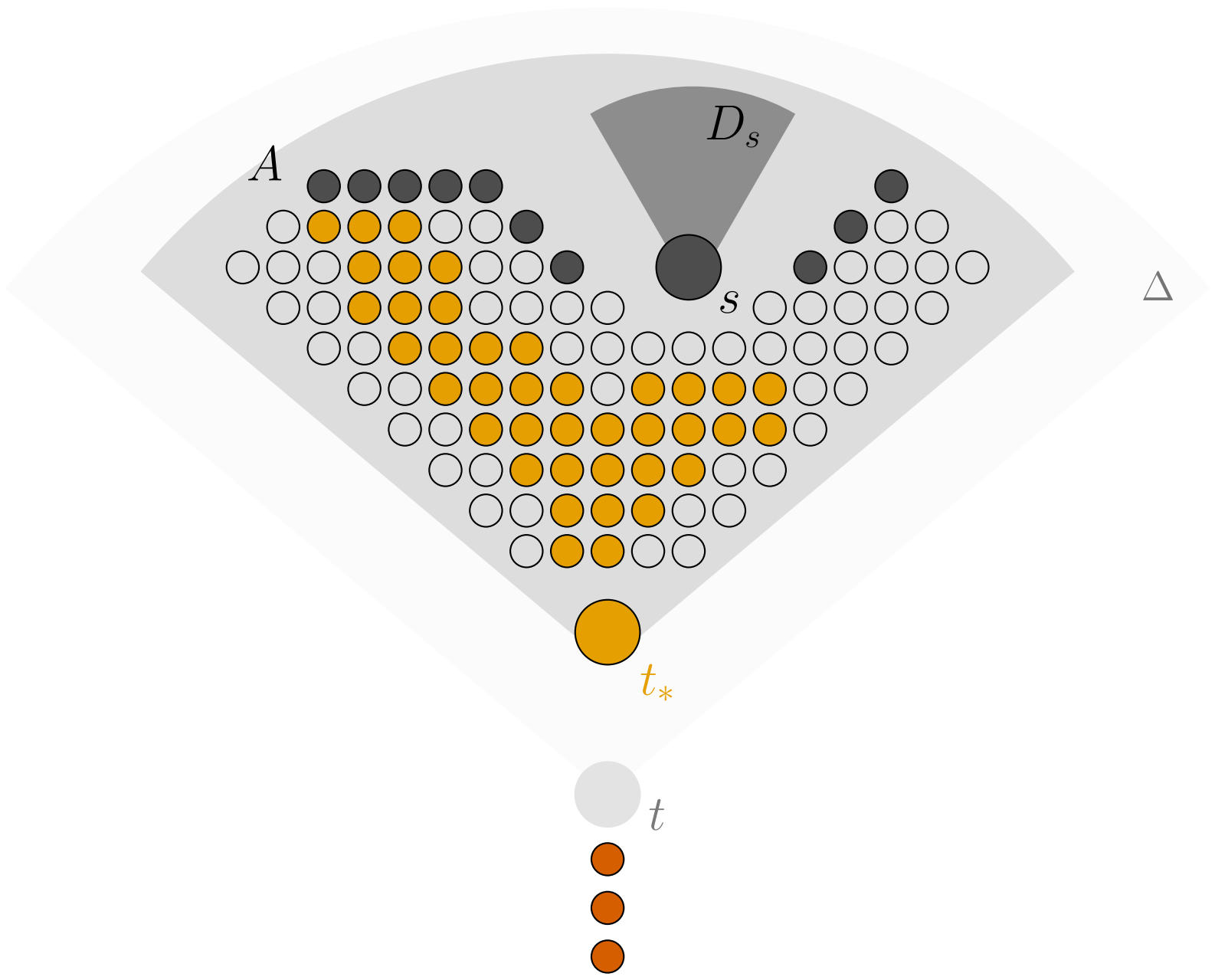
\begin{figure}
    \centering
    \begin{tikzpicture}[scale=0.7, transform shape]

 \node[shape=circular sector, rotate=-90, circular sector angle=100, minimum width=525pt, fill=white!92!black, fill opacity=0.2] at (20,15.93) {};
 \node[shape=circular sector, rotate=-90, circular sector angle=100, minimum width=176pt, minimum height=214pt, fill=black!21, fill opacity=0.6] at (20,16.2) {};
 \node[shape=circular sector, minimum height=72pt, rotate=-90, fill=gray!96!black, fill opacity=0.8] at (21.05,17.99) {};
 \draw[fill=okvermillion] (20,8) circle (0.2cm);
 \draw[fill=okvermillion] (20,8.6) circle (0.2cm);
 \draw[fill=okvermillion] (20,9.2) circle (0.2cm);
 \draw[fill=okvermillion, black!11] (20,10) circle (0.4cm);
 \draw (19,13) circle (0.2cm);
 \draw[fill=okorange] (19.5,13) circle (0.2cm);
 \draw[fill=okorange] (20,13) circle (0.2cm);
 \draw (20.5,13) circle (0.2cm);
 \draw (21,13) circle (0.2cm);
 \draw (19,13.5) circle (0.2cm);
 \draw[fill=okorange] (19.5,13.5) circle (0.2cm);
 \draw[fill=okorange] (20,13.5) circle (0.2cm);
 \draw[fill=okorange] (20.5,13.5) circle (0.2cm);
 \draw (21,13.5) circle (0.2cm);
 \draw (21.5,13.5) circle (0.2cm);
 \draw (18.5,13.5) circle (0.2cm);
 \draw[fill=okorange] (19,14) circle (0.2cm);
 \draw[fill=okorange] (19.5,14) circle (0.2cm);
 \draw[fill=okorange] (20,14) circle (0.2cm);
 \draw[fill=okorange] (20.5,14) circle (0.2cm);
 \draw[fill=okorange] (21,14) circle (0.2cm);
 \draw (21.5,14) circle (0.2cm);
 \draw (18.5,14) circle (0.2cm);
 \draw (22,14) circle (0.2cm);
 \draw (18,14) circle (0.2cm);
 \draw[fill=okorange] (19,14.5) circle (0.2cm);
 \draw[fill=okorange] (19.5,14.5) circle (0.2cm);
 \draw[fill=okorange] (20,14.5) circle (0.2cm);
 \draw[fill=okorange] (20.5,14.5) circle (0.2cm);
 \draw[fill=okorange] (21,14.5) circle (0.2cm);
 \draw[fill=okorange] (21.5,14.5) circle (0.2cm);
 \draw[fill=okorange] (18.5,14.5) circle (0.2cm);
 \draw[fill=okorange] (22,14.5) circle (0.2cm);
 \draw (18,14.5) circle (0.2cm);
 \draw (22.5,14.5) circle (0.2cm);
 \draw (17.5,14.5) circle (0.2cm);
 \draw[fill=okorange] (19,15) circle (0.2cm);
 \draw[fill=okorange] (19.5,15) circle (0.2cm);
 \draw (20,15) circle (0.2cm);
 \draw[fill=okorange] (20.5,15) circle (0.2cm);
 \draw[fill=okorange] (21,15) circle (0.2cm);
 \draw[fill=okorange] (21.5,15) circle (0.2cm);
 \draw[fill=okorange] (18.5,15) circle (0.2cm);
 \draw[fill=okorange] (22,15) circle (0.2cm);
 \draw[fill=okorange] (18,15) circle (0.2cm);
 \draw (22.5,15) circle (0.2cm);
 \draw (17.5,15) circle (0.2cm);
 \draw (23,15) circle (0.2cm);
 \draw (17,15) circle (0.2cm);
 \draw[fill=okorange] (19,15.5) circle (0.2cm);
 \draw (19.5,15.5) circle (0.2cm);
 \draw (20,15.5) circle (0.2cm);
 \draw (20.5,15.5) circle (0.2cm);
 \draw (21,15.5) circle (0.2cm);
 \draw (21.5,15.5) circle (0.2cm);
 \draw[fill=okorange] (18.5,15.5) circle (0.2cm);
 \draw (22,15.5) circle (0.2cm);
 \draw[fill=okorange] (18,15.5) circle (0.2cm);
 \draw (22.5,15.5) circle (0.2cm);
 \draw[fill=okorange] (17.5,15.5) circle (0.2cm);
 \draw (23,15.5) circle (0.2cm);
 \draw (17,15.5) circle (0.2cm);
 \draw (23.5,15.5) circle (0.2cm);
 \draw (16.5,15.5) circle (0.2cm);
 \draw[fill=black!70] (21,16.5) circle (0.4cm);
 \draw (19,16) circle (0.2cm);
 \draw (19.5,16) circle (0.2cm);
 \draw (20,16) circle (0.2cm);
 \draw (18.5,16) circle (0.2cm);
 \draw[fill=okorange] (18,16) circle (0.2cm);
 \draw[fill=okorange] (17.5,16) circle (0.2cm);
 \draw[fill=okorange] (17,16) circle (0.2cm);
 \draw (16.5,16) circle (0.2cm);
 \draw (16,16) circle (0.2cm);
 \draw (22,16) circle (0.2cm);
 \draw (22.5,16) circle (0.2cm);
 \draw (23,16) circle (0.2cm);
 \draw (23.5,16) circle (0.2cm);
 \draw (24,16) circle (0.2cm);
 \draw (19,16.5) circle (0.2cm);
 \draw[fill=black!70] (19.5,16.5) circle (0.2cm);
 \draw (18.5,16.5) circle (0.2cm);
 \draw[fill=okorange] (18,16.5) circle (0.2cm);
 \draw[fill=okorange] (17.5,16.5) circle (0.2cm);
 \draw[fill=okorange] (17,16.5) circle (0.2cm);
 \draw (16.5,16.5) circle (0.2cm);
 \draw (16,16.5) circle (0.2cm);
 \draw[fill=black!70] (22.5,16.5) circle (0.2cm);
 \draw (23,16.5) circle (0.2cm);
 \draw (23.5,16.5) circle (0.2cm);
 \draw (24,16.5) circle (0.2cm);
 \draw (24.5,16.5) circle (0.2cm);
 \draw (15.5,16.5) circle (0.2cm);
 \draw[fill=okorange] (17.5,17) circle (0.2cm);
 \draw[fill=okorange] (17,17) circle (0.2cm);
 \draw[fill=okorange] (16.5,17) circle (0.2cm);
 \draw (16,17) circle (0.2cm);
 \draw[fill=black!70] (23,17) circle (0.2cm);
 \draw (23.5,17) circle (0.2cm);
 \draw (24,17) circle (0.2cm);
 \draw (18,17) circle (0.2cm);
 \draw (18.5,17) circle (0.2cm);
 \draw[fill=black!70] (19,17) circle (0.2cm);
 \draw[fill=black!70] (17.5,17.5) circle (0.2cm);
 \draw[fill=black!70] (17,17.5) circle (0.2cm);
 \draw[fill=black!70] (16.5,17.5) circle (0.2cm);
 \draw[fill=black!70] (23.5,17.5) circle (0.2cm);
 \draw[fill=black!70] (18,17.5) circle (0.2cm);
 \draw[fill=black!70] (18.5,17.5) circle (0.2cm);
 \draw[fill=okorange] (20,12) circle (0.4cm);

  \node[draw=none, node font=\LARGE, text=black!52] at (20.61,9.75) {$t$};
 \node[draw=none, node font=\LARGE] at (21.5,16.07) {$s$};
   \node[draw=none, node font=\LARGE, text=okorange] at (20.62,11.37) {$t_*$};
  \node[draw=none, minimum width=45pt, node font=\Large, text=black!54] at (26.79,16.27) {$\Delta$};
  \node[draw=none, node font=\LARGE] at (21.55,18.24) {$D_s$};
  \node[draw=none, node font=\LARGE] at (15.77,17.78) {$A$};
\end{tikzpicture}

\caption{The bags in red define the set $A_{t}$, while the bags in orange denote those meeting the corresponding $N(A_t)$--rank reducing set $X_{t,t_*}$ for the component $D_{t_*}$. 
In particular, the well-ordering $\preceq_{t,t_*}$ ensures that there is no such bag whose index node has height greater than $\mathrm{ht}(t_*) +\cardinality{X_{t,t_*}}$. 
Similarly, if $A$ is an antichain in $\Delta$ as represented by the black circles, at most $\lambda = \cardinality{X_{t,t_*}}$ of their corresponding components intersect $X_{t,t_*}$.}
        \label{fig:ProofShortneighborhood}
    \end{figure}

\begin{proof}
    Let $t \in T$ be any node.
    In order to first verify that $\Dot{G}_{\kappa}$ satisfies the down-closed neighborhood property, fix $s \in \upcl{t}_{N(t) \cup \Set{t}}$ and some $s' \in T$ such that $t < s' < s$.
    By the above \cref{remark:DescriptionNodes}, we can consider the corresponding indices $i',i < \kappa$ for which there are components $D'$ and $D$ of $G - G_{i'}$ and $G - G_i$, respectively, such that $s' \coloneq t_{D'}$ and $s \coloneq t_D$.
    Moreover, the inequality $t_{D'} = s'<s = t_D$ claims further that $i'<i$ and $D\subseteq D'$.
    Since $s$ is a neighbor of $t$ in $\Dot{G}_{\kappa}$, the intersection $N(V_t)\cap V_s$ is non-empty, which implies now that $N(V_t)\cap D \neq \emptyset$ and, hence, $N(V_t) \cap D'\neq \emptyset$.
    In particular, $t \in \Dot{N}(D')$, ensuring by construction the existence of some $v_t^{D'} \in Z_{t_{D'}} \subseteq V_{t_{D'}}$ which is a neighbor of $V_t$ in $G$.
    The edge connecting $v_t^{D'}$ to $V_t$ in $G$ witnesses, in the contraction minor $\Dot{G}_{\kappa}$, that $s'=t_{D'} \in \upcl{t}_{N(t)}$.
    In other words, $\Dot{G}_{\kappa}$ verifies the down-closed neighborhood property as a $T$--graph.

    To verify the short neighborhood property, fix a convex $\Delta \subseteq \upcl{t}_{N(\odowncl{t})}$ of size at least $\kappa$. 
    If the up-closure $\upcl{s}_{\Delta}$ has size $< \kappa$ for every node $s \in \Delta$, then the set of all $\leq$--minimal elements in $\Delta$ provides an antichain satisfying the short neighborhood property.
    Otherwise, there exists at least one $t_* \in \Delta$ with $\cardinality{\upcl{t_*}_{\Delta}}\geq \kappa$, and due to \cref{remark:DescriptionNodes} there is some $i_{t_*}<\kappa$ such that $t_*= t_{D_{t_*}}$ for some component $D_{t_*}$ of $G - G_{i_{t_*}}$.
    In particular, the component $D_{t_*}$ has some $N(A_t)$--rank $> 0$ since $\cardinality{N(A_t) \cap D_{t_*}} \geq \kappa$ due to the definition of $A_t$ as $A_t = \bigcup_{s < t}V_s$ and the hypothesis that $\Delta \subseteq N(\odowncl{t})$.
    Then let $t_* \in \Delta$ be a node for which $\upcl{t_*}_\Delta$ has size at least $\kappa$ and which minimizes the $N(A_t)$--rank $\alpha > 0$ of $D_{t_*}$ among all possible choices.
    Recall that the corresponding separator $X_{t,t_*} \subseteq D_{t_*}$ is endowed with a well-ordering $\preceq_{t,t_*}$ of order type $\lambda \coloneq \cardinality{X_{t,t_*}} < \kappa$.
    The minimality of $\alpha$ now leads to the following observation.

    \begin{claim}\label{claim:neighborshort:1}
        There is no node $s \in \upcl{t_*}_{\Delta}$ of height strictly greater than $\mathrm{ht}(t_*)+ \max{\Set{\omega, \lambda}}$ (in $T$) and with $\cardinality{\upcl{s}_{\Delta}} \geq \kappa$.
    \end{claim}

    \begin{proof}[Proof of~\cref{claim:neighborshort:1}]\renewcommand{\qedsymbol}{$\blacksquare$}
        Suppose for a contradiction that such an $s \in \upcl{t_*}_{\Delta}$ exists. 
        In this case, it holds that $X_{t,t_*} \cap D_{s}=\emptyset$: otherwise, $X_{t,t_*}\cap D_{p}\neq \emptyset$ for every $t_*<p<s$ (since $D_{s}\subseteq D_p\subseteq D_{t_*}$), and in particular there exists a vertex $x_p \coloneq \min_{\preceq_{t,t_*}}D_p\cap X_{t,t_*}$ by definition of $X_p\subseteq V_p$ as in \cref{item:VstarSet}.
        Moreover, the minimal choice of $x_p$ ensures that the set $\set{x_p}:{t_*<p<s}$ is ordered via $\preceq_{t,t_*}$ in the same order type as the interval $(t_*,s) \subseteq \Delta$.
        However, this order type is strictly greater than $\lambda$ by the chosen height for $s$, contradicting the order type of $\preceq_{t,t_*}$.

        Summarized, it holds that $X_{t, t_*} \cap D_{s} = \emptyset$ and in fact $D_{s}$ is contained in some component of $D_{t_*} - X_{t,t_*}$. 
        \cref{lemma:RankReducing} then claims that the $N(A_t)$--rank of $D_{s}$ is strictly smaller than $\alpha$, contradicting the minimal choice of the node $t_*$.
    \end{proof}

    In its turn, the set $L \subseteq \upcl{t_*}_{\Delta}$ comprising the nodes of height $\mathrm{ht}(t_*) + \max{\Set{\omega, \lambda}} + 1$ in $T$ (and which could potentially be the empty set) clearly defines an antichain in $\Delta$.
    Hence, \cref{claim:neighborshort:1} finishes the verification of the short neighborhood property in the case that $\cardinality{ \bigcup_{a \in L} \upcl{a}_{\Delta} } \geq \kappa$. 
    If $\cardinality{ \bigcup_{a \in L} \upcl{a}_{\Delta} } < \kappa$, however, then $\upcl{t_*}_{\Delta}$ can be regarded as an order tree of size $\geq \kappa$ but of height at most $\lambda + \lambda'$, where $\lambda'$ is an ordinal of cardinality $\cardinality{ \bigcup_{a \in L} \upcl{a}_{\Delta}}$.
    In particular, $\lambda + \lambda' < \kappa$ and, by \cref{claim:antichainKappa}, the order tree $\upcl{t_*}_{\Delta}$ (of size at least $\kappa$ but height $< \kappa$) admits an antichain $A$ of size $\kappa$.

    Then, the proof of the short neighborhood property of $\Dot{G}_{\kappa}$ is finished if we show that $\cardinality{\upcl{a}_{\Delta}} \geq \kappa$ for fewer than $\kappa$ many nodes $a \in A$. 
    The formal argument below is also illustrated by the picture on the right side of \cref{fig:ProofShortneighborhood}.

    \begin{claim}\label{claim:neighborshort:2}
        The set $A' \coloneq \set{a\in A}:{\cardinality{\upcl{a}_{\Delta}} \geq \kappa}$ has size at most $\lambda < \kappa$.
    \end{claim}

    \begin{proof}[Proof of~\cref{claim:neighborshort:2}]\renewcommand{\qedsymbol}{$\blacksquare$}
        Since $A'$ is an antichain, the set $\set{D_s}:{s\in A'}$ comprises pairwise disjoint connected subgraphs of $G_{\kappa}$.
        Suppose towards a contradiction that $\cardinality{A'} > \lambda$, then there exists some $s \in A'$ such that $X_{t,t_*} \cap D_s = \emptyset$, because $\lambda = \cardinality{X_{t,t_*}}$. 
        As in the proof of \cref{claim:neighborshort:1}, $D_s$ is now contained in some component of $D_{t_*} - X_{t,t_*}$, thus having smaller $N(A_t)$--rank than $D_{t_*}$ due to \cref{lemma:RankReducing}.
        Since $s > t_* > t$ and $\cardinality{\upcl{s}_{\Delta}} \geq \kappa$ by definition of $A'$, this once more contradicts the minimal choice of $t_*$.
    \end{proof}
\end{proof}

As a consequence, the above \cref{lemma:DownClosedNeighborhoods} establishes the equality $\bigcup \mathcal{V} = G$ and, hence, that $(T,\mathcal{V})$ is a normal partition tree for $G = G_{\kappa}$ whose corresponding contraction minor $\Dot{G}\coloneq \Dot{G}_{\kappa}$ verifies the down-closed and short neighborhood properties. 
In fact, if there is a component $D$ of $G-G_{\kappa}$, \cref{item:NPTcomponent} ensures that $\Dot{N}(D)$ is contained in a branch $B$ of $T_{\kappa}$ of order type $\kappa$. 
Given $t\in \Dot{N}(D)$ and $t^+ \in B$ the unique successor of $t$ in $B$, the down-closed neighborhood property of $\Dot{G}_{\kappa}$ concludes that almost all nodes of $B$ are contained in the neighborhood of $t$; in particular $\upcl{t^+}_{N(\odowncl{t^+})}$ contains a subinterval $B' \subseteq B$ of size $\kappa$; but $B'$ contains no antichains satisfying the short neighborhood property \cref{prop:shortNeighborhoodProp}, a contradiction. 
It remains to verify that $G / \mathcal{V}$ is $\kappa$--edge-connected, thus addressing the last remaining requirement of the second bullet point from \cref{lem:reductionToTgraphs}.

\begin{lemma}\label{lemma:kEdgeConnected}
The contraction minor $\Dot{G} \coloneq G / \mathcal{V}$ is $\kappa$--edge-connected.
\end{lemma}

\begin{proof}
    First, $\Dot{G}$ contains at least two vertices because $G$ has size at least $\kappa$ but $\cardinality{V_t}<\kappa$ for every $t\in T$.
    Since $\kappa = \sup I_{\kappa}$, it suffices to show that, for every $t\in T$ and $\lambda_0 \in I_{\kappa}$, the node $t$ is $\lambda_0$--edge-connected to the root $r\in T$.
    Similarly, by induction on $i \coloneq \mathrm{ht}(t)$, this follows if we are able to show that there is some node $s < t$ which is $\lambda$--edge-connected to $t$ for some $\lambda \in I_{\kappa}$ with $\lambda \geq \lambda_0$.
    As in the previous proofs, recall that $t = t_D$ for some component $D$ of $G - G_i$ and that  $\cardinality{A_t} \leq \sum_{s < t} \cardinality{V_s} \leq \cardinality{i} + \theta(\kappa)$ by item \cref{item:NPTnarrow}.

    Indeed, since $t$ has height $i < \kappa$ due to \cref{item:NPTheight}, let $\lambda\in I_{\kappa}$ be any cardinal for which $\lambda > \cardinality{i} + \theta(\kappa)$ and $\lambda \geq \lambda_0$.
    We recall from the definition of $Y_t \subseteq V_t$ as in \cref{item:VstarSet} that this set contains a finite sequence $S(u_{t},v_*,\lambda)$ with the properties claimed by \cref{lemma:AuxiliaryPaths}. 
    In particular, since $u_t \in V_t$ but $v_* \in V_r$, there are consecutive vertices $x,y \in S(u_t,v_*,\lambda)$ with vertex-connectivity at least $\lambda$ but $x \in D$ and $y \in V(G) - D$. 
    Moreover, there must exist a node $s < t$ such that $y \in V_s$: otherwise, $A_t = \bigcup_{s < t} V_s$ defines a set of fewer than $\lambda$ vertices such that $x$ and $y$ belong to distinct components of $G - A_t$, by the normality of $(T_i,\mathcal{V}_i)$ in $G$ (see item \cref{item:NPTcomponent} in the definition of $t=t_D$).
    This contradicts the $\lambda$--vertex-connectivity of $x$ and $y$ in $G$.

    We shall now prove that $\cardinality{\Dot{N}(t)\cap \Dot{N}(s)}\geq \lambda$ for the node $s<t$ with $y \in V_s$ just found.
    In particular, this will conclude that $s$ and $t$ are connected via at least $\lambda$ internally disjoint paths of length two in $\Dot{G}$, and hence that $s$ and $t$ are $\lambda$--vertex-connected.

    \begin{figure}[h!]
    \centering
   \begin{tikzpicture}
 
  \draw[fill=okvermillion] (0.4,5) circle (0.4cm);
  \node[shape=circular sector, minimum height=161pt, rotate=-90, fill=black!61, fill opacity=0.2, circular sector angle=116, minimum width=199pt] at (2,9.7) {};
  \draw[fill=okvermillion, fill opacity=1] (2,6.8) circle (0.4cm);
  \draw[fill=okvermillion, fill opacity=0.5] (1.5,7.5) circle (0.2cm);
  \draw (1,7.5) circle (0.2cm);
  \draw[fill=okvermillion, fill opacity=0.5] (2,7.5) circle (0.2cm);
  \draw[fill=okvermillion, fill opacity=0.5] (2.5,7.5) circle (0.2cm);
  \draw[fill=okvermillion, fill opacity=0.5] (3,7.5) circle (0.2cm);
  \draw[fill=okvermillion, fill opacity=0.5] (1.5,8) circle (0.2cm);
  \draw[fill=okvermillion, fill opacity=0.5] (1,8) circle (0.2cm);
  \draw[fill=okvermillion, fill opacity=0.5] (2,8) circle (0.2cm);
  \draw[fill=okvermillion, fill opacity=0.5] (2.5,8) circle (0.2cm);
  \draw[fill=okvermillion, fill opacity=0.5] (3,8) circle (0.2cm);
  \draw[fill=okvermillion, fill opacity=0.5] (3.5,8) circle (0.2cm);
  \draw (0.5,8) circle (0.2cm);
  \draw (0,8) circle (0.2cm);
  \draw (4,8) circle (0.2cm);
  \draw[fill=okvermillion, fill opacity=0.5] (1.5,8.5) circle (0.2cm);
  \draw[fill=okvermillion, fill opacity=0.5] (1,8.5) circle (0.2cm);
  \draw[fill=okvermillion, fill opacity=0.5] (2,8.5) circle (0.2cm);
  \draw[fill=okvermillion, fill opacity=0.5] (2.5,8.5) circle (0.2cm);
  \draw[fill=okvermillion, fill opacity=0.5] (3,8.5) circle (0.2cm);
  \draw[fill=okvermillion, fill opacity=0.5] (3.5,8.5) circle (0.2cm);
  \draw[fill=okvermillion, fill opacity=0.5] (0.5,8.5) circle (0.2cm);
  \draw (0,8.5) circle (0.2cm);
  \draw[fill=okvermillion, fill opacity=0.5] (4,8.5) circle (0.2cm);
  \draw (4.5,8.5) circle (0.2cm);
  \draw (-0.5,8.5) circle (0.2cm);
  \draw[fill=okvermillion, fill opacity=0.5] (1.5,9) circle (0.2cm);
  \draw (1,9) circle (0.2cm);
  \draw[fill=okvermillion, fill opacity=0.5] (2,9) circle (0.2cm);
  \draw[fill=okvermillion, fill opacity=0.5] (1,9) circle (0.2cm);
  \draw[fill=okvermillion, fill opacity=0.5] (0.5,9) circle (0.2cm);
  \draw[fill=okvermillion, fill opacity=0.5] (0,9) circle (0.2cm);
  \draw (-0.5,9) circle (0.2cm);
  \draw (-1,9) circle (0.2cm);
  \draw (1.5,9.5) circle (0.2cm);
  \draw[fill=okvermillion, fill opacity=0.5] (1,9.5) circle (0.2cm);
  \draw[fill=okvermillion, fill opacity=0.5] (0.5,9.5) circle (0.2cm);
  \draw[fill=okvermillion, fill opacity=0.5] (0,9.5) circle (0.2cm);
  \draw (-0.5,9.5) circle (0.2cm);
  \draw (-1,9.5) circle (0.2cm);
  \draw[fill=okvermillion, fill opacity=0.5] (4,9) circle (0.2cm);
  \draw[fill=okvermillion, fill opacity=0.5] (4.5,9) circle (0.2cm);
  \draw (5,9) circle (0.2cm);
  \draw (4.5,9.5) circle (0.2cm);
  \draw (5,9.5) circle (0.2cm);
  \draw (5.5,9.5) circle (0.2cm);
  \draw (5,9.5) circle (0.2cm);
  \node[shape=circular sector, minimum height=77pt, rotate=-90, circular sector angle=89, fill=black!61, fill opacity=0.4] at (3,10.65) {};
  \draw[fill=black] (0.28,5.1) circle (0.05cm);
  \draw[fill=black] (2.15,6.95) circle (0.05cm);
  \draw[semithick] (0.28,5.1) .. controls (-0.5,6.5) and (-1.5,7.3) .. (1.8,10.5);
  \draw[semithick] (2.15,6.95) .. controls (4.7,6.95) and (5.3,7.5) .. (4,10.5);
  \draw[decorate, decoration=snake, /pgf/decoration/amplitude=2pt, semithick] (1.8,10.5) .. controls (2.47,10.5) and (3.13,10.5) .. (4,10.5);
  \draw[dashed, semithick] (0.38,5.4) .. controls (-0.2,6.77) and (0.6,7.7) .. (2.6,9.3);
  \draw[semithick, dashed] (2.39,6.8) .. controls (2.99,6.8) and (3.7,9) .. (3.39,9.2);
  \draw[fill=okvermillion, fill opacity=1] (3,9.3) circle (0.4cm);
  \node[draw=none, minimum width=20pt, node font=\Large, text=okvermillion] at (2.3,6.25) {$t$};
  \node[draw=none, minimum width=20pt, node font=\Large, text=okvermillion] at (3.6,9.9) {$t_{D'}$};
  \node[draw=none, minimum width=20pt, node font=\Large, text=okvermillion] at (0.97,4.71) {$s$};
  \node[draw=none] at (0.42,4.89) {$y$};
  \node[draw=none] at (2,6.77) {$x$};
  \node[draw=none, minimum width=24pt, node font=\large] at (3.77,11.16) {$D'$};
\end{tikzpicture}

        \caption{The subtree $T'$ from the proof of \cref{lemma:kEdgeConnected} is displayed in red. The high vertex-connectivity between the vertices $x$ and $y$ implies that $t_{D'}\in T'$.}
        \label{fig:kEdgeConnected}
\end{figure}
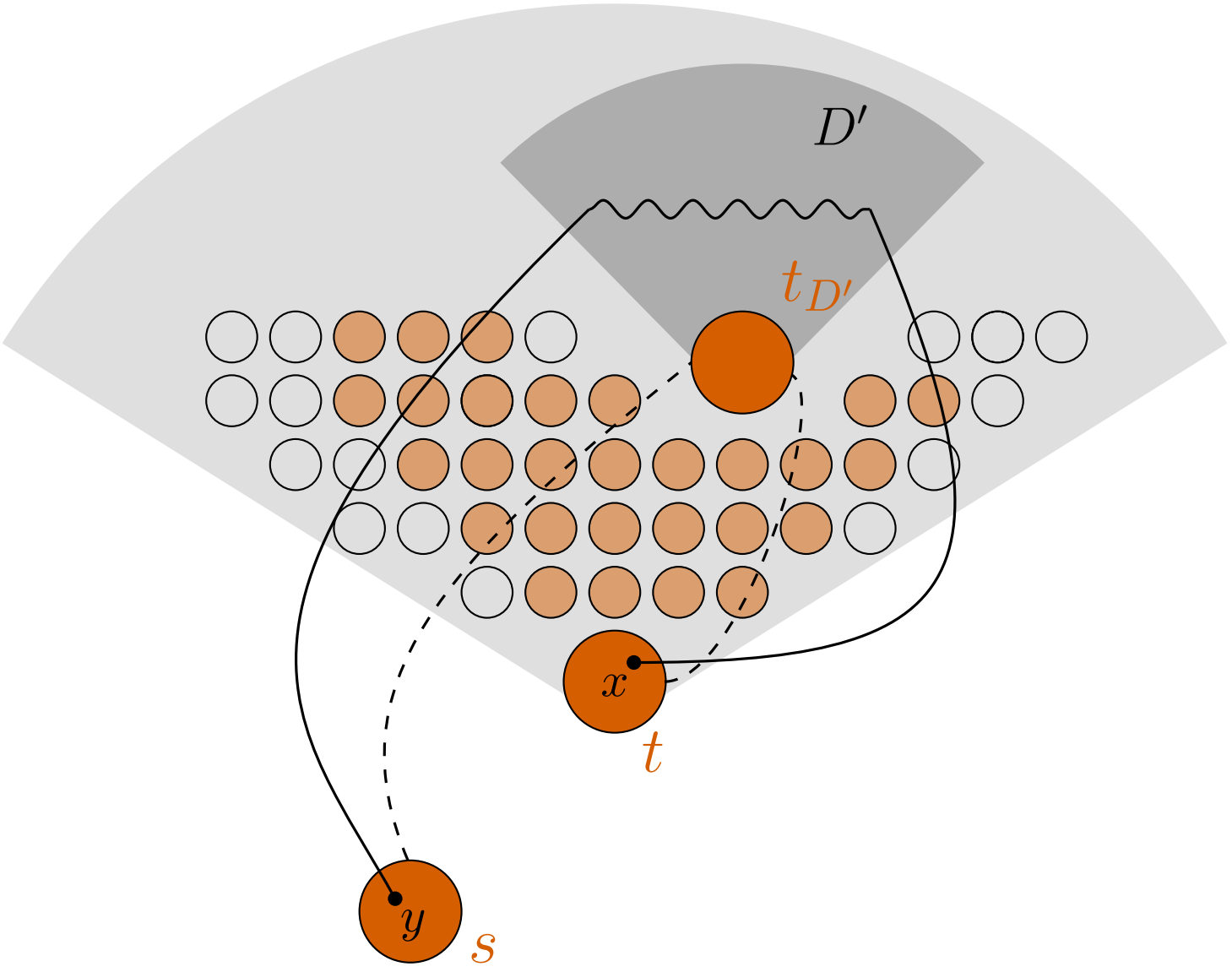

    In fact, suppose for a contradiction that the subtree $T' \coloneq (\Dot{N}(t)\cap \Dot{N}(s))\cup \downcl{t}$ has size $<\lambda$, hence comprising only nodes of height $<\lambda$ due to the down-closed neighborhood property of $\Dot{G}$ (and because $\mathrm{ht}(s)<\mathrm{ht}(t)< \lambda$). 
    By \cref{item:NPTnarrow} it follows that $\cardinality{V_{s'}}\leq \mathrm{ht}(s')+\theta(\kappa)< \lambda$ for every $s'\in T'$, and the regularity of $\lambda$ then ensures that $X' \coloneq \bigcup_{s'\in T'}V_{s'}$ also has size $< \lambda$. 
    Therefore, the $\lambda$--vertex-connectivity in $G$ between the vertices $x$ and $y$ claims the existence of some component $D'$ of $G - X'$ containing neighbors of both $x$ and $y$.
    Due to $\Dot{G}$ being a $T$--graph and $T'$ being a down-closed subset of $T$, \cref{lemma:PropertiesTGraphs} provides a branch $C'$ of $T'$ such that $\Dot{N}(D')\subseteq C'$ or, equivalently, $N(D')\subseteq \bigcup_{s'\in C'}V_{s'}$.
    In particular, there is some $j < \kappa$ (namely, the order type of $C'$) such that $C'$ is a branch of $T_j$, and the normality of $(T_j,\mathcal{V}_j)$ in $G$ then ensures that $D'$ is a component of $G - G_{j}$, in which both $x$ and $y$ have neighbors.
    Hence, $s,t \in \Dot{N}(D')$ and, thus, there are vertices $v_{s}^{D'}, v_t^{D'} \in V_{t_{D'}}$ which have neighbors in $V_s$ and $V_t$, respectively, due to the definition of the set $Z_{t_{D'}} \subseteq V_{t_{D'}}$.
    Then, as suggested by \cref{fig:kEdgeConnected}, $t_{D'}$ provides a node in the contraction minor $\Dot{G}$ that is a neighbor of both $t$ and $s$.
    Moreover, by the definition of $T'$ it follows that $t_{D'} \in T'$.
    This contradicts the fact that $D' \cap X' = \emptyset$ (as $V_{t_{D'}} \subseteq D'$), thus concluding that $\cardinality{\Dot{N}(t) \cap \Dot{N}(s)} = \cardinality{(\Dot{N}(t) \cap \Dot{N}(s)) \cup \downcl{t}} \geq \lambda$ (because $\cardinality{\odowncl{t}} = \cardinality{\mathrm{ht}(t)} < \lambda$).
\end{proof}

\section{Constructing the interleaved graph}\label{sec:sufficientCondition}

Let $\kappa$ be some fixed uncountable cardinal.
As suggested in \cref{sec:OrganizingPaper}, the rank functions introduced in \cref{sec:RankFunctions} provide a first tool for embedding $\TkStar$ in a $\kappa$--edge-connected graph $G$. 
Indeed, if $A\subseteq V(G)$ is a connected set of size $<\kappa$ for which some component of $G-A$ admits no $N(A)$--rank, then $A$ can be defined as the branch set of the apex vertex in a minor embedding of $\TkStar$ in $G$, where the branch sets corresponding to the tree $T_{\kappa}$ are now given by the second bullet point of \cref{lemma:Arank}. 
Of course, the connectedness assumption on $A$ can be dropped here: after all, this set could be replaced by a connected extension $A'\supseteq A$ of size $\cardinality{A'} \leq \cardinality{A}+\aleph_0$, since any $N(A)$--rankless component of $G-A$ must then contain an $N(A')$--rankless component of $G-A'$.

Unfortunately, as shown by \cref{example:LackOfSpreadingPair}, this sufficient condition is not necessary.
The goal of this section is then to deduce that, if $\kappa$ is a regular cardinal, then $\TkM$ is found as a minor of every $\kappa$--edge-connected graph to which the above approach fails for embedding $\TkStar$. 
In this sense, $\TkM$ is essentially the minor-minimal obstruction as to why \cref{lem:reductionToTgraphs} cannot be stated only in terms of its first bullet point.

Indeed, when relying on this technical result, any $\kappa$--edge-connected graph in which this sufficient condition is not satisfied admits a normal partition tree $(T,\mathcal{V})$ whose arising $T$--graph minor, which we also call $G$, is again $\kappa$--edge-connected and satisfies the down-closed and short neighborhood properties.
Such a graph, as well as the corresponding normal order tree $(T,\leq)$, will then be fixed from now on in this section, where we aim to detail a further reduction of \cref{thm:mainThm}.

First, let us observe the following consequence of the down-closed neighborhood property: given two nodes $t_0<t_1 \in T$, the set $N([t_0,t_1)) \cap \upcl{t_1}$ is down-closed in $\upcl{t_1}$ and, hence, a subtree of $T$.
More generally, if $B \subseteq V(T)$ is any set and $t \in T$ a node for which $B$ is disjoint from $\upcl{t}$, the set $\upcl{t}_{N(B)}$ is again a convex subset of $\upcl{t}$ and, if non-empty, contains $t$. 
Then, when $B$ is connected and disjoint from $\upcl{t}$, we call the tuple $(B,t)$ a \defn{spreading pair} if $\upcl{t}_{N(B)}$ is of size at least $\kappa$. 
Due to the normality of $T$ in $G$, we highlight that $\upcl{t}_{N(B)} = \upcl{t}_{N(B\cap \downcl{t})}$.
In the case in which $G$ is $\kappa$--edge-connected, the cut \defn{$\partial \upcl{t}$}${} \coloneq \set{uv\in E(G)}:{u\geq t, v\notin \upcl{t}}$ has size at least $\kappa$, if $t$ is not the root of $T$, and, if $\mathrm{ht}(t) < \kappa$, then $\cardinality{N(\odowncl{t})\cap \upcl{t}}\geq \kappa$. 
When $(B,t)$ is a spreading pair satisfying $\cardinality{N(\odowncl{t}\setminus B)\cap \upcl{t}}<\kappa$, meaning that all but fewer than $\kappa$ edges from $\partial \upcl{t}$ have an endpoint in $B$, we then say that $(B,t)$ is a \defn{concentrated spreading pair}.

In particular, it is an immediate consequence of the short neighborhood property that, when $\kappa$ is a regular cardinal, the antichain claimed by \cref{prop:shortNeighborhoodProp} must have size $\kappa$. 
For the reader interested in the combinatorics of uncountable order trees, the paragraph below argues that the existence of large antichains in $\upcl{t}_{N(B)}$ is essentially the only explanation for this set to have size $\geq \kappa$.

\vspace{0.5cm}

\noindent {\huge $\ulcorner$}

The following paragraph, up to the end of \cref{obs:hasARank}, is not needed to follow the paper, but might be interesting to readers curious about the combinatorics of uncountable order trees, and of $\upcl{t}_{N(B)}$ in particular.
Now, the existence of large antichains in  $\upcl{t}_{N(B)}$ does not follow just from the order-theoretical properties of this tree:
after all, $\upcl{t}_{N(B)}$ could be a complex order tree in the following sense. It could either contain a $\kappa$--sized antichain, a $\kappa$--sized chain or none of those, in which case it is a \textit{$\kappa$--Suslin tree}.
Curiously, the short neighborhood property has the further advantage of excluding the latter two cases. 
First, it is immediate from the particular antichain claimed by \cref{prop:shortNeighborhoodProp} that $\upcl{t}_{N(B)}$ cannot contain an interval (and hence no chain) of size $\kappa$.
On the other hand, while the existence of Suslin trees is independent of the axioms of $\mathrm{ZFC}$, when they exist, they can be pruned to a convex subtree with the property that all their up-closures are again Suslin trees (see Lemma III.5.27 in \cite{KennethKunen}).
In particular, $\upcl{t}_{N(B)}$ has no $\kappa$--rank: otherwise, choose $s \in \upcl{t}_{N(B)}$ such that $\upcl{s}_{N(B)}$ contains a Suslin tree and has minimal rank among all such $s$.
But for any rank reducing set $X$ of $\upcl{s}_{N(B)}$ there exists some $s' \in \upcl{s}_{N(B)}$ such that $\upcl{s'}_{N(B)}$ is disjoint from $X$ and again contains an order-theoretic Suslin tree, contradicting the choice of $s$.
However, the following application of \cref{prop:shortNeighborhoodProp} shows that $\upcl{t}_{N(B)}$ does have a $\kappa$--rank:

\begin{observation}\label{obs:hasARank}
    If $(B,t)$ is a spreading pair, then $\upcl{t}_{N(B)}$ has a $\kappa$--rank.
\end{observation}

\begin{proof}
    Let us assume for a contradiction that $\upcl{t}_{N(B)}$ induces a subgraph of $G$ with no $\kappa$--rank.
    In particular, there exists a $T_\kappa$ minor embedding $\set{V_s}:{s \in T_\kappa}$ in $\upcl{t}_{N(B)}$ by \cref{thm:kappaRank}.
    Then let $v_s \coloneq \min V_s$ for every $s\in T_{\kappa}$ and define $B' \coloneq \bigcup_{s \in T_\kappa} [t,v_s)$. 
    On the one hand, for each $s \in T_{\kappa}$, the set $V_s$ is joined via $\kappa$ internally disjoint paths to the nodes from $\set{v_{t'}}:{t' \in \mathrm{succ}_{\preceq}(s)}$, where $\preceq$ denotes here the graph-theoretical order of $T_{\kappa}$. 
    On the other hand, together with the minimality of $v_s$, and the fact that its height in $T$ is $<\kappa$, the normality of $(T,\mathcal{V})$ implies that $\kappa$ many of the nodes of $\set{v_{t'}}:{t' \in \mathrm{succ}_{\preceq}(s)}$ belong to $\upcl{v_s}$.  
    All in all, $\set{v_s}:{s\in T_{\kappa}}\subseteq B'$ and, thus, $B'$ is a $\kappa$--sized convex subset of $\upcl{t}_{N(B)} = \upcl{t}_{N(B\cap \odowncl{t})}\subseteq\upcl{t}_{N(\odowncl{t})} $ by the down-closed neighborhood property of $G$.

    However, due to the short neighborhood property of this $T$--graph, there exists a node $x \in B'$ whose up-closure $\upcl{x}_{B'}$ is of size $< \kappa$.
    But this contradicts the definition of $B'$: given any $t \in T_\kappa$, only $< \kappa$ many nodes from $\set{v_{t'}}:{t' \in \mathrm{succ}_{\preceq}(t)}$ are contained in $T \setminus \upcl{v_t}$.
\end{proof}

\hfill {\huge $\lrcorner$}

\vspace{0.5cm}

To summarize, the order tree $\upcl{t}_{N(B)}$ must contain an antichain of size $\kappa$ if $\kappa$ is a regular cardinal.
In fact, for every convex set $\Delta \subseteq \upcl{t}_{N(B)}$ of size $\geq \kappa$, \cref{prop:shortNeighborhoodProp} can be restated as

\begin{enumerate}[label=$(\diamond_{2.1})$]
    \item\label{prop:shortNeighborhoodProp:case1} there is an antichain $A \subseteq \Delta$ of size $\kappa$ with $\cardinality{\upcl{a}_{\Delta}} < \kappa$ for every $a \in A$.
\end{enumerate}

If $\kappa$ is singular and $\Set{\kappa_{\alpha}}_{\alpha < \cof{\kappa}}$ is a cofinal sequence comprising regular cardinals, then the short neighborhood property \cref{prop:shortNeighborhoodProp} splits into two cases: either the above item \cref{prop:shortNeighborhoodProp:case1} holds or, after applying \cref{claim:antichainKappa} to cofinally many up-closures from $\set{\upcl{a}_{\Delta}}:{a\in A}$,

\begin{enumerate}[label=$(\diamond_{2.2})$]
    \item\label{prop:shortNeighborhoodProp:case2} there is a $\cof{\kappa}$--sized antichain $A = \Set{a_\alpha}_{\alpha < \cof{\kappa}} \subseteq \Delta$ such that $\upcl{a_\alpha}_{\Delta}$ has cardinality $< \kappa$ but height at least $\kappa^+_{\alpha}$.
\end{enumerate}

Inspired by this discussion, given a cardinal $\lambda$ and a connected set $B \subseteq T$, we say that $A \subseteq N(B)$ is a \defn{$(B,\lambda)$--strong antichain} if

\begin{enumerate}
    \item $A$ is an antichain of size $\lambda$;
    
    \item\label{item:StrongAntichainCondition} $\upcl{a}$ is disjoint from $B$ and $\upcl{a}_{N(B)}$ is of size $<\kappa$ for every $a \in A$.
\end{enumerate}

Lastly, we say that a spreading pair $(B,t)$ is \defn{extendable} if there exists a connected set $B_* \subseteq V(T)$ with $B \subseteq B_* \subseteq B \cup \oupcl{t}$ and for which there exists a $(B_*,\kappa)$--strong antichain $A \subseteq \oupcl{t}$ such that $\left( \bigcup_{a\in A} \odowncl{a} \setminus B_* \right) \cap \upcl{t}$ is connected (see also \cref{fig:extension}).
Incidentally, this latter set must contain the node $t$.

The next theorem concludes that if every concentrated spreading pair is extendable, one finds a minor embedding of $\TkM$ in $G$, thus obtaining $\TkStar$ as a byproduct minor due to the discussion carried out after \cref{example:interleaved}.

More precisely, formalizing the second step from the approach outlined in \cref{sec:OrganizingPaper}, the statement below is the main result of this section and provides the central criterion for addressing the proof of \cref{thm:mainThm}.

\begin{theorem}\label{thm:sufficientCondition}
    Let $T$ be an order tree of height at most $\kappa$ and $G$ be a $\kappa$--edge-connected $T$--graph that satisfies the down-closed neighborhood property.
    If every concentrated spreading pair is extendable, then $G$ contains $\TkM$ as a minor.
\end{theorem}

\begin{proof}
    Let $t_0$ denote the root of $T$ and note that there is an antichain $A_0 \subseteq \oupcl{t_0}$ of size $\kappa$: after all, due to the $\kappa$--edge-connectivity of $G$, the tuple $(\Set{t_0},s)$ provides a concentrated spreading pair for any given $s \in \successor{t_0}$, thus admitting a corresponding strong antichain by the extendability hypothesis.
    Set then $T_0 \coloneq \Set{t_0}$, $X_0 \coloneq \emptyset$ and $B_0 \coloneq \bigcup_{a\in A_0} \odowncl{a}$; in particular, $T_1\coloneq \Set{t_0} \cup A_0$ is order-isomorphic to the first level of $T_{\kappa}$ endowed with its usual graph-theoretical ordering $\preceq$.
    By the normality of $T$ and the $\kappa$--edge-connectivity of $G$, it follows that $(B_0,a)$ is a concentrated spreading pair for every $a \in A_0$, since $\mathrm{ht}(a) < \kappa$.

    By induction on $n\in\mathbb{N}$, we will now describe an order-isomorphic copy $T_n$ in $T$ of the tree $T_{\kappa}^n$ comprising the nodes of height at most $n$ in $T_{\kappa}$, whose standard graph-theoretical order tree we will denote by $\preceq$.
    This is done so that $T_{n-1} \subseteq T_{n}$ for $n \geq 1$. 
    Moreover, we refer to the predecessor of $t \in T_n$ as $t^-$.
    
    In particular, for every $t \in T_{n-1}$ we assign two sets $B_t, X_t \subseteq \upcl{t}$ with the properties below (see also \cref{fig:interleavedConstruction}):

    \begin{enumerate}
        \item $\Set{B_{t},X_{t}}$ is a disjoint cover of $\bigcup_{s \in \mathrm{succ}_{\preceq}(t)}[t,s)$;
        
        \item\label{item:BagIndex} $B_t$ is connected and $t$ is the $\leq$--minimal element of $B_t$;
        
        \item\label{item:ConnectedPiece} $X_t \cup B_{t^-}$ is connected;

        \item\label{item:connectedToV--} if $t \notin T_1$ then $t$ neighbors some node in $X_{t^-} \cup B_{t^{--}}$;
        
        \item\label{item:ConcentratedPair} $(B_{t},s)$ is a concentrated spreading pair for every $s \in \mathrm{succ}_{\preceq}(t)$.
    \end{enumerate}

    \begin{figure}[ht]
        \centering
        \begin{minipage}{0.5\textwidth}
            \centering
            \begin{tikzpicture}[scale=2]

                \tikzset{
                    vertex/.style={circle, fill, inner sep=1.7pt},
                    coverV/.style={circle, fill, inner sep=1.7pt, color=okorange},
                    edge/.style={draw=black, line width=1pt},
                    smallEdge/.style={draw=black, line width=0.8pt},
                    every picture/.append style={transform shape},
                    directed/.style={->,edge},
                    big_cover/.style={line width=5pt, color=okorange, cap=round, opacity=0.4},
                    med_cover/.style={line width=3pt, cap=round},
                    hidden_edge/.style={draw=black, line width=0.4pt, dashed},
                }

                \node[vertex] (r) at (0,-1) {};
            
                \coordinate (b_start) at (0,-0.8);
                \coordinate (b_Smid) at (0,-0.7);
            
                \coordinate (b_1) at ($(b_start) !0.5! (b_Smid)$);

                \coordinate (b_Emid) at (0,-0.4);
                \coordinate (b_end) at (0,-0.2);
            
                \coordinate (b_2) at ($(b_Emid) !0.5! (b_end)$);

                \node[anchor=east, color=okorange] at (0,-0.5) {$B$};
            
                \coordinate (lMiddle) at (120:1.5);
                \coordinate (rMiddle) at (60:1.5);
            
                \coordinate(lShiftMiddle) at ($(lMiddle) + (-0.2,0)$);
                \coordinate(rShiftMiddle) at ($(rMiddle) + (0.2,0)$);

                \fill[color=oktintblack, opacity=0.8, on background layer] (0,0) -- ($(0,0) !0.98! (lShiftMiddle)$) -- ($(0,0) !0.98! (rShiftMiddle)$);
                \node[vertex] (t) at (0,0) {};
            
                \draw[edge] (t) -- (lShiftMiddle);
                \draw[edge] (t) -- (rShiftMiddle);
                
                \coordinate (cf1) at ($(lMiddle)!0!(rMiddle)$);
                
                \begin{scope}[shift=(cf1)]
                    \coordinate (cf1LTrEnd) at (120:1);
                    \coordinate (cf1RTrEnd) at (60:1);
            
                    \fill[color=oktintblack, opacity=0.8, on background layer] (cf1.center) -- ($(cf1.center) !0.95! (cf1LTrEnd)$) -- ($(cf1.center) !0.95! (cf1RTrEnd)$);
            
                    \node[vertex, color=okgreen] (cf1E1) at ($(cf1LTrEnd) !0.2! (cf1RTrEnd)$) {};
                    \node[vertex, color=okgreen] (cf1E2) at ($(cf1LTrEnd) !0.5! (cf1RTrEnd)$) {};
                    \node[vertex, color=okgreen] (cf1E3) at ($(cf1LTrEnd) !0.8! (cf1RTrEnd)$) {};
            
                    \node[vertex, color=okorange] (cf1Et) at ($(cf1E2) !0.5! (cf1)$) {};
                    \draw[smallEdge,dash pattern=on 1pt off 1pt on 1pt off 1pt] (cf1E1) to[bend right=25] (cf1Et);
                    \draw[smallEdge,dash pattern=on 1pt off 1pt on 1pt off 1pt] (cf1E2) to[bend left=25] (cf1Et);
                    \draw[smallEdge,dash pattern=on 1pt off 1pt on 1pt off 1pt] (cf1E3) to[bend left=25] (cf1Et);

                    \coordinate (a11Prime) at ($(cf1E1) + (-0.1,0)$);
                    \coordinate (a31Prime) at ($(cf1E3) + (0.1,0)$);  
                    \draw[med_cover, opacity=0.3, color=okgreen] (a11Prime) -- (a31Prime);

                    \node[vertex, color=okskyblue] (cf1) at ($(lMiddle)!0!(rMiddle)$) {};
                \end{scope}

                \node[vertex, color=okskyblue] (cf2) at ($(lMiddle)!0.5!(rMiddle)$) {};
                
                \coordinate (cf3) at ($(lMiddle)!0.95!(rMiddle)$);
            
                \draw[med_cover, opacity=0.3, color=okskyblue] ($(lMiddle) + (-0.2,0)$) -- ($(rMiddle) + (0.2,0)$);

                \begin{scope}[shift=(cf3)]
                    \coordinate (cfLTrEnd) at (120:1.5);
                    \coordinate (cfRTrEnd) at (60:1.5);
            
                     \fill[color=oktintblack, opacity=0.8, on background layer] (cf3.center) -- ($(cf3) !0.25! (cfLTrEnd)$) -- ($(cfRTrEnd) !0.25! (cfLTrEnd)$) -- (cfRTrEnd);
                     \node[vertex, color=okskyblue] (cf3) at ($(lMiddle)!0.95!(rMiddle)$) {};
                
                    \coordinate (ext_start) at (0,0.3);
                    \coordinate (ext_end_1) at ($(cfRTrEnd) + (ext_start)$);
                    \coordinate (ext_end) at ($(ext_start) !0.7! (ext_end_1)$);
            
                    \node[coverV] (a0_B) at (ext_start) {};
                    
                    \node[coverV] (a1_B) at ($(ext_start) !0.5! (ext_end)$) {};
                    \node[vertex, color=okgreen] (a1) at ($(cfLTrEnd) !0.2! (cfRTrEnd)$) {};
            
                    \node[coverV] (a2_B) at ($(ext_start) !0.8! (ext_end)$) {};
                    \node[vertex, color=okgreen] (a2) at ($(cfLTrEnd) !0.5! (cfRTrEnd)$) {};

                    \node[vertex, color=okgreen] (a3) at ($(cfLTrEnd) !0.7! (cfRTrEnd)$) {};
            
                    \coordinate (a1Prime) at ($(a1) + (-0.1,0)$);
                    \coordinate (a3Prime) at ($(a3) + (0.1,0)$);  
                    \draw[med_cover, opacity=0.3, color=okgreen] (a1Prime) -- (a3Prime);

                    \draw[smallEdge,dash pattern=on 1pt off 1pt on 1pt off 1pt] (a0_B) to[bend left=25] (a1);
                    \draw[smallEdge,dash pattern=on 1pt off 1pt on 1pt off 1pt] (a1_B) to[bend left=25] (a2);
                    \draw[smallEdge,dash pattern=on 1pt off 1pt on 1pt off 1pt] (a2_B) to[bend left=25] (a3);
                    \draw[smallEdge] (a0_B) -- (a1_B) -- (a2_B);
                    \draw[smallEdge, ->] (a2_B) -- (ext_end);
                \end{scope}

                \draw[edge] (cf1) -- (cf1LTrEnd);
                \draw[edge] (cf1) -- (cf1RTrEnd);

                \draw[edge] (cf3) -- (cfLTrEnd);
                \draw[edge] (cf3) -- (cfRTrEnd);

                \draw[edge] (0,0) to[in=-90, out=90] (cf3);
            
                \path (0,0) to[in=-90, out=90] coordinate[pos=0.5] (intersec1) (cf3);

                \draw[edge, dashed, color=okorange] (a2_B) to[bend left=50] (a1_B);
                \draw[edge, dashed, color=okorange] (a2_B) to [bend left=60] (a0_B);
                \draw[edge, dashed, color=okorange] (a0_B) to[bend left=60] (intersec1);

                \draw[edge, color=okorange] (b_1) to[bend right=35] (b_2);
                \draw[edge, dashed,  color=okorange] (b_2) to[bend right=35] (intersec1);
                \node[vertex, color=okorange] at (intersec1) {};

                \draw[edge] (0,0) to[in=-80, out=90] (cf1);
                \path (0,0) to[in=-90, out=90] coordinate[pos=0.3] (intersec2) (cf1);

                \draw[edge, dashed, color=okorange] (intersec2) to[bend left=70] (cf1Et); 
                \draw[edge, dashed, color=okorange] (b_2) to[bend left=35] (intersec2); 
                
                \node[vertex, color=okorange] at (intersec2) {};
                
                \draw[edge] (r) -- (t);
            
                \node[anchor=east] at (t) {$t$};
                \node[anchor=east] at (r) {$r$};
            
                 \node[anchor=west, color=okgreen] at (cfRTrEnd) {$A$};

                \node[anchor=west] at ($(0.2,0) !0.4! (rShiftMiddle)$) {$\upcl{t}_{N(B)}$};
            
                \draw[med_cover, color=okorange] (b_start) -- (b_Smid);
                \draw[med_cover, color=okorange] (b_Emid) -- (b_end);
                 
            \end{tikzpicture}
            \caption{The set $B$ can be extended to a connected superset $B_*$ via the dashed orange edges.
            The gray shaded area corresponds to $\upcl{t}_{N(B)}$.}\label{fig:extension}
        \end{minipage}\hfill
        \begin{minipage}{0.5\textwidth}
            \centering
                \begin{tikzpicture}[scale=0.7]
                    \tikzset{
                        vertex/.style={circle, fill, inner sep=1.7pt},
                        edge/.style={draw=black, line width=1pt},
                        every picture/.append style={transform shape},
                    }
                
                    \tikzmath{
                        \Xset = 2.5; \Bset = 5; \circleDia = 0.15; 
                    }

                    \coordinate (vertexHeight) at (0,0.5);
                    \coordinate (vertexShift) at (0.5,0);

                    \coordinate (ranVertexHeight) at (0,0.3);
                    \coordinate (ranVertexShift) at (0.2,0);

                    \node[vertex,okpurple] (t-) at (0,0) {};
                    \node (t-lh) at (120:\Xset) {};
                    \node (t-rh) at (60:\Xset) {}; 
                    \node (t-l) at (120:\Bset) {};
                    \node (t-r) at (60:\Bset) {};
                
                    \node (t-circleP1) at (120:\circleDia * 1.2) {};
                    \node (t-circleP2) at (60:\circleDia * 1.2) {};
                
                    \filldraw[fill=okgreen] (t-lh.center) -- (t-rh.center) -- (t-circleP2.center) arc[start angle=60, end angle=120, radius=\circleDia * 1.2] (t-circleP1.center) -- (t-lh.center);
                        
                    \filldraw[fill=okpurple] (t-lh.center) -- (t-l.center) -- (t-r.center) -- (t-rh.center) -- (t-lh.center);

                    \node[vertex,okskyblue] (t) at ($(t-l) + (vertexHeight) + (vertexShift)$) {};
                    
                    \begin{scope}[shift=(t)]
                        \node (tlh) at (120:\Xset) {};
                        \node (trh) at (60:\Xset) {};
                        \node (tl) at (120:\Bset) {};
                        \node (tr) at (60:\Bset) {};

                        \node (tcircleP1) at (120:\circleDia*1.2) {};
                        \node (tcircleP2) at (60:\circleDia*1.2) {};
                
                        \begin{scope}[on background layer]
                            \filldraw[fill=okpurple] (tlh.center) -- (trh.center) -- (tcircleP2.center) arc[start angle=60, end angle=120, radius=\circleDia*1.2] (tcircleP1.center) -- (tlh.center);
                        
                            \filldraw[fill=okskyblue] (tlh.center) -- (tl.center) -- (tr.center) -- (trh.center) -- (tlh.center);
                        \end{scope}
                    \end{scope}

                    \node[vertex] (tBedge) at ($(t-)!0.6!(t)$) {};
                    \node[vertex] (tXedge) at ($(t-)!0.3!(t)$) {};
                
                    \draw[edge] (t) to[bend right=25] (tBedge);
                    \draw[edge] (t) to[bend left=40] (tXedge);

                    \node[vertex] (t1) at ($(tl) + (vertexHeight) + (vertexShift)$) {};
                    \node[vertex] (t4) at ($(tr) + (vertexHeight) + -1*(vertexShift)$) {};
                    \node[vertex] (t2) at ($(t1)!0.33333!(t4)$) {};
                    \node[vertex] (t3) at ($(t1)!0.66666!(t4)$) {};

                    \node[vertex,okorange] (rant) at ($(t-r) + (vertexHeight) + -1*(vertexShift)$) {};
                
                    \node[vertex] (rantBedge) at ($(t-)!0.7!(rant)$) {};
                    \node[vertex] (rantXedge) at ($(t-)!0.4!(rant)$) {};
                
                    \draw[edge] (rant) to[bend right=20] (rantBedge);
                    \draw[edge] (rant) to[bend left=15] (rantXedge);
                    
                    \begin{scope}[scale=0.4]
                        
                        \begin{scope}[shift=(rant)]
                            \node (rantlh) at (120:\Xset) {};
                            \node (rantrh) at (60:\Xset) {};
                            \node (rantl) at (120:\Bset) {};
                            \node (rantr) at (60:\Bset) {};
                
                            \node (rantcircleP1) at (120:\circleDia*3) {};
                            \node (rantcircleP2) at (60:\circleDia*3) {};
                
                            \begin{scope}[on background layer]
                                \filldraw[fill=okpurple] (rantlh.center) -- (rantrh.center) -- (rantcircleP2.center) arc[start angle=60, end angle=120, radius=\circleDia*3] (rantcircleP1.center) -- (rantlh.center);
                        
                                \filldraw[fill=okorange] (rantlh.center) -- (rantl.center) -- (rantr.center) -- (rantrh.center) -- (rantlh.center);
                            \end{scope}
                        \end{scope}
                    \end{scope}

                    \begin{scope}[scale=0.3]
                        \node[vertex] (rant1) at ($(rantl) + (ranVertexHeight) + (ranVertexShift)$) {};
                        
                        \begin{scope}[shift=(rant1)]
                            \node (rant1lh) at (120:\Xset) {};
                            \node (rant1rh) at (60:\Xset) {};
                            \node (rant1l) at (120:\Bset) {};
                            \node (rant1r) at (60:\Bset) {};
                
                            \node (rant1circleP1) at (120:4*\circleDia) {};
                            \node (rant1circleP2) at (60:4*\circleDia) {};
                
                            \begin{scope}[on background layer]
                                \filldraw[fill=okorange] (rant1lh.center) -- (rant1rh.center) -- (rant1circleP2.center) arc[start angle=60, end angle=120, radius=4*\circleDia] (rant1circleP1.center) -- (rant1lh.center);
                
                                \filldraw[fill=oktintblack] (rant1lh.center) -- (rant1l.center) -- (rant1r.center) -- (rant1rh.center) -- (rant1lh.center);
                            \end{scope}
                        \end{scope}
                    \end{scope}

                    \begin{scope}[scale=0.3]
                        \node[vertex] (rant2) at ($(rantr) + (ranVertexHeight) + -1*(ranVertexShift)$) {};
                        
                        \begin{scope}[shift=(rant2)]
                            \node (rant2lh) at (120:\Xset) {};
                            \node (rant2rh) at (60:\Xset) {};
                            \node (rant2l) at (120:\Bset) {};
                            \node (rant2r) at (60:\Bset) {};
                
                            \node (rant2circleP1) at (120:\circleDia*4) {};
                            \node (rant2circleP2) at (60:\circleDia*4) {};
                
                            \begin{scope}[on background layer]
                                \filldraw[fill=okorange] (rant2lh.center) -- (rant2rh.center) -- (rant2circleP2.center) arc[start angle=60, end angle=120, radius=\circleDia*4] (rant2circleP1.center) -- (rant2lh.center);
                        
                                \filldraw[fill=oktintblack] (rant2lh.center) -- (rant2l.center) -- (rant2r.center) -- (rant2rh.center) -- (rant2lh.center);
                            \end{scope}
                        \end{scope} 
                    \end{scope}

                    \node[anchor=east] at ($(t-)!0.5!(t-lh)$) {$X_{t^-}$};
                    \node[anchor=east] at ($(t-lh)!0.5!(t-l)$) {$B_{t^-}$};
                
                    \node[anchor=east] at ($(t)!0.5!(tlh)$) {$X_{t}$};
                    \node[anchor=east] at ($(tlh)!0.5!(tl)$) {$B_{t}$};
                
                    \node[anchor=east] at (t-) {$t^-$};
                    \node[anchor=east] at (t) {$t$};

                    \node[anchor=east] at (t1) {$s$};
                    \node[anchor=west] at ($(t4)+(0.3,0)$) {$A$};
                    
                \end{tikzpicture}
        \caption{The sets shaded in the same color will define the bags corresponding to their minimal elements.
        In particular, the node $t$ which belongs to the set $B_t$ will witness the two edges of $V_t$ to its neighbors $V_{t^-}$ and $V_{t^{--}}$.}\label{fig:interleavedConstruction}
        \end{minipage}
    \end{figure}

    Indeed, at the end of this recursive process, the family $\set{V_t \coloneq B_t \cup \bigcup_{s \in \mathrm{succ}_{\preceq}(t)}}X_s:{t \in T_\kappa}$ will comprise connected sets by item \cref{item:ConnectedPiece} and provide the desired minor embedding of $\TkM$ in $G$. 
    After all, relying on the down-closed neighborhood property of $G$, an edge between $V_t$ and $V_{t^-}$ is ensured by items \cref{item:BagIndex} and \cref{item:ConcentratedPair}, while \cref{item:connectedToV--} highlights an edge between $V_t$ and $V_{t^{--}}$.

    Fix then $0 < n < \omega$, and let us assume we have constructed an order-isomorphic copy $T_n$ of $T^n_\kappa$ satisfying the above five conditions.
    Let $t \in T_n$ be a node at height $n$.
    By the inductive assumption \cref{item:ConcentratedPair} it follows that $(B_{t^-}, t)$ is a concentrated spreading pair, thus being extendable by assumption of the theorem.
    Then let $X_t \subseteq \oupcl{t}$ be a set such that $B_{t^-} \cup X_t$ is connected and let $A \subseteq \oupcl{t}$ be a $(B_{t^-} \cup X_t, \kappa)$--strong antichain for which $B_t \coloneq \bigcup_{a \in A} [t, a) \setminus X_t$ is connected.
    We extend $t$ by declaring the $\kappa$ many elements from $A$ as its successors in $T_{n+1}$: indeed, since $A \subseteq \oupcl{t}$, this is a compatible extension with respect to the ordering of $T$.

    Given $s \in A$, it now remains to verify \cref{item:ConcentratedPair}, i.e., that $(B_t, s)$ is a concentrated spreading pair. 
    Since $G$ is $\kappa$--edge-connected, this requires showing $\cardinality{N(\odowncl{s} \setminus B_t) \cap \upcl{s}} < \kappa$.
    We first highlight that $N(\odowncl{s} \setminus B_t) \cap \upcl{s} = \big( N(B_{t^-} \cup X_t) \cap \upcl{s} \big) \cup \big( N(\odowncl{t} \setminus B_{t^-}) \cap \upcl{s} \big)$.
    The first set has size $< \kappa$ since $A$ is a $(B_{t^-} \cup X_t, \kappa)$--strong antichain, and the second set has size $< \kappa$ since $(B_{t^-}, t)$ is concentrated.    
\end{proof}

Our main strategy for proving \cref{thm:mainThm} is now to find sufficient conditions under which a given spreading pair $(B,t)$ is extendable.
For regular cardinals, this is an immediate consequence of the short neighborhood property \cref{prop:shortNeighborhoodProp:case1}:

\begin{corollary}\label{cor:FinishRegularCase}
    If $\kappa$ is an uncountable regular cardinal, every $\kappa$--edge-connected graph $G$ contains $\TkStar$ as a minor.
\end{corollary}

\begin{proof}
    As previously outlined, it suffices to consider the case in which $G$ is a $T$--graph with the down-closed and short neighborhood properties for some suitable order tree $T$ of height at most $\kappa$.
    Due to \cref{thm:sufficientCondition}, it suffices to prove that any spreading pair $(B,t)$ in $G$ is extendable.
    Indeed, $\Delta\coloneq \upcl{t}_{N(B)}$ contains $t$, is convex due to the down-closed neighborhood property, and of size $\geq \kappa$ by definition of spreading pair.
    By the short neighborhood property \cref{prop:shortNeighborhoodProp:case1}, there exists then an antichain $A\subseteq \Delta$ of size $\kappa$ such that $\upcl{a}_\Delta$ is of size $<\kappa$ for every $a \in A$.
    In particular, $t\notin A$.
    Hence, $A$ is a $(B,\kappa)$--strong antichain and witnesses that $(B,t)$ is extendable, since $\left(\bigcup_{a\in A}\odowncl{a}\setminus B \right)\cap \upcl{t} = \bigcup_{a\in A}[t,a)$ is connected by \cref{lemma:PropertiesTGraphs}.
\end{proof}

In the remainder of this paper, we introduce the machinery for concluding \cref{thm:mainThm} when $\kappa$ is a singular cardinal.
In contrast with the regular case above, where the short neighborhood property readily ensured that any spreading pair is extendable, this can fail in the singular setting.
In particular, it might be that the fixed order tree $T$, even when satisfying the short neighborhood property, contains no antichains of size $\kappa$ at all:

\begin{example}[Revisited tree of cofinal cliques]
    Let $\set{\kappa_\alpha}:{\alpha < \cof{\kappa}} \subseteq \kappa$ be a strictly increasing sequence of cardinals converging towards $\kappa$.
    Similarly to the construction of the tree of cofinal cliques in \cref{example:TreeOfCliques}, one can define a $\kappa$--edge-connected supergraph $\hat{\mathbb{T}}_{\kappa}$ of $\mathbb{T}_{\kappa}$ by inflating $p \in \cof{\kappa}^{<\omega}$ by a clique of size $\kappa_{\ell(p)}$ where $\ell(p)$ is the last element of $p$ or $1$ if $p$ is the empty sequence, and then adding all edges of the form $\set{uv}:{u \in K_t, v \in K_s, \text{ where }t \sqsubseteq s, \mathrm{dist}_{\cof{\kappa}^{<\omega}}(t,s) \leq 2}$.

    However, if $\mathcal{V}\coloneq \set{V_t}:{t\in T}$ is any normal partition tree of $\hat{\mathbb{T}}_{\kappa}$, the index tree $T$ contains no antichain of size $\kappa$: after all, for any subset $A\subseteq T$ of cardinality greater than $\cof{\kappa}$, there must exist some $s \in \cof{\kappa}^{<\omega}$ such that $V_{t} \cap K_s \neq \emptyset \neq V_{t'} \cap K_s$ for two distinct $t,t' \in A$. 
    In this case, a $V_t$--$V_{t'}$ edge in $E(K_s) \subseteq E(\hat{\mathbb{T}}_{\kappa})$ ensures that $t$ and $t'$ are comparable in $T$ by normality of $(T,\mathcal{V})$. 
    
    Despite that, $\hat{\mathbb{T}}_{\kappa}$ has a $T$--graph structure that satisfies the down-closed neighborhood and short neighborhood properties: consider the normal order tree $(T,\preceq)$ obtained from the order on $\cof{\kappa}^{<\omega}$ after inflating each node $s \in \cof{\kappa}^{<\omega}$ by an interval of order type $\kappa_{\ell(s)}$ corresponding to the nodes of $K_s$. 
    In this setting, $\hat{\mathbb{T}}_{\kappa}$ satisfies the short neighborhood property: let $s \in \cof{\kappa}^{<\omega}$ and $v \in K_s$.
    If $v = \min_{\preceq} K_s$, then
    $N(\odowncl{v}_{\preceq}) \cap \upcl{v}$ consists of $K_s$ and the cliques $K_{s'}$ with
    $s' \in \successor{s}$; otherwise $N(\odowncl{v}_{\preceq}) \cap \upcl{v}$ additionally contains the cliques $K_{s'}$ with
    $s \sqsubseteq s'$ and $\mathrm{dist}_{\cof{\kappa}^{<\omega}}(s,s') = 2$.
    In either case, choosing one element from each clique on the last level provides an antichain of size $\cof{\kappa}$ satisfying \cref{prop:shortNeighborhoodProp:case2}.
\end{example}

Even worse, there are order trees whose existence is independent of $\mathrm{ZFC}$ which do not admit antichains of size at least $\cof{\kappa}$, when $\cof{\kappa}$ is an uncountable cardinal.
Indeed, let $\mathrm{S}$ be a $\cof{\kappa}$--Suslin tree defined on a cofinal set of $\kappa$ comprising regular cardinals. Similar to the above example, after replacing every node $\lambda\in \mathrm{S}$ by an interval of order type $\lambda$, one obtains an order tree of size $\kappa$, but with neither a branch of size $\kappa$ nor an antichain of size $\cof{\kappa}$.

On the other hand, when inspecting the details of \cref{cor:FinishRegularCase}, the reader might notice that the definition of extendability of a spreading pair was not required in full: in fact, it was just necessary to find a strong antichain in the neighborhood of the underlying connected set.
As it turns out, the current general definition of extendable spreading pairs is suitable precisely for including $\mathbb{T}_{\kappa}$ in a tailored generalization of \cref{thm:sufficientCondition} to singular cardinals. 
Indeed, the last section of this paper is devoted to the details of the following extendability lemma:

\begin{restatable}{lemma}{mainSingularLemma}\label{lem:NoCliqueTree}
    Let $\kappa$ be a singular cardinal and let $G$ be a $\kappa$--edge-connected $T$--graph with the down-closed and short neighborhood properties.
    Let $(B, t)$ be a concentrated spreading pair. 
    Then one of the following two cases holds: 
    
    \begin{enumerate}[label=$(L_{\arabic*})$]
        \item\label{item:NoCliqueTree:CliqueTree} $G$ admits a minor embedding of $\mathbb{T}_{\kappa}$, or
        
        \item $(B,t)$ is extendable.
    \end{enumerate}
\end{restatable}

Once the above result is established, the same proof as for \cref{cor:FinishRegularCase} would conclude \cref{thm:mainThm} for the singular cardinal $\kappa$ as well. 
In fact, \cref{example:TreeOfCliques} immediately describes a $\TkStar$ minor in $G$ if the first item \cref{item:NoCliqueTree:CliqueTree} holds. 
Otherwise, every spreading pair as in the premise of \cref{thm:sufficientCondition} is extendable, thus providing a minor embedding of $\TkM$ (and hence of $\TkStar$) in $G$ via \cref{thm:sufficientCondition}. 
To summarize, once \cref{lem:NoCliqueTree} is proven, we can also conclude \cref{thm:mainThm} for all cardinals:

\begin{corollary}
    If $\kappa$ is a singular cardinal, then every $\kappa$--edge-connected graph admits a minor embedding of $\TkStar$.
\end{corollary}

\section{Constructing trees of cofinal cliques}\label{sec:singCard}

From now on in this paper, we work towards proving the aforementioned \cref{lem:NoCliqueTree}, thus also establishing the main \cref{thm:mainThm} for all cardinals.
We fix throughout this section a singular cardinal $\kappa$ and a $\kappa$--edge-connected $T$--graph $G$ satisfying the down-closed and short neighborhood properties.
Furthermore, we also consider an increasing cofinal sequence $\Set{\kappa_{\alpha}}_{\alpha < \cof{\kappa}} \subseteq \kappa$ of regular cardinals greater than $\cof{\kappa}$.

Given a spreading pair $(B,t)$ satisfying the premises of \cref{lem:NoCliqueTree}, we shall assess the extendability of $(B,t)$ through the combinatorial structure of $\Delta = \upcl{t}_{N(B)}$, which is an order tree of size at least $\kappa$ due to the down-closed neighborhood property of $G$.
In particular, we will distinguish whether $\Delta$ admits \cref{prop:shortNeighborhoodProp:case1} or \cref{prop:shortNeighborhoodProp:case2}.

In the former scenario, where $\Delta$ admits an antichain $A$ satisfying \cref{prop:shortNeighborhoodProp:case1}, the proof of \cref{cor:FinishRegularCase} works verbatim to establish the extendability of $(B,t)$ as predicted by \cref{lem:NoCliqueTree}, since $A$ is then already a $(B,\kappa)$--strong antichain.

The challenge in proving \cref{lem:NoCliqueTree} lies in the case where $\Delta$ admits an antichain $A$ satisfying \cref{prop:shortNeighborhoodProp:case2}.
This is the setting where we intend to use the definition of extendability in full, by providing a suitable connected superset $B_*\supseteq B$ that extends $B$ along a cofinal choice of these up-closures from $\set{\upcl{a}_{N(B)}}:{a\in A}$ in order to reach a strong antichain of size $\kappa$ in $N(B_*)$.
As predicted by item \cref{item:NoCliqueTree:CliqueTree} within \cref{lem:NoCliqueTree}, this approach will be successful if we know in advance that $G$ contains no minor embedding of $\mathbb{T}_{\kappa}$.
See also \cref{fig:extension} where the blue nodes are the antichain as in \cref{prop:shortNeighborhoodProp:case2}.

In order to detail this argument, we fix some auxiliary notation for the rest of this section.
First, let $t < s \in T$ and let $B \subseteq T$ be a connected set with $B \cap \upcl{s} = \emptyset$.
The pair $(B,s)$ is then called a \defn{small spreading pair} if $\upcl{s}_{N(B)}$ is of cardinality $< \kappa$. 
For a cardinal $\lambda < \kappa$, we say that a small spreading pair $(B,s)$ is \defn{$\lambda$--extendable from $t$ to $s$} if there exists a connected set $B_*$ with $B \subseteq B _* \subseteq B \cup \oupcl{s}$ and a corresponding $(B_*,\lambda)$--strong antichain $A \subseteq \oupcl{s}$ such that $\left( \bigcup_{a \in A} \odowncl{a} \setminus B_* \right) \cap \upcl{t}$ is connected. 

\noindent
\begin{minipage}[c]{0.67\textwidth}
    On the other hand, let $C \subseteq T$ be an interval of order type $\lambda$ for some uncountable regular cardinal $\lambda$ such that $C$ induces a complete graph in $G$. 
    We say that a sequence $\mathcal{B} \coloneq \set{B_{\alpha}}:{\alpha < \cof{\kappa}}$ of other intervals in $T$ is a \defn{successor sequence of cliques} for $C$ (see also \cref{fig:successorSequence} on the right) if there exists a corresponding disjoint family $\mathcal{W} \coloneq \set{W_{i}}:{i<\lambda}$ of subsets of $C$ satisfying the following conditions:

    \begin{enumerate}[label=$(\dagger_{\arabic*})$, itemsep=3pt]
        \item\label{item:ConvergingSize} $\sup_{\alpha < \cof{\kappa}} \cardinality{B_{\alpha}} = \kappa$;
        
        \item\label{item:MinimumAboveChain} $b_{\alpha} \coloneq \min B_{\alpha}$ belongs to $\Above C$ for every $\alpha < \cof{\kappa}$; 
        
        \item\label{item:MinimumAntichain} $\set{b_{\alpha}}:{\alpha < \cof{\kappa}}$ is an antichain in $T$;
        
        \item\label{item:DominatingBags} for every $\alpha < \cof{\kappa}$ and every $i < \lambda$, there exists a node $y \in W_i$ such that $B_\alpha \subseteq \upcl{b_{\alpha}}_{N(y)}$;
        
        \item\label{item:smallSpreadingTreeClique} there exists a node $c_* \in C \setminus \Set{\min{C}}$, such that $B_\alpha \subseteq N(c_*)$ and $\cardinality{\upcl{b_{\alpha}}_{N(c_*)}} < \kappa$, for every $\alpha < \cof{\kappa}$.
    \end{enumerate}

    \vspace{0.5em}
    Since $C$ is assumed to be complete, the elements of $\mathcal{W}$ define branch sets of a $K_{\lambda}$ minor in $G$.
    Moreover, the last item \cref{item:smallSpreadingTreeClique} implies that $\cardinality{B_{\alpha}} < \kappa$ for every $\alpha < \cof{\kappa}$.
\end{minipage}\hfill
\begin{minipage}[c]{0.4\textwidth}
    \centering
    \captionsetup{type=figure}
    \begin{tikzpicture}[scale=1]
        \tikzset{
            vertex/.style={circle, fill, inner sep=1.7pt},
            edge/.style={draw=black, line width=1pt},
            smallEdge/.style={draw=black, line width=0.6pt},
            every picture/.append style={transform shape},
            directed/.style={->,edge},
            big_cover/.style={line width=5pt, cap=round},
        }
        
        \node[vertex] (C_start) at (0,0) {};
        \node[vertex] (C_w) at (0,1.5) {};
        \coordinate (C_end) at (0,3);
        
        \node[vertex, color=okskyblue] (C_alpha0) at (0,0.5) {};
        \node[vertex, color=okgreen] (C_alpha2) at (0,2) {};
        \node[vertex, color=okpurple] (C_alpha3) at (0,2.6) {};
         
        \draw[edge] (C_start) -- (C_alpha0);
        \draw[edge] (C_alpha0) -- (C_w);
        \draw[edge] (C_w) -- (C_alpha2);
        \draw[edge] (C_alpha2) -- (C_alpha3);
        \draw[directed] (C_alpha3) -- (C_end);

        \begin{scope}[shift=(C_end)]
            \coordinate (r1_0) at (2,0.4);
            \coordinate (tmp1) at (2,0.5);
            \begin{scope}[shift=(tmp1)]
                
                \coordinate (r1_2) at (0,1);
                \coordinate (r1_3) at (0,2);
                \coordinate (r1_4) at (0,3);
                       
            \end{scope}

            \draw[big_cover, opacity=0.5, color=okskyblue] (r1_0) -- (r1_4);
            \node[vertex] (r1_1)  at (2,0.5) {};
            \draw[directed] (r1_1) -- (r1_4);

            \coordinate (r0_0) at (0,0.4);
            \coordinate (r0_2) at (0,1.5);
            \coordinate (r0_3) at (0,2);

            \draw[big_cover, opacity=0.5, color=okpurple] (r0_0) -- (r0_3);
            \node[vertex] (r0_1) at (0,0.5) {};
            \draw[directed] (r0_1) -- (r0_3);

            \coordinate (r00_0) at (0.7,0.4);
            \coordinate (r00_2) at (0.7,2);
            \coordinate (r00_3) at (0.7,2.5);

            \draw[big_cover, opacity=0.5, color=okgreen] (r00_0) -- (r00_3);
            \node[vertex] (r00_1) at (0.7,0.5) {};
            \draw[directed] (r00_1) -- (r00_3);

            \node[fill=white] at (1.3,1.25) {${\cdots}$};
            \node[fill=white] at (2.8,1.25) {${\cdots}$};
        \end{scope}
        
        \draw[smallEdge, opacity=0.6] (C_alpha0) to[bend right=40] (C_alpha2);
        \draw[smallEdge, opacity=0.6] (C_alpha2) to[bend right=40] (C_alpha3);

        \node[anchor=east] at (r1_1) {$b_\alpha$};
        \node[anchor=east] at (r0_1) {$b_{0}$};
        \node[anchor=east] at (r00_1) {$b_{1}$};

        \node[anchor=east] at (C_w) {$c_*$};
        \node[anchor=west] at ($(C_start) + (0.3,0)$) {$C$};
        
        \node[anchor=west] at ($(C_alpha2) + (0.1,0)$) {$W_i$};
        
        \node[anchor=west, color=okskyblue] at (r1_3) {$B_\alpha$};
            
    \end{tikzpicture}
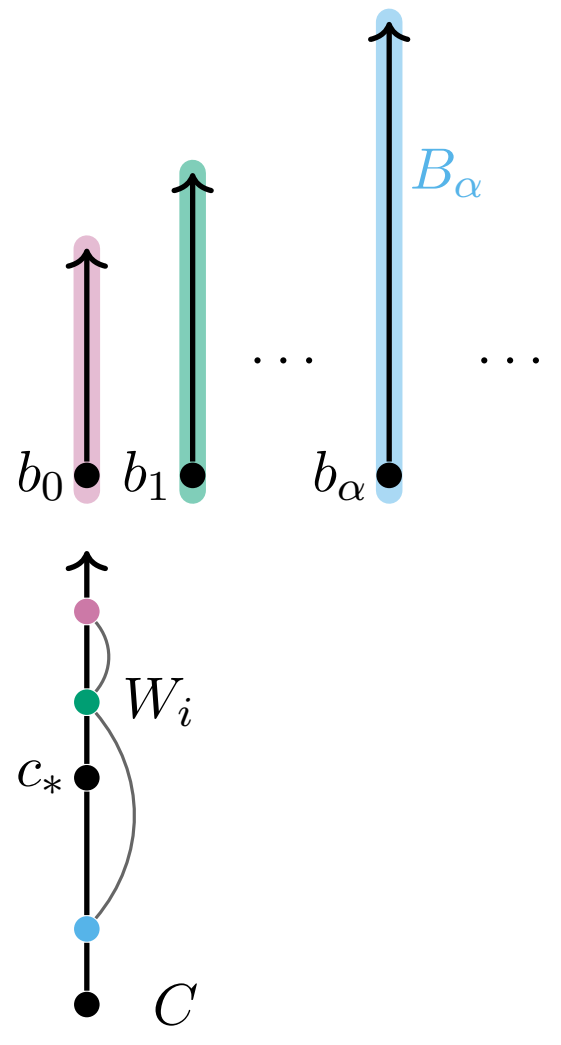
    \caption{Each colored vertex on $C$ is apex to the interval of same color; together these vertices form $W_i$.}\label{fig:successorSequence}
    \vspace{1em}
\end{minipage}

Given a spreading pair $(B,t)$ which admits an antichain $A$ as in item \cref{prop:shortNeighborhoodProp:case2}, a standard example of a small spreading pair has the form $(B,s)$ for $s \in A$. 
In its turn, a $\lambda$--extension of $(B,s)$ from $t$ to $s$ will mostly have the form $B_s\subseteq B\cup \upcl{s}_{N(B)}$, from which we might hope to find a set $A'\subseteq A$ of size $\cof{\kappa}$ and such that $B_*\coloneq B\cup \bigcup_{s\in A'}B_s$ provides an extension of the pair $(B,t)$. 
The following routine lemma will be useful for exhibiting the corresponding strong antichains:

\begin{lemma} \label{lem:extendingToStrongAntichain}
    Suppose that $B$ is a connected set and $A \subseteq N(B)$ is an antichain satisfying $B \cap \bigcup_{a \in A}\upcl{a} = \emptyset$.
    Then there exists $\Phi(a) \in \upcl{a}$ such that $\Phi(A)\coloneq \set{\Phi(a)}:{a \in A}$ defines a $(B, \cardinality{A})$--strong antichain.
\end{lemma}

\begin{proof}
    Set $A' \coloneq \set{a \in A}:{\cardinality{\upcl{a}_{N(B)}}\geq \kappa}$, so that the tuple $(B,a)$ is a spreading pair for every $a\in A'$.
    Define $\Delta_a \coloneq \upcl{a}_{N(B)}$, which is then an order tree of size at least $\kappa$ rooted in $a$.

    By the short neighborhood property with respect to $a$ and the inclusion $\Delta_a \subseteq \upcl{a}_{N(\odowncl{a})}$ there exists a node $\Phi(a) \in \Delta_a$ with $\cardinality{\upcl{\Phi(a)}_{\Delta_a}} < \kappa$.
    If we then set $\Phi(a) = a$ for every $a \in A \setminus A'$, the antichain $\Phi(A) \coloneq \set{\Phi(a)}:{a\in A}$ is as desired.
\end{proof}

On the other hand, if we are given a small spreading pair $(B,s)$ and an interval $C \subseteq \upcl{s}_{N(B)}$ which admits a successor sequence of cliques, we would like to pass to another successor sequence $\mathcal{B}$ for $C$ in which $N(B)$ is disjoint from each interval in $\mathcal{B}$:

\begin{observation}\label{obs:passingToDisjointSuccessorCliques}
    Let $C$ and $(B'_\alpha)_{\alpha < \cof{\kappa}}$ be as in the definition of a successor sequence of cliques.
    If $(B, \min C)$ is a small spreading pair, then there exists another successor sequence of cliques $\mathcal{B}$ for $C$ such that  $\bigcup \mathcal{B}$ is disjoint from $N(B)$.
\end{observation}

\begin{proof}
    We can assume without loss of generality that $\cardinality{B'_\alpha} \geq \kappa_\alpha$, for every $\alpha < \cof{\kappa}$ (due to \cref{item:ConvergingSize}).
    Moreover, if $B_\alpha \subseteq B'_\alpha$ denotes the subinterval of $B'_\alpha$ rooted in $\min B'_\alpha$ and consisting of all nodes up to height $\kappa_\alpha$ in $B'_\alpha$, then $B_\alpha$ has order type exactly $\kappa_\alpha$; in particular, $(B_\alpha)_{\alpha < \cof{\kappa}}$ is again a successor sequence of cliques for $C$.

    Now, since $\upcl{\min C}_{N(B)}$ has cardinality $< \kappa$ by assumption, there exists some $\alpha < \cof{\kappa}$ such that $\cardinality{\upcl{\min C}_{N(B)}} < \kappa_\alpha$.
    Given $\beta \geq \alpha$, it follows by the down-closed neighborhood property that $N(B) \cap B_\beta$ is bounded in the order type of $B_\beta$ (due to the regularity of $\kappa_\beta$).
    In particular, there exists a node $t_\beta \in B_\beta$ such that $\upcl{t_\beta}_{B_\beta}$ is disjoint from $N(B)$.
    By regularity of $\kappa_\beta$ again, the interval $\upcl{t_\beta}_{B_\beta}$ still has order type $\kappa_\beta$ and $\set{\upcl{t_\beta}_{B_\beta}}:{\alpha \leq \beta < \cof{\kappa}}$ is again a successor sequence of cliques for $C$, where every interval is now disjoint from $N(B)$.
\end{proof}

As an intuition behind the next lemma, we highlight that the definition of successor sequence of cliques is in some sense dual to the extendability from $t$ to $s$:

\begin{lemma}\label{lemma:SendingAbove}
    Fix a regular cardinal $\lambda$ with $\cof{\kappa} < \lambda < \kappa$.
    Let $t < s \in T$ be given such that $(B,s)$ is a small spreading pair, $C \subseteq \upcl{s}_{N(B)}$ is an interval of order type $\lambda$ rooted in $s$, $[t,s] \setminus B$ is connected and $\cardinality{N(\odowncl{t} \setminus B) \cap \upcl{t}} < \kappa$. 
    Then, one of the items below holds:

    \begin{enumerate}[label=$(J_{\arabic*})$]
        \item \label{item:FindingAntichain} $(B,s)$ is $\lambda$--extendable from $t$ to $s$, or

        \item \label{item:FindingSuccessors} there exists an interval $C_* \subseteq C$ of order type $\lambda$ which is a clique and admits a successor sequence of cliques.
    \end{enumerate}
    
\end{lemma}

\begin{proof}
    We first argue that there exists a node $v_* \in C$ satisfying $\oupcl{v_*}_C\subseteq N(v_*)$: consider the connected induced subgraph $G[C]$ of size $\lambda$; since $\lambda$ is regular, there exists some node $v_* \in C$ with $\lambda$ many neighbors in $C$.
    But then $N(v_*)$ is cofinal in the order of $C$ and by the down-closed neighborhood property, $v_*$ neighbors every node on $\oupcl{v_*}_C$.
    Also, $v_* \notin B$ since $B$ is disjoint from $\upcl{s}$ by hypothesis.
    In that setting, we shall conclude \cref{item:FindingAntichain} or \cref{item:FindingSuccessors} via the analysis of the following cases:

    \begin{case}
        The set $\Gamma \coloneq \set{c \in C}:{\cardinality{N(c) \cap \Above C}<\kappa} \cap \oupcl{v_*}$ is of size $\lambda$.
    \end{case}

       Within this case, we denote by $c^+\in C$ the unique successor in $C$ of a given $c\in C$. 
       Then, define $\Gamma_1 \coloneq \set{c \in C \cap \oupcl{v_*}}:{\upcl{c^+}_{N(c)}\setminus (C \cup \Above C) \neq \emptyset}$.
       In the following, we further distinguish whether $\Gamma_1$ has size $\lambda$ or not, aiming in both cases to establish item \cref{item:FindingAntichain}.

       \begin{subcase}\label{subcase:11}
           $\Gamma_1$ has size $\lambda$.
       \end{subcase}

       Let $(c_i)_{i < \lambda}$ be an ordering of $\Gamma_1$ that agrees with the ordering of $(T, \leq)$.
       Then, for each $i < \lambda$, fix $a_i \in \upcl{c_i^+}_{N(c_i)} \setminus (C\cup \Above C)$ as ensured by the definition of $\Gamma_1$. 
       In particular, $a_i$ is comparable with $c_i$ but with only $<\lambda$ many vertices from $C$. 
       Since $\lambda$ is regular, there exists an unbounded index set $I \subseteq \lambda$ for which $(c_i,c_i^+,\sup_C(C \cap \downcl{a_i}))$ is distinct in every entry for pairwise distinct indices $i \in I$.
       In fact, a standard recursive construction of $I$ provides $c_i < c_i^+ \leq \sup_C(C \cap \downcl{a_i}) < c_j$ if $i < j \in I$.
       Hence, $A \coloneq \Set{a_i}_{i\in I}$ is a $\lambda$--sized antichain in the neighborhood of the connected set $B_* \coloneq B \cup \Set{c_i}_{i\in I} \subseteq B \cup \upcl{s}_{N(B)}$.
       \cref{lem:extendingToStrongAntichain} then admits a $(B_*,\lambda)$--strong antichain $\set{\Phi(a)}:{a\in A}$ with $\Phi(a) \in \upcl{a}$ for every $a\in A$.
       Recalling the choice of $v_*$ and $I$, the induced subgraph on $[s,v_*] \cup [c_i^+,\Phi(a_i))$ is disjoint from $B_*$ and connected, for every $i \in I$. 
       As sketched by the bold subgraph in \cref{fig:subcase11}, this witnesses the connectedness of $\left( \bigcup_{a\in A} \odowncl{\Phi(a)} \setminus B_* \right) \cap \upcl{s}$. 
       Together with the connectedness of $[t,s] \setminus B$ as well as the fact that $s$ is an element of this set, both statements resulting from the hypotheses, it follows that $(B,s)$ is $\lambda$--extendable from $t$ to $s$.

    \begin{subcase}\label{subcase:12}
        $\Gamma_1$ has size $<\lambda$.
    \end{subcase}

    Due to this subcase and due to the regularity of $\lambda$, we can fix a node $w_* \in \oupcl{v_*}_C$ with $\upcl{w_*}_C \cap \Gamma_1 = \emptyset$.
    Then let $c \in \Gamma \cap \upcl{w_*}_C$ be given.
    It follows that $c$ has no neighbors in $\upcl{c^+} \setminus (C \cup \Above C)$ since $c \notin \Gamma_1$, and that $c$ has $<\kappa$ many neighbors in $\Above C$ by definition of $ \Gamma$.
    In particular, since $C$ has size $\lambda < \kappa$, it follows that the node $c$ has $<\kappa$ many neighbors in the set $\upcl{c^+}$.
    But as a consequence of $\kappa$--edge-connectivity of $G$, the node $c$ has at least $\kappa$ many neighbors and, hence, almost all neighbors of $c$ are contained in the set $\Delta_c \coloneq \upcl{c} \setminus \upcl{c^+}$.

    Recalling that both $\upcl{s}_{N(B)}$ and $N(\odowncl{t} \setminus B) \cap \upcl{t}$ have cardinality $< \kappa$ by hypothesis, the inequality $\cardinality{N(\odowncl{t} \cup B) \cap \upcl{s}} < \kappa$ now holds.
    Then let $\Delta_c^{\circ} \coloneq \Delta_c\setminus \Set{c}$.
    Since $\Delta_c$ induces an edge cut of size at least $\kappa$, and since $N(\upcl{w_*}_C \setminus \Set{c}) \cap \Delta_c^{\circ} = \emptyset$ by choice of $w_*$, it follows that $\cardinality{N([t,w_*) \setminus B) \cap \Delta_c^{\circ}} \geq \kappa$.
    This leads then to the following conclusion:

    \begin{enumerate}[label=$(\star)$]
        \item \label{lemma:SendingAbove:StarProp} $([t,w_*)\setminus B) \cap N(\Delta_c^{\circ}) \neq \emptyset$ for every $c \in \Gamma \cap \upcl{w_*}_C$. 
    \end{enumerate}

    Fix again an enumeration $(c_i)_{i<\lambda}$ for $\Gamma \cap \upcl{w_*}_C$ which agrees with the ordering inherited from $(T, \leq)$.
    By \cref{lemma:SendingAbove:StarProp} we can fix a pair $(d_i, a_i)$ where $d_i \in [t,w_*) \setminus B$ and $a_i \in N(d_i) \cap \Delta_{c_i}^{\circ}$ for each $i < \lambda$.
    Let $B_* \coloneq B \cup \bigcup_{i < \lambda} [c_i,a_i)$ and $A \coloneq \Set{a_i}_{i < \lambda}$ and note that $B_*$ is connected since $c_i \in N(B)$.
    Furthermore, $A \subseteq N(B_*)$ by definition of $B_*$, while $A$ defines an antichain of size $\lambda$.
    Thus, \cref{lem:extendingToStrongAntichain} returns a $(B_*,\lambda)$--strong antichain $\set{\Phi(a)}:{a\in A}$ with $\Phi(a) \in \upcl{a}$ for every $a \in A$.
    It remains to verify that $\left( \bigcup_{a \in A} \odowncl{\Phi(a)} \setminus B_* \right) \cap \upcl{t}$ is connected, as suggested by the bold subgraph sketched in \cref{fig:subcase12}.
    
    We first claim that $[t,w_*) \setminus B_*$ is connected.
    In fact, it holds that $[t,w_*) \setminus B_* = ([t,s] \setminus B) \cup [s,w_*)$ due to $B_* \subseteq B \cup \upcl{w_*}$ and since $B$ is disjoint from $\upcl{s}$ by hypothesis.
    The claim then follows by the connectedness assumption on $[t,s] \setminus B$.

    Now, let $a_i \in A$ be given and write the set $D_{a_i} \coloneq (\odowncl{\Phi(a_i)} \setminus B_*) \cap \upcl{t}$ as $([t, w_*) \setminus B_*) \cup ([w_*, c_i) \setminus B_*) \cup [a_i, \Phi(a_i))$.
    Let us recall that $d_i \in [t,w_*) \setminus B = [t,w_*) \setminus B_*$ and that $d_i$ neighbors $a_i$.
    In particular, due to the down-closed neighborhood property, $d_i$ neighbors every node along the interval $(d_i, a_i]$, which strictly contains $[w_*, c_i) \setminus B_*$.
    All in all, $D_{a_i}$ is connected and, since any two sets $D_a$ and $D_{a'}$ for distinct $a, a' \in A$ contain the connected set $[t,w_*) \setminus B_*$ in common, it follows that $(B,s)$ is $\lambda$--extendable from $t$ to $s$.

    \begin{figure}
    \centering
    \scalebox{0.8}{
    \begin{minipage}{0.5\textwidth}
        \centering
        \begin{tikzpicture}[
        scale=0.5,
        TreeNode/.style = {circle, fill=black, inner sep=2pt},
        Special/.style = {circle, fill=black, inner sep=3pt},
        Cnode/.style = {circle, fill=okorange, inner sep=3pt},
        Anode/.style = {circle, fill=okgreen, inner sep=3pt},
        directed/.style = {->, black, line width=1pt},
        ApexEdge/.style = {gray, line width=0.4pt}
    ]
    \foreach \i in {0,...,13}{
        \node[TreeNode] (x\i) at (0,\i) {};
    }
        \draw[directed] (x0) -- (0,14);
        \draw[dashed, okorange, line width=0.6pt] (-4,-6.25) to[in=-60, out=20] (0,8);
         \filldraw[decorate, okorange, fill=okorange, line width=0.6pt] 
        (-5,-12) to[in=210, out=170] (-3.5,-4.25) decorate[decoration={snake, amplitude=2mm, segment length=8mm}]{--(-4.5,-6.75)} to[in=170, out=210] (-5,-12);

        \node[Cnode] (ci) at (0,3) {};
        \node[Anode] (ai) at (3.7,8.3) {};
        \draw[decorate, decoration={snake, amplitude=0.5mm}, line width=1pt] (ai) -- (0,6);  
        \draw[dashed, line width=0.6pt] (ci) -- (ai);

        \node[Cnode] (cj) at (0,8) {};
        \node[Anode] (aj) at (3.7,13.3) {};
        \draw[decorate, decoration={snake, amplitude=0.5mm}, line width=1pt] (aj) -- (0,11);  
        \draw[dashed, line width=0.6pt] (cj) -- (aj);
        
        \node[Special] (v) at (0,0) {};

        \node[Special] (s) at (-3,-3) {};
        \node[Special] (t) at (-5,-8) {};

        \draw[decorate, decoration={snake, amplitude=0.5mm}, line width=1pt] (s) -- (v);
        \draw[decorate, black, decoration={snake, amplitude=0.5mm}, line width=1pt] (t) -- (s);
        
        \filldraw[decorate, okorange, fill=okorange, decoration={snake, amplitude=0.5mm}, line width=0.6pt] (-5,-12) ellipse (1cm and 2cm);

        \foreach \i in {2,...,13}{
        \draw[dashed, line width=0.6pt] (v) to[in=240, out=140] (x\i);
    }
         \draw[dashed, okorange, line width=0.6pt] (-4,-12) to[in=-60, out=30] (ci);

         \node (t1) at (-4.3,-8) {$t$};
         \node (s1) at (-2.3,-3) {$s$};
         \node (v1) at (0.9,-0.3) {$v_*$};
         \node (ci1) at (0.9,3) {\textcolor{okorange}{$c_i$}};
         \node (cj1) at (0.9,8) {\textcolor{okorange}{$c_j$}};
        \node (b1) at (-7.5,-10) {\textcolor{okorange}{$B$}};
        \node (ai1) at (3.7,9.3) {\textcolor{okgreen}{$\Phi(a_i)$}};
        \node (aj1) at (3.7,14.3) {\textcolor{okgreen}{$\Phi(a_j)$}};
        \end{tikzpicture}
        \caption{The sequence $\Set{c_i}_{i\in I}$ and the corresponding antichain $\Set{\Phi(a_i)}_{i\in I}$ as constructed within \cref{subcase:11}.}
        \label{fig:subcase11}
    \end{minipage}
    }
    \hfill
    \scalebox{0.8}{
    \begin{minipage}{0.5\textwidth}
        \centering
        \begin{tikzpicture}[
        scale=0.5,
        TreeNode/.style = {circle, fill=black, inner sep=2pt},
        Special/.style = {circle, fill=black, inner sep=3pt},
        Cnode/.style = {circle, fill=okorange, inner sep=3pt},
        Anode/.style = {circle, fill=okgreen, inner sep=3pt},
        directed/.style = {->, black, line width=1pt},
        ApexEdge/.style = {gray, line width=0.4pt}
    ]

    \draw[dashed, line width=0.6pt] (0.95,6.2) to[in=270, out=290] (3.9,10.2);
    \draw[dashed, line width=0.6pt] (0.95,11.2) to[in=350, out=290] (3.9,15.2);
    \draw[dashed, line width=0.6pt] (-1,-1) to[in=320, out=0] (1.9,7.2);
    \draw[dashed, line width=0.6pt] (-2,-2) to[in=150, out=150] (1.9,12.2);
    \foreach \i in {1,...,5}{
         \draw[dashed, line width=0.6pt] (-2,-2) to[in=180, out=150] (0,\i);
     }

     \foreach \i in {7,...,10}{
        \draw[dashed, line width=0.6pt] (-2,-2) to[in=180, out=150] (0,\i);
    }
    \foreach \i in {0,...,13}{
        \node[TreeNode] (x\i) at (0,\i) {};
    }
        \draw[directed] (x0) -- (0,14);
        \draw[dashed, okorange, line width=0.6pt] (-4,-6.25) to[in=-60, out=20] (0,11);
         \filldraw[decorate, okorange, fill=okorange, line width=0.6pt] 
        (-5,-12) to[in=210, out=170] (-3.5,-4.25) decorate[decoration={snake, amplitude=2mm, segment length=8mm}]{--(-4.5,-6.75)} to[in=170, out=210] (-5,-12);

        \node[Cnode] (ci) at (0,6) {};
        \node[TreeNode] (ai) at (1.9,7.2) {};
        \node[Anode] (pai) at (3.9,10.2) {};
        \draw[decorate, okorange, decoration={snake, amplitude=0.5mm}, line width=1pt] (ai) -- (0,6);  
        \draw[decorate, decoration={snake, amplitude=0.5mm}, line width=1pt] (ai) -- (pai);  
         \node[TreeNode] (di) at (-1,-1) {};

     \foreach \i in {1,...,5}{
        \draw[dashed, line width=0.6pt] (di) to[in=0, out=0] (0,\i);
    }

         \node[Cnode] (cj) at (0,11) {};
        \node[TreeNode] (aj) at (1.9,12.2) {};
        \node[Anode] (paj) at (3.9,15.2) {};
        \draw[decorate, okorange, decoration={snake, amplitude=0.5mm}, line width=1pt] (aj) -- (0,11);  
         \draw[decorate, decoration={snake, amplitude=0.5mm}, line width=1pt] (aj) -- (paj);  
         \node[TreeNode] (dj) at (-2,-2) {};

        \node[Special] (v) at (0,0) {};

        \node[Special] (s) at (-3,-3) {};
        \node[Special] (t) at (-5,-8) {};

        \draw[decorate, decoration={snake, amplitude=0.5mm, segment length=1.5mm}, line width=1pt] (s) -- (dj);
        \draw[decorate, decoration={snake, amplitude=0.5mm, segment length=1.5mm}, line width=1pt] (dj) -- (di);
        \draw[decorate, decoration={snake, amplitude=0.5mm, segment length=1.5mm}, line width=1pt] (di) -- (v);
        \draw[decorate, black, decoration={snake, amplitude=0.5mm}, line width=1pt] (t) -- (s);

        \filldraw[decorate, okorange, fill=okorange, decoration={snake, amplitude=0.5mm}, line width=0.6pt] (-5,-12) ellipse (1cm and 2cm);

         \draw[dashed, okorange, line width=0.6pt] (-4,-12) to[in=-60, out=30] (ci);

         \node (t1) at (-4.3,-8) {$t$};
         \node (s1) at (-3.7,-3) {$s$};
         \node (dj1) at (-1.3,-2.5) {$d_j$};
         \node (ai1) at (1.3,7.5) {$a_i$};
         \node (aj1) at (1.5,13) {$a_j$};
         \node (di1) at (-0.3,-1.5) {$d_i$};
         \node (paj1) at (5.2,10.5) {\textcolor{okgreen}{$\Phi(a_i)$}};
         \node (pai1) at (5.5,15) {\textcolor{okgreen}{$\Phi(a_j)$}};
         \node (v1) at (-0.9,0.2) {$w_*$};
         \node (ci1) at (-0.9,6) {\textcolor{okorange}{$c_i$}};
         \node (cj1) at (-0.9,11) {\textcolor{okorange}{$c_j$}};
        \node (b1) at (-3.5,-11) {\textcolor{okorange}{$B$}};
        \end{tikzpicture}
        \caption{The sequence $\Set{c_i}_{i < \lambda}$ and the corresponding antichain $\Set{\Phi(a_i)}_{i < \lambda}$ as constructed within \cref{subcase:12}.}
        \label{fig:subcase12}
    \end{minipage}
    }
    \end{figure}

    \begin{case}
        The set $\Gamma \coloneq \set{c \in C}:{\cardinality{N(c) \cap \Above C}<\kappa}\cap \oupcl{v_*}$ is of size $<\lambda$.
    \end{case}

    Once more due to the regularity of $\lambda$, we can fix a node $u_* \in \oupcl{v_*}_C$ such that $\upcl{u_*}_C\cap \Gamma = \emptyset$.
    In this case, the chain $C_* \coloneq \upcl{u_*}_C$ induces a clique in $G$: after all, for every $c \in C_*$, the definition of $\Gamma$ provides a neighbor of $c$ in $\Above C$ and, thus, the down-closed neighborhood property of $G$ implies that $\oupcl{c}_C \subseteq N(c)$.
    Now, if there exists $c \in C_*$ whose neighborhood $N(c) \cap \Above C$ contains an antichain $A$ of size $\lambda$, we conclude item \cref{item:FindingAntichain}.
    Indeed, an application of \cref{lem:extendingToStrongAntichain} with respect to the connected set $B_* \coloneq B \cup \Set{c}$ provides a $(B_*,\lambda)$--strong antichain $\Phi(A) = \set{\Phi(a)}:{a \in A}$ with $\Phi(a) \in \upcl{a}$ for every $a\in A$. 
    In this case, the connectedness of $\left( \bigcup_{a \in A} \odowncl{\Phi(a)} \setminus B_* \right) \cap \upcl{t}$ follows from the connectedness of $[t,v_*] \setminus B_* = ([t,s] \setminus B) \cup [s,v_*]$ (since $[t,s] \setminus B$ is connected and $B$ is disjoint from $\upcl{s}$ by hypothesis), as well as from the choice of $v_*$ and the fact that $\leq_T$--minimal nodes of $\Above C$ have cofinally many neighbors in $C_*$.
    Hence, in the following we can assume that $\Delta'_c \coloneq N(c) \cap \Above C$ is of size at least $\kappa$, but contains no antichains of size $\lambda$, for every $c \in C_*$.
    Under this assumption we are able to conclude the following:

    \begin{claim}\label{claim:ExistenceSequenceChains}
        For every $c\in C_*$ there is a sequence $(B^c_\alpha)_{\alpha < \cof{\kappa}}$ of pairwise disjoint intervals in $\Delta_c'$ such that the minimal elements $\Set{\min B_{\alpha}^c}_{\alpha < \cof{\kappa}}$ define an antichain in $T$, each interval $B^c_\alpha$ is of order type $\kappa_\alpha$ and $(\Set{c}, \min B_{\alpha}^c)$ is a small spreading pair.
    \end{claim}

    \begin{proof}[Proof of~\cref{claim:ExistenceSequenceChains}]\renewcommand{\qedsymbol}{$\blacksquare$}
        It is enough to define an antichain $\Set{a_\alpha}_{\alpha < \cof{\kappa}} \subseteq \Delta'_c$ such that $\upcl{a_\alpha}_{N(c)}$ has cardinality $< \kappa$ and height at least $\kappa_\alpha^+$: letting $B^c_\alpha \subseteq \upcl{a_\alpha}_{N(c)}$ be the strict down-closure (that is, excluding the maximal element) of any node at height $\kappa_{\alpha}$ in $\upcl{a_\alpha}_{N(c)}$  defines an interval of order type $\kappa_\alpha$ rooted in $a_\alpha$.
        In particular, the sequence $(B^c_\alpha)_{\alpha < \cof{\kappa}}$ would then satisfy the claim.

        Now, by the down-closed neighborhood property it follows that $\Delta'_c$ is a disjoint union of order trees of the form $\Delta'_c \cap \upcl{l} = N(c) \cap \upcl{l}_{N(c)}$, where $l$ is a $\leq$--minimal element of $\Delta'_c$.
        Let us first assume the existence of some such $l$ for which $\upcl{l}_{N(c)}$ is of size at least $\kappa$.
        In this case, recall that the paragraph right before \cref{claim:ExistenceSequenceChains} excluded antichains of size $\lambda$ in $\Delta'_c \supseteq \upcl{l}_{N(c)}$.
        In particular, the short neighborhood property of $G$ implies that $\upcl{l}_{N(c)}$ satisfies \cref{prop:shortNeighborhoodProp:case2}, from which one obtains an antichain $\Set{a_\alpha}_{\alpha < \cof{\kappa}} \subseteq \upcl{l}_{N(c)}$ such that $\upcl{a_\alpha}_{N(c)}$ has cardinality $< \kappa$ and height at least $\kappa^+_{\alpha}$.

        Hence, in the following we may assume that $\cardinality{\upcl{l}_{N(c)}}<\kappa$ for every $\leq$--minimal element $l$ of $\Delta_c'$.
        In this case, and by the fact that $\cardinality{\Delta_c'} \geq \kappa$, there exists for each $\alpha < \cof{\kappa}$ some $\leq$--minimal node $a^c_{\alpha} \in \Delta'_c$ for which $\upcl{a^c_{\alpha}}_{N(c)}$ has cardinality at least $\max{\Set{\lambda^+, \kappa_{\alpha}^{++}}}$. 
        In particular, we can choose the nodes $a^c_\alpha$ to be pairwise distinct and, thus, to induce an antichain in $T$.
        On the other hand, since $\upcl{a^c_\alpha}_{N(c)} \subseteq \Delta'_c$ contains no antichains of size $\lambda$ by assumption, the height of $\upcl{a^c_\alpha}_{N(c)}$ is at least $\theta \coloneq \max{\Set{\lambda, \kappa_{\alpha}^+}}$: indeed, relying on the regularity of $\theta^+ \leq \cardinality{\upcl{a^c_\alpha}_{N(c)}}$ and the fact that each level is of cardinality $< \lambda$, it follows that $\upcl{a^c_\alpha}_{N(c)}$ has height at least $\theta$.
        Hence, $\upcl{a^c_\alpha}_{N(c)}$ has height at least $\kappa^+_{\alpha}$ and is of cardinality $< \kappa$.
    \end{proof}

    The proof of the current lemma is then finished by the following last case distinction, where we work towards obtaining \cref{item:FindingSuccessors} as an output.

    \begin{subcase}\label{lemma:SendingAbove:Case21}
       For every $c \in C_*$ there exists $x_c \in \bigcup_{\alpha < \cof{\kappa}} B^c_\alpha$ such that $\cardinality{N(x_c) \cap C_*} < \lambda$.
    \end{subcase}

    In this subcase, let $A$ be the set of $\leq$--minimal elements of $\set{x_c}:{c \in C_*}$.
    We first claim that $\cardinality{A} = \lambda$.
    Suppose for a contradiction that $\cardinality{A} < \lambda$, which allows us to fix some node $c \in C_* \setminus (N(A) \cap C_*)$, since $C_*$ has size $\lambda$ but $N(A) \cap C_*$ now has size $< \lambda$ by the premises of the current \cref{lemma:SendingAbove:Case21}.
    The node $c$ is then associated with a corresponding node $x_c \in N(c)$, from which it follows that $x_c \notin A$ and $a < x_c$ for some $a \in A$.
    The down-closed neighborhood property then implies that $c$ neighbors $a$, so that $c \in N(A) \cap C_*$, a contradiction.
    
    One can then fix an unbounded subset $I \subseteq C_*$ such that $A = (x_c)_{c \in I}$. 
    In this case, for any $c \in I$, fix $l_c \coloneq \min \downcl{x_c} \cap \Above C$ and note that $l_c \neq x_{c'}$ if $c,c'\in I$ are distinct.
    Due to the normality property \cref{item:CofinalNeighborhoods} of $T$ and the regularity of $\lambda = \cardinality{C_*}$, the neighborhood of $l_c$ in $C_*$ defines a $\lambda$--sized set $N_c$.
    Hence, unless by passing $I$ and $N_c$ to suitable $\lambda$--sized subsets, we may assume that the elements from $\Set{I} \cup \Set{N_{c}}_{c \in I}$ are pairwise disjoint.
    As in the previous cases, we can now apply \cref{lem:extendingToStrongAntichain} with respect to the set $B_* \coloneq B \cup I$ to obtain a $(B_*,\lambda)$--strong antichain $\set{\Phi(a)}:{a \in A}$ with $\Phi(a) \in \upcl{a}$.
    It remains to check that $\left( \bigcup_{a\in A}\odowncl{\Phi(a)}\setminus B_* \right) \cap \upcl{t}$ is connected.
    Indeed, given $x_c \in A$, we can write $D_c \coloneq (\odowncl{\Phi(x_c)} \setminus B_*) \cap \upcl{t}$ as the disjoint union $([t, u_*) \setminus B) \cup (C_* \setminus B_*) \cup [l_c, x_c) \cup [x_c, \Phi(x_c))$, because $B_* \setminus B \subseteq C_*$ and $B \cap \upcl{s} = \emptyset$ by assumption.
    From the hypotheses it also follows that $[t, u_*) \setminus B$ is connected.
    We further have $v_* \in [t,u_*) \setminus B$ since $C_*$ was chosen not to contain $v_*$, while $[l_c,x_c] \cap B_* = \emptyset$ since $I \subseteq C_*$ and $[l_c,x_c] \subseteq \Above C$.
    In its turn, $l_c$ has a neighbor in $C_*\setminus B_*$ due to the disjointness of $N_c$ and $I$.
    By the choice of $v_*$ it follows that $v_*$ is adjacent to such a neighbor of $l_c$, proving thus the connectedness of $D_c$.
    Hence, since $D_c$ and $D_{c'}$ have the connected set $[t, u_*) \setminus B$ in common for $c, c' \in I$, the extendability of $(B,s)$ from $t$ to $s$ follows. 
    
    The proof is then finished after analyzing the following remaining scenario.

    \begin{subcase}\label{lemma:SendingAbove:Case22}
        There exists $c_* \in C_*$ such that 
        $\cardinality{N(x) \cap C_*} = \lambda$ for every $x \in \bigcup_{\alpha < \cof{\kappa}} B_\alpha^{c_*}$.
    \end{subcase}

    Firstly, we can assume without loss of generality that there exists such a $c_*$ that is not the minimum element of $C_*$: otherwise, the same proof as in \cref{lemma:SendingAbove:Case21} with respect to $\oupcl{c_*}_{C_*}$ would obtain \cref{item:FindingAntichain} as an output.

    Now, if $\alpha < \cof{\kappa}$ is given with $\kappa_\alpha \geq \lambda^+$, then the set $\Gamma_\alpha \coloneq \set{c \in C_*}:{B_\alpha^{c_*} \subseteq N(c)}$ has cardinality $\lambda$: otherwise, there would exist $z_* \in C_*$ such that every node in $\upcl{z_*}_{C_*}$ has fewer than $\kappa_{\alpha}$ many neighbors in $B_{\alpha}^{c_*}$, due to the regularity of the order type of $C_*$ and the down-closed neighborhood property of $G$. 
    But since $\kappa_{\alpha}> \lambda$, this would imply $N(x) \cap \upcl{z_*}_{C_*} = \emptyset$ (and hence $\cardinality{N(x) \cap C_*} < \lambda$) for some $x \in B_{\alpha}^{c_*}$, which is a contradiction.

    Denote then $B_{\alpha} \coloneq B_{\alpha}^{c_*}$ and $b_{\alpha} \coloneq \min B_{\alpha}$ for every $\alpha < \cof{\kappa}$ with $\kappa_\alpha \geq \lambda^+$.
    Since $\cardinality{\Gamma_{\alpha}} = \lambda$, let $\Set{W_i}_{i<\lambda}$ be any system of disjoint subsets of $C_*$ with $W_i \cap \Gamma_{\alpha} \neq \emptyset$ for every $i < \lambda$ and $\alpha < \cof{\kappa}$ with $\kappa_\alpha \geq \lambda^+$.
    In particular, the family $\Set{W_i}_{i<\lambda}$ meets item \cref{item:DominatingBags} for the family of intervals $\mathcal{B}\coloneq \set{B_{\alpha}}:{\alpha < \cof{\kappa}, \kappa_\alpha \geq \lambda^+}$.
    The node $c_*$ then satisfies all conditions of \cref{item:smallSpreadingTreeClique}, where $\cardinality{\upcl{b_\alpha}_{N(c_*)}}<\kappa$ holds by the assumption that $(\Set{c_*}, b_{\alpha})$ is a small spreading pair (as in \cref{claim:ExistenceSequenceChains}).
    All in all, $\mathcal{B}$ is indeed a successor sequence of cliques for $C_*$, where items \cref{item:ConvergingSize}--\cref{item:MinimumAntichain} are inherited from \cref{claim:ExistenceSequenceChains}.
\end{proof}

Regarding the definition of successor sequence of cliques, the statement of the above \cref{lemma:SendingAbove} might be seen as an inductive step towards constructing a minor embedding of $\mathbb{T}_{\kappa}$ in $G$, or towards finishing the proof of \cref{lem:NoCliqueTree}. 
More precisely, iterating \cref{lemma:SendingAbove} yields the following conclusion.

\begin{lemma}\label{lemma:TreeOfCliquesConstruction}
    Let $\lambda$ be a regular cardinal with $\cof{\kappa} < \lambda < \kappa$ and let $(B, t)$ be a concentrated spreading pair.
    Suppose further that there exists $s \in T$ with $t < s$ such that $(B,s)$ is a small spreading pair and $\upcl{s}_{N(B)}$ has height at least $\lambda^+$.
    Then, one of the following two cases holds:

    \begin{enumerate}[label=$(L_{\arabic*})$]
        \item\label{item:OutputTreeOfCliques} $G$ contains $\mathbb{T}_{\kappa}$ as a minor;
        
        \item \label{item:TreeOfCliquesConstruction:Extension} $(B,s)$ is $\mu$--extendable from $t$ to $s$, for some $\lambda \leq \mu < \kappa$.
    \end{enumerate}
\end{lemma}

\begin{proof}
    Let us assume $\kappa_0>\lambda$ without loss of generality and recall the graph-theoretical tree structure of $(\cof{\kappa}^{<\omega},\sqsubseteq)$ as in the definition of $\mathbb{T}_{\kappa}$ given by \cref{example:TreeOfCliques}. 
    We denote by $\ell(p)$ the last element of any given sequence $p \in \cof{\kappa}^{<\omega}\setminus \Set{\emptyset}$.
    Set $\cof{\kappa}^{\leq 0} \coloneq \Set{\emptyset}$ and, for each $n \geq 1$, we let $\cof{\kappa}^{\leq n}$ be the set of finite increasing sequences in $\cof{\kappa}$ of length at most $n$. 
    Equivalently, $\cof{\kappa}^{\leq n}$ comprises all nodes of height at most $n$ in $\cof{\kappa}^{<\omega}$.
    Under the assumption that \cref{item:TreeOfCliquesConstruction:Extension} does not hold, our aim is then to define a minor embedding $\set{B_p}:{p \in \cof{\kappa}^{<\omega}}$ of $\cof{\kappa}^{<\omega}$ in $\upcl{s}$ such that, for every $p \in \cof{\kappa}^{<\omega}$, the family $\set{B_q}:{q \in \mathrm{succ}_{\sqsubseteq}(p)}$ is a successor sequence of cliques for $B_p$.
    
    Since $\upcl{s}_{N(B)}$ has height at least $\lambda^+$, there exists some node at height $\lambda$ and, in particular, some interval $C \subseteq \upcl{s}_{N(B)}$ of order type $\lambda$ rooted in $s$.
    The construction of $\set{B_p}:{p\in \cof{\kappa}^{<\omega}}$ is then recursive and, since the triple $(B,t,s)$ does not address condition \cref{item:TreeOfCliquesConstruction:Extension}, we start by considering $B_{\emptyset} \subseteq C$ as the clique of order type $\lambda$ provided by item \cref{item:FindingSuccessors} from \cref{lemma:SendingAbove} and which admits a successor sequence of cliques $\mathcal{B}_{\emptyset} \coloneq \set{C_q}:{q<\cof{\kappa}}$.
    In particular, by \cref{obs:passingToDisjointSuccessorCliques} we can assume that $C_q$ is disjoint from $N(B)$ for every $q <\cof{\kappa}$.
    For inductive purposes, set also $X_{\emptyset}\coloneq \emptyset$.

    For some $n \in \mathbb{N}$, suppose now that $\set{B_p}:{p \in \cof{\kappa}^{\leq n}}$ is already defined and that $B_p$ is a clique of order type $\kappa_{\ell(p)}$ for each $p\in \cof{\kappa}^{\leq n}\setminus \Set{\emptyset}$. 
    For distinct $p,p'\in \cof{\kappa}^{\leq n}$ with $p \sqsubseteq p'$, it is also ensured by induction that $b_{p'} \coloneq \min B_{p'}$ belongs to $\Above B_p$, i.e., that every element of $B_{p'}$ is greater than every element of $B_p$.
      
    Moreover, if $p \neq \emptyset$ is a leaf of $\cof{\kappa}^{\leq n}$ and $p_0 p_1 p_2 \dots p_n$ describes the unique path in $\cof{\kappa}^{\leq n}$ of length $n$ connecting $p_0 \coloneq \emptyset$ to $p_n \coloneq p$, we assume that there is a corresponding path $X_p$ of length $n - 1$ in $G$ and such that the following conditions are met:

    \begin{enumerate}
        \item \label{item:XpPaths:i} $\cardinality{X_p\cap B_{p_i}} = 1$ for every $0\leq i < n$;
        
        \item \label{item:XpPaths:ii} $[t,b_p]\setminus X_p$ is connected;
        
        \item \label{item:XpPaths:iii} $N(w_p) \cap \upcl{b_p} \supseteq B_p$, where $w_p$ is the unique element from $X_p\cap B_{p_{n-1}}$.  
    \end{enumerate}

    Finally, if $p\in\cof{\kappa}^{\leq n}$ is a leaf node, we further assume that $B_p$ admits a family $\mathcal{B}_{p}\coloneq \set{C_q^p}:{q \in \mathrm{succ}_{\sqsubseteq}(p)}$ of a successor sequence of cliques and with $N(B \cup X_p) \cap \upcl{b_q^p} = \emptyset$ for every $q \in \successor{p}$, where $b_q^p \coloneq \min C_q^p$.
    Due to item \cref{item:ConvergingSize}, each $C_q^p$ can be chosen so that $\cardinality{C_q^p} \geq \kappa_{\ell(q)}$ for every $q\in \successor{p}$. 
    In addition, let $\mathcal{W}_p \coloneq \set{W_i^p}:{i< \cardinality{B_p}}$ be the family of disjoint subsets of $B_p$ verifying properties \cref{item:ConvergingSize}--\cref{item:DominatingBags} from the definition of successor sequence of cliques.
    This completes the list of requirements we will ensure at any given step $n$.

    To detail the inductive construction, let $q \in \mathrm{succ}_{\sqsubseteq}(p)$ be given and, relying on item \cref{item:smallSpreadingTreeClique}, fix a node $x_q \in B_p$, $x_q \neq \min{B_p}$, such that $\upcl{b_q^p}_{N(x_q)} \supseteq C_q^p$ and such that $(\Set{x_q},b_q^p)$ is a small spreading pair. 
    Then define $X_q \coloneq X_p \cup \Set{x_q}$, which we will later show to meet the requirements \cref{item:XpPaths:i}--\cref{item:XpPaths:iii}.

    If $p \neq \emptyset$, then the element $w_p$ (defined in \cref{item:XpPaths:iii}) is adjacent to $x_q$ and witnesses that $X_q$ induces a path of length $n$.
    In the following we address the premises of \cref{lemma:SendingAbove}. 
    More precisely, we will verify that $(B\cup X_q, b_q^p)$ defines a small spreading pair and that $[t, b_q^p] \setminus (B \cup X_q)$ is connected.
    Since $C^p_q \subseteq \upcl{b_q^p}_{N(B \cup X_q)}$ by the choice of $x_q$, we are even able to pass $C^p_q$ to a subinterval rooted in $b_q^p$ and of order type $\kappa_{\ell(q)}$.

    Firstly, the set $B \cup X_q$ is connected because $B$ and $X_q$ are connected and $B_{\emptyset} \subseteq N(B)$, while $X_q \cap B_{\emptyset} \neq \emptyset$ by \cref{item:XpPaths:i}. 
    Moreover, $(B \cup X_q) \cap \upcl{b_q^p} = \emptyset$ because $B \cap \upcl{t} = \emptyset$ and $X_q \subseteq \bigcup_{i=0}^{n} B_{p_i} \subseteq \odowncl{b_q^p}$ by construction and \cref{item:MinimumAboveChain}.
    Then $N(B \cup X_q) \cap \upcl{b_q^p} = N(x_q) \cap \upcl{b_q^p} \supseteq C^p_q$ given that $N(B \cup X_p) \cap \upcl{b_q^p} = \emptyset$ by the inductive assumption on $\mathcal{B}_p$.
    Since $(\Set{x_q}, b_q^p)$ is a small spreading pair by construction, it follows that $(B \cup X_q, b_q^p)$ is also a small spreading pair.
    It remains to verify that $[t, b^p_q] \setminus (B \cup X_q)$ is connected.
    By hypothesis \cref{item:XpPaths:ii} on $X_p$ it follows that $[t, b^p_q] \setminus (B \cup X_p)$ is connected (since $X_p$ does not contain $b_p$); and since $x_q$ is distinct from $\min{B_p}$ while $B_p$ induces a clique in $G$, the connectedness of $[t, b^p_q] \setminus (B \cup X_q)$ follows.
    
    All in all, the triple $(B \cup X_q, t, b_q^p)$ together with the subinterval of $C^p_q$ rooted in $b_q^p$ of order type $\kappa_{\ell(q)}$ satisfies the premises of \cref{lemma:SendingAbove} and, thus, either item \cref{item:FindingAntichain} or item \cref{item:FindingSuccessors} can be concluded.
    The following claim then implies that item \cref{item:FindingAntichain} does not hold:

    \begin{claim}\label{claim:inductiveExtendable}
        If $(B \cup X_q, b_q^p)$ is $\kappa_{\ell(q)}$--extendable from $t$ to $b_q^p$, then $(B,s)$ is $\kappa_{\ell(q)}$--extendable from $t$ to~$s$.
    \end{claim}

    \begin{proof}[Proof of~\cref{claim:inductiveExtendable}]\renewcommand{\qedsymbol}{$\blacksquare$}
        As in the definition of extendability from $t$ to $b_q^p$, let $B_q^*$ be a connected set with $B \cup X_q \subseteq B_q^* \subseteq B \cup X_q \cup \oupcl{b_q^p}$ and let $A \subseteq \oupcl{b_q^p}$ be a corresponding $(B_q^*, \kappa_{\ell(q)})$--strong antichain such that $\left( \bigcup_{a \in A} \odowncl{a} \setminus B_q^* \right) \cap \upcl{t}$ is connected.
        Since $B_q^*$ in particular contains $B$, it also witnesses that $(B,s)$ is $\kappa_{\ell(q)}$--extendable from $t$ to $s$.
    \end{proof}

    To summarize, \cref{lemma:SendingAbove} concludes item \cref{item:FindingSuccessors}: there exists an interval $B_q \subseteq N(B \cup X_q) \cap \upcl{b_q^p}$ which is a clique of order type $\kappa_{\ell(q)}\geq \lambda$ and which admits a successor sequence of cliques $\mathcal{B}_q$. 
    Again, by \cref{obs:passingToDisjointSuccessorCliques}, we may assume without loss of generality that $N(B \cup X_q) \cap \upcl{\min B} = \emptyset$ for every $B \in \mathcal{B}_q$.
    It thus remains to verify that $X_q$ meets \cref{item:XpPaths:i}--\cref{item:XpPaths:iii}.

    Firstly, item \cref{item:XpPaths:i} is satisfied by definition of $X_q$ and the hypothesis \cref{item:XpPaths:i} on $X_p$.
    Towards item \cref{item:XpPaths:ii}, recall that $[t, b^p_q] \setminus (B \cup X_q) = [t, b^p_q] \setminus X_q$ is connected and, since $X_q$ is disjoint from $[b^p_q, \min{B_q}]$, also $[t, \min{B_q}] \setminus X_q$ is connected.
    Item \cref{item:XpPaths:iii} then follows from the inductive hypothesis that $N(B \cup X_p) \cap \upcl{b^p_q} = \emptyset$ and $B_q \subseteq N(B \cup X_q) \cap \upcl{b^p_q}$.

    After applying this procedure to all leaves $p$ of $\cof{\kappa}^{\leq n}$, we obtain a minor embedding of $\cof{\kappa}^{\leq n+1}$ into $\upcl{s}$ which meets the inductive assumptions.

    At the end of this recursive process and since $\cof{\kappa}^{<\omega} = \bigcup_{n \in \mathbb{N}}\cof{\kappa}^{\leq n}$, the branch sets $\bigcup_{p \in \cof{\kappa}^{<\omega}} \bigcup \mathcal{W}_p$ define a minor model of $\mathbb{T}_{\kappa}$, thus addressing item \cref{item:OutputTreeOfCliques}. 
    After all, for every $p \in \cof{\kappa}^{<\omega}$ and by property \cref{item:DominatingBags} from the definition of successor sequence of cliques, we can regard each element from $\mathcal{W}_p$ as a bag of a $K_{\kappa_{\ell(p)}}$ minor and which contains a vertex that, for every $q \in \mathrm{succ}_{\sqsubseteq}(p)$, neighbors every vertex in the $\kappa_{\ell(q)}$--sized clique $B_q \subseteq C_q^p$.
    In addition, $\sup \set{\cardinality{B_q}}:{q\in \mathrm{succ}_{\sqsubseteq}(p)} = \sup \set{\kappa_{\ell(q)}}:{q \in \mathrm{succ}_{\sqsubseteq}(p)} = \sup \set{\kappa_{\alpha}}:{\alpha > \ell(p)} = \kappa$ by construction.
\end{proof}

Now we can prove \cref{lem:NoCliqueTree}, which we restate for convenience.

\mainSingularLemma*

\begin{proof}
    Suppose that $G$ contains no minor embedding of $\mathbb{T}_{\kappa}$.
    Due to the short neighborhood property of $G$, the order tree $\Delta \coloneq N(B) \cap \upcl{t}$ verifies item \cref{prop:shortNeighborhoodProp:case1} or item \cref{prop:shortNeighborhoodProp:case2}.
    The first item \cref{prop:shortNeighborhoodProp:case1} claims that there is an antichain $A \subseteq \Delta$ of size $\kappa$ with $\cardinality{\upcl{a}_\Delta} < \kappa$ for every $a \in A$.   
    In particular, $A$ is a $(B,\kappa)$--strong antichain by definition and the set $B_* \coloneq B$ witnesses that $(B,t)$ is extendable.

    Thus, it remains to consider the case where $\Delta$ satisfies item \cref{prop:shortNeighborhoodProp:case2}, that is, $\Delta$ admits an antichain $\Set{a_\alpha}_{\alpha < \cof{\kappa}} \subseteq \Delta$ such that $\upcl{a_{\alpha}}_{\Delta}$ has cardinality $< \kappa$ and height at least $\kappa^+_{\alpha}$.
    By \cref{lemma:TreeOfCliquesConstruction}, it follows that $(B, a_\alpha)$ is $\mu_{\alpha}$--extendable from $t$ to $a_{\alpha}$ for some $\kappa_\alpha \leq \mu_{\alpha} < \kappa$ for every $\alpha < \cof{\kappa}$.
    Hence, let $B_{\alpha}$ be a connected set with $B \subseteq B_{\alpha} \subseteq B \cup \oupcl{a_{\alpha}}$ and let $A_{\alpha} \subseteq \oupcl{a_{\alpha}}$ be a corresponding $(B_{\alpha},\mu_{\alpha})$--strong antichain such that $D_{\alpha} \coloneq \left(\bigcup_{a\in A_{\alpha}}\odowncl{a}\setminus B_{\alpha}\right)\cap \upcl{t}$ is connected. 
    In particular, as a union of connected sets containing $B$, the set $B_* \coloneq \bigcup_{\alpha < \cof{\kappa}}B_{\alpha}$ is also connected and contains $B$. 
    Similarly, $A \coloneq \bigcup_{\alpha < \cof{\kappa}}A_{\alpha}$ is an antichain since both $A_{\alpha} \subseteq \upcl{a_{\alpha}}$ and $\set{a_{\alpha}}:{\alpha < \cof{\kappa}}$ are antichains.
    Moreover, by normality of $T$, the definition of $B_*$ implies that 
    \[
        N(B_*)\cap \upcl{a} = N(B_*\cap \downcl{a})\cap \upcl{a} = N(B\cup (B_{*}\cap \upcl{a_{\alpha}}))\cap \upcl{a} = N(B_{\alpha})\cap \upcl{a}
    \]
    for every $a\in A$ and the unique corresponding index $\alpha < \cof{\kappa}$ such that $a\in A_{\alpha}$. 
    Hence, $A$ is in particular a $(B_*,\kappa)$--strong antichain, because $A_{\alpha}$ is a $(B_{\alpha},\mu_{\alpha})$--strong antichain for every $\alpha < \cof{\kappa}$ and $\cardinality{A} = \sum_{\alpha < \cof{\kappa}} \mu_{\alpha} = \kappa$.
    Again relying on the definition of $B_*$, we can write $D \coloneq \left( \bigcup_{a\in A} \odowncl{a} \setminus B_* \right) \cap \upcl{t}$ as 
    
    $$D = \bigcup_{\alpha < \cof{\kappa}} \left[ \bigcup_{a\in A_{\alpha}}(\odowncl{a} \setminus B_*) \cap \upcl{t} \right]  =  \bigcup_{\alpha < \cof{\kappa}} \left[ \bigcup_{a\in A_{\alpha}}(\odowncl{a}\setminus B_{\alpha}) \cap \upcl{t}\right] = \bigcup_{\alpha < \cof{\kappa}}D_{\alpha}.$$

    Then, as a union of connected sets containing $t$, the set $D$ is also connected. 
    All in all, the connected set $B_*$ and the antichain $A$ verify that $(B,t)$ is extendable. 
\end{proof}

\bibliographystyle{amsplain}
\bibliography{collective.bib}

@article{enddefiningsequences,
  title={Characterising {$k$}-connected sets in infinite graphs},
  author={Gollin, J. Pascal and Heuer, Karl},
  journal={Journal of Combinatorial Theory, Series B},
  volume={157},
  pages={451--499},
  year={2022},
  publisher={Elsevier}
}

@book{bibel,
    author = {Diestel, Reinhard},
	edition = {6th},
	publisher = {Springer},
	title = {{Graph Theory}},
	year = {2025},
    series = {Graduate Texts in Mathematics}, 
    volume = {173}
}

@book{pitzbook,
  author    = {Pitz, Max},
  title     = {Infinite Graphs},
  series    = {Graduate Texts in Mathematics},
  publisher = {Springer},
  note      = {In preparation},
}

@article{knappe2024immersion,
  title={The immersion-minimal infinitely edge-connected graph},
  author={Knappe, Paul and Kurkofka, Jan},
  journal={Journal of Combinatorial Theory, Series B},
  volume={164},
  pages={492--516},
  year={2024},
  publisher={Elsevier}
}

@article{ErdosTarski1943,
  author  = {Erd{\H o}s, P. and Tarski, A.},
  title   = {On families of mutually exclusive sets},
  journal = {Annals of Mathematics},
  series  = {2},
  volume  = {44},
  number  = {2},
  pages   = {315--329},
  year    = {1943},
  doi     = {10.2307/1968767}
}

@article{seymour1993binarytree,
  title={Excluding infinite trees},
  author={Seymour, Paul D. and Thomas, Robin},
  journal={Transactions of the American Mathematical Society},
  volume={335},
  number={2},
  pages={597--630},
  year={1993}
}

@article{OrderCompatible,
  title={On order compatible paths in infinite graphs},
  author={Pitz, Max and Real, Lucas and Schaut, Roman},
  journal={arXiv preprint arXiv:2603.08454},
  year={2026}
}

@article{halin78,
    title={Simplicial decompositions of infinite graphs},
    author={Rudolf Halin},
    journal={Annals of Discrete Mathematics}, 
    volume={3}, 
    pages={93--109},
    year={1978}
}

@article {BrochetDiestel,
    AUTHOR = {Brochet, J.-M. and Diestel, R.},
     TITLE = {Normal tree orders for infinite graphs},
   JOURNAL = {Trans. Amer. Math. Soc.},
  FJOURNAL = {Transactions of the American Mathematical Society},
    VOLUME = {345},
      YEAR = {1994},
    NUMBER = {2},
     PAGES = {871--895},
      ISSN = {0002-9947,1088-6850},
   MRCLASS = {05C05},
MRREVIEWER = {H.\ A.\ Jung},
       DOI = {10.2307/2155004},
       URL = {https://doi.org/10.2307/2155004},
}

@article{RSTExcludingSubdivisions,
    title={Excluding subdivisions of infinite cliques},
    author={Robertson, N. and Seymour, P.D. and Thomas, R.},
    journal={Transactions of the American Mathematical Society},
    volume={332}, 
    number={1},
    pages={211--223},
    year={1992}
}

@book {KennethKunen,
    AUTHOR = {Kunen, Kenneth},
     TITLE = {Set theory},
    SERIES = {Studies in Logic (London)},
    VOLUME = {34},
 PUBLISHER = {College Publications, London},
      YEAR = {2011},
     PAGES = {viii+401},
      ISBN = {978-1-84890-050-9},
   MRCLASS = {03-01 (03E05 03E10 03E35 03E40 03E45 03E50)},
MRREVIEWER = {Klaas\ Pieter\ Hart},
}

@article {Oporowski,
    AUTHOR = {Oporowski, Bogdan and Oxley, James and Thomas, Robin},
     TITLE = {Typical subgraphs of {$3$}- and {$4$}-connected graphs},
   JOURNAL = {J. Combin. Theory Ser. B},
  FJOURNAL = {Journal of Combinatorial Theory. Series B},
    VOLUME = {57},
      YEAR = {1993},
    NUMBER = {2},
     PAGES = {239--257},
      ISSN = {0095-8956,1096-0902},
   MRCLASS = {05C40 (05C75)},
MRREVIEWER = {Kevin\ Lawrence\ McAvaney},
       DOI = {10.1006/jctb.1993.1019},
       URL = {https://doi.org/10.1006/jctb.1993.1019},
}

@article {halin1967,
    AUTHOR = {Halin, R.},
     TITLE = {Unterteilungen vollst\"andiger {G}raphen in {G}raphen mit
              unendlicher chromatischer {Z}ahl},
   JOURNAL = {Abh. Math. Sem. Univ. Hamburg},
  FJOURNAL = {Abhandlungen aus dem Mathematischen Seminar der Universit\"at
              Hamburg},
    VOLUME = {31},
      YEAR = {1967},
     PAGES = {156--165},
      ISSN = {0025-5858,1865-8784},
   MRCLASS = {05.55},
MRREVIEWER = {W.\ Imrich},
       DOI = {10.1007/BF02992395},
       URL = {https://doi.org/10.1007/BF02992395},
}

@article{lorenz2025obstructions,
  title   = {Obstructions for normally spanned sets of vertices},
  author  = {Lorenz, Nicola and Pitz, Max},
  journal = {arXiv preprint arXiv:2510.04367},
  year    = {2025},
  url     = {https://arxiv.org/abs/2510.04367}
}

@article{BuergerKurkofka2022,
author = {B\"urger, Carl and Kurkofka, Jan},
title = {End-faithful spanning trees in graphs without normal spanning trees},
journal = {Journal of Graph Theory},
volume = {101},
number = {1},
pages = {95-105},
doi = {https://doi.org/10.1002/jgt.22813},
url = {https://onlinelibrary.wiley.com/doi/abs/10.1002/jgt.22813},
year = {2022}
}

@article {kurkofka2022,
    AUTHOR = {Kurkofka, Jan},
     TITLE = {Every infinitely edge-connected graph contains the {F}arey
              graph or {$T_{\aleph _0}*t$} as a minor},
   JOURNAL = {Math. Ann.},
  FJOURNAL = {Mathematische Annalen},
    VOLUME = {382},
      YEAR = {2022},
    NUMBER = {3-4},
     PAGES = {1881--1900},
      ISSN = {0025-5831,1432-1807},
   MRCLASS = {05C63 (05C10 05C40 05C55 05C83)},
MRREVIEWER = {Tsz\ Lung\ Chan},
       DOI = {10.1007/s00208-021-02259-7},
       URL = {https://doi.org/10.1007/s00208-021-02259-7},
}

@article{HalinJung1963,
    AUTHOR = {Halin, R. and Jung, H. A.},
     TITLE = {{\"U}ber {M}inimalstrukturen von {G}raphen, insbesondere von $n$-fach zusammenh{\"a}ngenden {G}raphen},
   JOURNAL = {Math. Ann.},
    VOLUME = {152},
      YEAR = {1963},
     PAGES = {75--94},
}

@article{Dirac1977,
    AUTHOR = {Dirac, G. A.},
     TITLE = {Cardinal-determining subgraphs of infinite graphs},
   JOURNAL = {Contributi del Centro Linceo Interdisciplinare di Scienze Matematiche e loro Applicazioni},
    VOLUME = {34},
      YEAR = {1977},
     PAGES = {3--14},
}

\end{document}